\documentclass{article}

\usepackage{amsmath,amsfonts,amssymb,amsthm}
\usepackage{slashbox}
\usepackage{graphics}
\newenvironment{tabud}[3]{\setbox2=\hbox{\bf #2}%
\renewcommand{\arraystretch}{#3}\begin{tabular}{#1}%
}{\end{tabular}\renewcommand{\arraystretch}{1}\nopagebreak%
\vskip6pt{\bf Table \thesubsection.\box2}}
\newcommand{\IF}{\text{ if }}
\newcommand{\maparrow}[1]{\underset{#1}{\stackrel{\sim}{\longrightarrow}}}
\newcommand{\ind}{\operatorname{index}}
\newcommand{\twomatrix}[4]{\left(
\begin{smallmatrix} #1 & #2\\ #3&#4\end{smallmatrix}\right)}

\makeatletter
\@addtoreset{equation}{subsection}
\def\@makefnmark{\empty}
\def\@seccntformat#1{\csname the#1\endcsname%
\expandafter\ifx\csname#1\endcsname\subsection.\fi\quad}
\makeatother

\theoremstyle{plain}
\newtheorem{thm}{Theorem}[subsection]
\newtheorem{cor}[thm]{Corollary}
\newtheorem{prop}[thm]{Proposition}
\newtheorem{conj}[thm]{Conjecture}
\newtheorem{la}[thm]{Lemma}
\newtheorem{fact}[thm]{Fact}

\theoremstyle{definition}
\newtheorem{dfn}[thm]{Definition}

\newtheorem{ex}[thm]{Example}

\theoremstyle{remark}
\newtheorem{rem}[thm]{Remark}

\def\tr#1{\mathord{\mathopen{{\vphantom{#1}}^t}#1}}
\def\st#1{{#1}^*}

\def\ritem[#1]{\item[\rm #1]}

\newcommand{\N}{\mathbb{N}}
\newcommand{\Z}{\mathbb{Z}}

\newcommand{\R}{\mathbb{R}}
\newcommand{\C}{\mathbb{C}}
\newcommand{\F}{\mathbb{F}}
\newcommand{\bbH}{\mathbb{H}}

\newcommand{\g}{\mathfrak{g}}
\newcommand{\h}{\mathfrak{h}}
\newcommand{\p}{\mathfrak{p}}
\newcommand{\fa}{\mathfrak{a}}
\newcommand{\fl}{\mathfrak{l}}
\newcommand{\fk}{\mathfrak{k}}
\newcommand{\fo}{\mathfrak{o}}
\newcommand{\ft}{\mathfrak{t}}

\newcommand{\spin}{\mathfrak{spin}}

\newcommand{\z}{\{0\}}
\newcommand{\set}[2]{\{ #1 : #2\}}
\newcommand{\floor}[1]{\lfloor #1\rfloor}

\newcommand{\id}{\operatorname{id}}

\newcommand{\Ad}{{\operatorname{Ad}}}

\newcommand{\Image}{\operatorname{Image}}
\newcommand{\Ker}{\operatorname{Ker}}
\newcommand{\Span}{\operatorname{Span}}

\newcommand{\End}{\operatorname{End}}
\newcommand{\Hom}{\operatorname{Hom}}
\newcommand{\Aut}{\operatorname{Aut}}
\newcommand{\rank}{\operatorname{rank}}
\newcommand{\diag}{\operatorname{diag}}

\newcommand{\rSpan}{\Span_\R}
\newcommand{\rrank}{\R\text{-}\rank}

\newcommand{\tilf}{\tilde{f}}
\newcommand{\simarrow}{\stackrel{\sim}{\to}}

\newcommand{\proper}{\pitchfork}
\newcommand{\similar}{\sim}
\newcommand{\ssimilar}{{\,\similar_s\,}}

\newcommand{\properin}[3]{#1\proper #2\text{ in $#3$}}
\newcommand{\similarin}[3]{#1\similar #2\text{ in $#3$}}
\newcommand{\ssimilarin}[3]{#1\ssimilar #2\text{ in $#3$}}

\newcommand{\codd}{C_{\operatorname{odd}}}
\newcommand{\ceven}{C_{\operatorname{even}}}
\newenvironment{tabu}[1]{\renewcommand{\arraystretch}{1.3}
\begin{tabular}{#1}}{\end{tabular}\renewcommand{\arraystretch}{1}%
\vskip6pt{\bf Table \thesubsection}}
\newenvironment{tabu2}[2]{\renewcommand{\arraystretch}{#2}
\begin{tabular}{#1}}{\end{tabular}\renewcommand{\arraystretch}{1}%
\vskip6pt{\bf Table \thesubsection}}

\title{Compact Clifford-Klein forms\\
of symmetric spaces -- revisited\\
{\normalsize \bf In memory of Professor Armand Borel}
}
\author{Toshiyuki Kobayashi and Taro Yoshino}
\date{}

\begin{document}
\maketitle
\begin{abstract}
This article discusses the existence problem of
a compact quotient of a symmetric space by
a properly discontinuous
group with emphasis on the non-Riemannian case.
Discontinuous groups are not always abundant in a homogeneous
space $G/H$ if $H$ is non-compact. The
first half of the article elucidates general machinery to study
discontinuous groups for $G/H$, followed by the most update and
complete list of
symmetric spaces with/without compact quotients. In the second
half, as applications of general theory, we prove:
(i) there exists
a $15$ dimensional compact pseudo-Riemannian manifold of
signature $(7,8)$ with constant curvature, (ii) there
exists a compact quotient of the complex sphere of dimension
$1$, $3$ and $7$, and (iii) there exists a compact quotient
of the tangential space form of signature $(p,q)$ if and only
if $p$ is smaller than the Hurwitz-Radon number of $q$.
\end{abstract}
\tableofcontents
\section{Introduction}
\setcounter{subsection}{1}
A Clifford-Klein form of a symmetric space $M=G/H$ is the
double coset manifold $\Gamma\backslash G/H$, where
$\Gamma\subset G$ is a discontinuous group for $M$.
Geometrically, we are dealing with
a complete, locally symmetric space.
The aim of this article is to discuss recent progress on the
following problem
after Borel \cite{borel} for the Riemannian case and
Kobayashi \cite{koba88} for the non-Riemannian case.
\vskip10pt
\noindent
{\bf Problem A.} Does a symmetric space $M$ admit
compact Clifford-Klein forms?
\vskip10pt

In his celebrated paper
``Compact Clifford-Klein forms
of symmetric spaces'', Borel proved Problem A affirmatively,
but under the tacit assumption that $M$ is Riemannian.

\begin{thm}[{\cite{borel}}]\label{thm:borel}
Compact Clifford-Klein forms always exist if $M$ is Riemannian.
\end{thm}
Opposite extremal cases may happen for
general symmetric spaces $M$.
\begin{thm}[see Theorem~\ref{thm:paraH}]
A compact Clifford-Klein form never exists if $M$ is para-Hermitian.
\end{thm}
Among symmetric spaces, irreducible ones were classified infinitesimally by
Berger \cite{berger}. Even for irreducible symmetric spaces, Problem A
is far from being solved.

In this paper, we shall discuss general machinery
to study Problem A,
and then apply it to specific
symmetric spaces. Then we shall
provide a table of (irreducible) symmetric spaces
that are proved to admit/not to admit compact Clifford-Klein forms.
It includes:
\begin{thm}[see Section~\ref{subsec:complex}]\label{thm:intro:complex_sphere}
The $7$ dimensional complex sphere admits a compact Clifford-Klein form.
\end{thm}

\begin{thm}[\cite{calabi},\cite{kulkarni},\cite{klingler}; see
Fact~\ref{fact:lorentz2}]
\label{thm:lorentz}
 There exists an $n$-dimensional compact
Lorentz manifold 
of constant curvature $\kappa$ if and only if
$\kappa=0$, or $\kappa<0$ and $n$ is odd.
\end{thm}
\begin{thm}[see Section~\ref{sec:spherical}]\label{thm:spherical78}
There exists a $15$ dimensional compact complete pseudo-Riemannian manifold
 of signature $(7,8)$ with constant positive sectional
curvature.
\end{thm}

These latter two theorems are regarded as a partial answer to
the {\it space form problem}, which asks whether or not there exists
a compact pseudo-Riemannian manifold of general signature $(p,q)$ with
constant sectional curvature (see Conjecture~\ref{conj:spherical}).
Furthermore, its `tangential' version for the Cartan motion group
is formulated and completely solved. We recall from \cite{hurwitz} or 
\cite{radon} that the
Hurwitz-Radon number $\rho(q)$ is defined to be $8\alpha+2^\beta$ if we write
$q=u\cdot 2^{4\alpha+\beta}$ $(u$ is odd,
$0\leq\beta\leq 3)$.

\begin{thm}[see Section~\ref{sec:tangential}]\label{thm:tangential}
The tangential symmetric space of \hbox{$O(p+1,q)/O(p,q)$}
admits a compact Clifford-Klein form if and only if
$p<\rho(q)$.
\end{thm}
This article focuses on Problem A for symmetric spaces.
For non-symmetric spaces such as $SL(n,\R)/SL(m,\R)\ (n>m)$,
in addition to the methods that we shall explain
in Section~\ref{sec:methods}, various other approaches have been also
made including ergodic theory, symplectic
geometry and unitary representations 
(\cite[96b]{koba92}, \cite{benoistlabourie}, \cite{corlette},
\cite{zimmer}, \cite{lmz}, \cite{lz},
\cite{margulis97}, \cite{oh}, \cite{shalom}, \cite{ohwitte02}),
though
some of these approaches are not applicable to symmetric spaces.
See also recent survey papers \cite{margulis00}, \cite[02]{koba01},
\cite{iozzi}  and references therein.

This article has a survey nature largely in
the first half, and gives new results on Problem A in the latter half.
It is organized as follows:
Section 2 exhibits the most update knowledge on the
existence problem of compact Clifford-Klein forms of symmetric spaces.
Various examples given in Section~\ref{sec:example}
 are obtained by 
a general theory to Problem A based on
group theoretic study of discontinuous groups and
geometric consideration of Clifford-Klein forms.
An introduction to the general theory
will be the main topic of Section 3.
In Section 4, we apply the general theory to
the space form problem, and in particular, give a proof of
Theorem~\ref{thm:spin_ck_form} (this includes
Theorems \ref{thm:intro:complex_sphere}
and~\ref{thm:spherical78}) using
Clifford algebras
associated to indefinite quadratic forms.
Section 5 studies Problem A for the tangential symmetric spaces,
and in particular, gives a proof of
Theorem~\ref{thm:tangential}.

Professor Armand Borel invited
the first author to the Institute 
for Advanced Study at Princeton
with warmest encouragement on
 the initial stage of this work,
right after his visit to Kyoto in 1990, 
 where the first author gave a presentation
 at the RIMS workshop organized by Professor Satake.
It is an honor to dedicate this paper
to the memory of Professor Borel.

\section{Examples}\label{sec:example}
\subsection{Notation and Definition}
First, suppose $H$ is a closed subgroup of a Lie group $G$.
Assume that $\Gamma$ is a discrete subgroup of $G$
acting properly discontinuously and freely
on the homogeneous space $G/H$.
Then, the double coset space $\Gamma\backslash G/H$ is
Hausdorff in the quotient topology and the quotient map
$$\pi: G/H \rightarrow \Gamma\backslash G/H$$
becomes a covering map, so that $\Gamma\backslash G/H$ carries
naturally a manifold structure.
\begin{dfn}\label{dfn:clifford}
The resulting manifold $\Gamma\backslash G/H$ is
said to be a {\it Clifford-Klein form} of $G/H$, 
and the discrete subgroup $\Gamma$ is a {\it discontinuous group} for $G/H$.
If the double coset space $\Gamma\backslash G/H$ is furthermore compact,
we say $\Gamma$ is a {\it cocompact discontinuous group} for $G/H$,
and we say $G/H$ has a {\it tessellation}. A cocompact discontinuous
group is also referred to as a {\it uniform lattice} for $G/H$
or a {\it crystallographic group} for $G/H$.
\end{dfn}
Second, suppose $\sigma$ is an automorphism of a Lie group
$G$ such that
$\sigma^2=\id$. Then, the set of fixed points
$G^\sigma:=\set{g\in G}{\sigma(g)=g}$ is
a closed subgroup of $G$. We write ${(G^\sigma)}_0$ for the identity
component of $G^\sigma$. The homogeneous space $G/H$ is called
a {\it symmetric space} if 
${(G^\sigma)}_0 \subset H \subset G^\sigma$. A symmetric space is
said to be {\it semisimple} [{\it reductive,\dots}] if $G$ is a semisimple 
[reductive,\dots] Lie group. 

Any Clifford-Klein form of a symmetric space (therefore,
a symmetric space itself) becomes a complete locally symmetric space
through its canonical affine connection.
Conversely, any complete locally symmetric space is
represented as a Clifford-Klein form $\Gamma\backslash G/H$
of some symmetric space $G/H$ (\cite{kn69}).

\subsection{Symmetric spaces with compact Clifford-Klein forms}
\label{subsec:admit_cck}

\newcounter{symcnt}
\newcommand{\symm}[1]{ $\thesymcnt$ \stepcounter{symcnt}  & $#1$\\ }
\setcounter{symcnt}{1}
\newcommand{\Hline}{\hline}

Any Riemannian symmetric space admits compact Clifford-Klein forms,
as we have seen in Theorem~\ref{thm:borel}.
This result was proved by finding a cocompact discrete subgroup $\Gamma$
of a linear semisimple Lie group $G'$. In turn,
the group manifold $G/H:=G'\times G'/\diag G'$ also
admits a cocompact discontinuous group $\Gamma\times \{e\}$.
Since $G'\times G'/\diag G'$ is a
symmetric space by the involution $\sigma(x,y):=(y,x)$,
this gives another example of symmetric spaces that have
compact Clifford-Klein forms.

Apart from these Riemannian symmetric spaces and
group manifold cases, there are some more
semisimple symmetric spaces that have compact Clifford-Klein
forms.
\begin{thm}[see Corollary~\ref{cor:admit_cck}]\label{thm:admit_cck}
The following semisimple symmetric spaces $G/H$ admit compact Clifford-Klein
forms ($n=1,2,\dots$).
\begin{center}
\begin{tabu}{c|c}
& $G/H$ \\ \Hline
\symm{SU(2,2n)/Sp(1,n)}
\symm{SU(2,2n)/U(1,2n)}
\symm{SO(2,2n)/U(1,n)}
\symm{SO(2,2n)/SO(1,2n)}
\symm{SO(4,4n)/SO(3,4n)}
\symm{SO(4,4)/SO(4,1)\times SO(3)}
\symm{SO(4,3)/SO(4,1)\times SO(2)}
\symm{SO(8,8)/SO(7,8)}
\symm{SO(8,\C)/SO(7,\C)}
\symm{SO(8,\C)/SO(7,1)}
\symm{SO^*(8)/U(3,1)}
\symm{SO^*(8)/SO^*(6)\times SO^*(2)}
\end{tabu}
\end{center}
\end{thm}

Clearly, the same conclusion still holds if we replace $G$ and $H$ by
locally isomorphic groups up to finite coverings and finitely
many connected components.
For example, the cases 11, 12 and $SO(2,6)/U(1,3)$ are
locally isomorphic to each other.
\subsection{Para-Hermitian symmetric spaces}\label{subsec:paraH}

We say a manifold $M^{2n}$ of even dimension has a {\it paracomplex
structure} if the tangent bundle $TM$ splits into a Whitney direct sum
$T^+M\oplus T^-M$ with equi-dimensional fibers and if $T^\pm M$ are
completely integrable. A pseudo-Riemannian metric $g$ over a paracomplex
manifold $M$ is said to be a {\it para-Hermitian metric} if $T_x^\pm M
\subset T_xM$ are maximal totally isotropic subspaces with respect to
$g$ for every point $x$ in $M$. 
A semisimple symmetric space equipped with
$G$-invariant para-Hermitian structure is said to be a {\it para-Hermitian
symmetric space}.
Kaneyuki and Kozai \cite{kaneyuki} characterized para-Hermitian
symmetric spaces by the property that
the center $C$ of $H$ is non-compact,
and gave a classification of
irreducible ones
based on the latter property.

Concerning Problem A, a para-Hermitian
symmetric space is an opposite extreme to a Riemannian
symmetric space:
\begin{thm}\label{thm:paraH}
 None of para-Hermitian symmetric spaces
admits a compact
Clifford-Klein form.
\end{thm}
\begin{proof}
If $G/H$ is a para-Hermitian symmetric space, then $H$ has
the non-compact center $C$ and its Lie algebra $\h$ coincides with
the centralizer of $C$ in $\g$. Since any maximal split abelian 
subspace $\fa$ of $\g$ contains $C$, $\fa$ is contained in
the centralizer of $C$ in $\g$, namely $\h$.
Therefore $\rrank G=\rrank H$.
Now, apply Theorem~\ref{thm:cr:calabi} (the Calabi-Markus
phenomenon).
\end{proof}
\setcounter{symcnt}{1}
\begin{ex} Here are some examples of para-Hermitian symmetric spaces $G/H$.
In particular, there are no compact Clifford-Klein forms of these $G/H$:
\begin{center}
\begin{tabu}{c|c}
& $G/H$ \\ \Hline
\symm{SL(p+q,\R)/S(GL(p,\R)\times GL(q,\R))}
\symm{Sp(n,\R)/GL(n,\R)}
\symm{SO^*(4n)/U^*(2n)}
\symm{SU^*(2(p+q))/S(U^*(2p)\times U^*(2q))}
\symm{Sp(n,n)/U^*(2n)}
\symm{SU(n,n)/SL(n,\C)\times\R}
\symm{SO(n,n)/GL(n,\R)}
\symm{SO(p,q)/SO(p-1,q-1)\times SO(1,1)}
\end{tabu}
\end{center}
\end{ex}

\begin{rem}
Any para-Hermitian symmetric space is realized as a 
hyperbolic orbit of the adjoint action
on the Lie algebra $\g$ (see \cite{koba98h}). Hence, Theorem~\ref{thm:paraH}
can be regarded as a special case of the following theorem
(see \cite[Theorem 1.3]{koba92n}):
{\it a semisimple orbit with a compact Clifford-Klein form must 
be an elliptic orbit and have
a $G$-invariant complex structure.}
Two different proofs are known for this theorem;
the original proof was based on the criterion of properness
(see Theorem~\ref{thm:non-cpt-dim}), while Benoist and
Labourie \cite{benoistlabourie} gave an alternative proof
based on symplectic geometry.
\end{rem}
\subsection{Complex symmetric spaces $G_\C/K_\C$}
\label{subsec:complex}
A symmetric space $G/H$ is a {\it complex symmetric space} if
$G$ is a complex Lie group and $\sigma$ is an involutive holomorphic
automorphism of $G$ such that $(G^\sigma)_0\subset H\subset G^\sigma$.

One of typical examples of complex symmetric spaces are
$SO(n+1,\C)/SO(n,\C)$, which is biholomorphic to the $n$-dimensional
complex sphere
$$S^n_{\C}:=\set{(z_1,\dots,z_{n+1})\in\C^{n+1}}{ z_1^2+\dots+z_{n+1}^2=1}.$$
\begin{thm}\label{thm:complex_sphere}
 The $n$-dimensional complex sphere $SO(n+1,\C)/SO(n,\C)$
admits a compact Clifford-Klein form if $n=1,3$ or~$7$.
\end{thm}
\begin{proof}
Clear in the case $n=1$ because $SO(2,\C)/SO(1,\C)\simeq \C^*/\{1\}$.
For $n=3$,
$SO(4,\C)/SO(3,\C)$ is
locally isomorphic to the group manifold $SO(3,\C)\times SO(3,\C)/
SO(3,\C)$ and therefore admits a compact Clifford-Klein form as we
saw in Section~\ref{subsec:admit_cck}. The proof for
the case $n=7$ will be postponed until Section~\ref{sec:spherical}
(see Theorem~\ref{thm:spin_ck_form}).
\end{proof}
We pose:
\begin{conj}\label{conj:complex_sphere}
 The $n$-dimensional complex sphere $SO(n+1,\C)/SO(n,\C)$
admits a compact Clifford-Klein form if and only if $n=1,3$ or~$7$.
\end{conj}
\begin{rem}
The condition $n=1,3$ or~$7$ appears as a necessary and sufficient
condition that there exist $n$ vector fields on the sphere
$S^n$ which are linearly independent at every point,
as was proved by Bott-Milnor, Kervaire, and Atiyah-Hirzebruch
(\cite{bottmilnor}, \cite{kervaire}, \cite{atiyahhirz}).
This is not just a coincidence. In fact, we can prove that the
tangential symmetric space (see Definition~\ref{dfn:tangential})
of $S^n_\C$ admits a compact Clifford-Klein form if and only if
$n=1,3$ or $7$ based on this result.
In Section~\ref{sec:tangential}, we shall see the relation between
Problem A and vector fields on the sphere in another setting arisen
from space form problem.
\end{rem}
The complex sphere $SO(n+1,\C)/SO(n,\C)$ is an {\it irreducible}
complex symmetric space if $n=2$ or $n\geq 4$.
The local classification of
irreducible complex symmetric spaces
is equivalent to the classification of real simple Lie
algebras, and there are 10 classical series (see Table~\thesubsection) and 22 exceptional ones
(\'{E}. Cartan 1914).
Locally, irreducible complex symmetric spaces
 are obtained as `complexification' of
irreducible Riemannian symmetric spaces.
Hereafter,
we shall write $G_\C/K_\C$ for the complex semisimple symmetric
space if it is
a complexification of the Riemannian symmetric space
$G/K$. For the reader's convenience, we present a list of
classical cases:
\newcommand{\symml}[2]{\symm{#1$ & $#2}}
\setcounter{symcnt}{1}
\begin{center}
\begin{tabu}{r|c|c}
&$G/K$ & $G_\C/K_\C$\\ \Hline
\symml{SL(n,\R)/SO(n)}{SL(n,\C)/SO(n,\C)}
\symml{SU^*(2n)/Sp(n)}{SL(2n,\C)/Sp(n,\C)}
\symml{SU(p,q)/S(U(p)\times U(q))}
      {SL(p+q,\C)/S(GL(p,\C)\times GL(q,\C))}
\symml{SO_0(p,q)/SO(p)\times SO(q)}
      {SO(p+q,\C)/SO(p,\C)\times SO(q,\C)}
\symml{SO^*(2n)/U(n)}{SO(2n,\C)/GL(n,\C)}
\symml{Sp(n,\R)/U(n)} {Sp(n,\C)/GL(n,\C)}
\symml{Sp(p,q)/Sp(p)\times Sp(q)}
      {Sp(p+q,\C)/Sp(p,\C)\times Sp(q,\C)}
\symml{SL(n,\C)/SU(n)}{SL(n,\C)\times SL(n,\C)/SL(n,\C)}
\symml{SO(n,\C)/SO(n)}{SO(n,\C)\times SO(n,\C)/SO(n,\C)}
\symml{Sp(n,\C)/Sp(n)}{Sp(n,\C)\times Sp(n,\C)/Sp(n,\C)}
\end{tabu}
\end{center}
We recall from Borel's result (Theorem~\ref{thm:borel}) that
any Riemannian symmetric space $G/K$ admits a compact Clifford-Klein
form (see the left-hand side). In contrast,
the following seems plausible:
\begin{conj}\label{conj:complex}
An irreducible complex symmetric space $G_\C/K_\C$
admits a compact Clifford-Klein form
if and only if 
$G_\C/K_\C$ is locally isomorphic to 
either $SO(8,\C)/SO(7,\C)$ or a group manifold.
\end{conj}
It should be noted:
\begin{prop}
Conjecture~\ref{conj:complex_sphere} is equivalent to
Conjecture~\ref{conj:complex}.
\end{prop}
\begin{proof}
What remains to prove is that an irreducible complex symmetric space
$G_\C/K_\C$ does not admit a compact Clifford-Klein form if $G_\C/K_\C$
is not locally isomorphic to $SO(n+1,\C)/SO(n,\C)$ or a group manifold.
By using
Theorem~\ref{thm:non-cpt-dim}, we proved this statement holds in
\cite[Example 1.9]{koba92}
except for the three cases $G_\C/K_\C\approx SO(2m+2,\C)/SO(2m+1,\C),
SL(2m,\C)/Sp(m,\C)$
and $E_{6,\C}/F_{4,\C}$ (see also Corollary~\ref{cor:not_exist}).
Furthermore, Benoist \cite{benoist} proved that
$SO(4m+2,\C)/SO(4m+1,\C)$ and the latter two cases
do not have compact Clifford-Klein forms
(see Corollary~\ref{cor:benoist}). Hence, Proposition is proved.
\end{proof}
As in the above proof, Conjecture~\ref{conj:complex_sphere} is
now reduced to showing the non-existence of
a compact Clifford-Klein form for 
$SO(4m,\C)/SO(4m-1,\C)\quad (m\geq 3)$.
We remark that $SO(4,\C)/SO(3,\C)$ and $SO(8,\C)/SO(7,\C)$
admit compact Clifford-Klein forms by Theorem~\ref{thm:complex_sphere}.

\subsection{Space form problem}\label{subsec:space_form}
A {\it space form} is a pseudo-Riemannian manifold $M$ with constant sectional
curvature.
This is another example of locally symmetric space.
It is an unsolved problem for which signature of indefinite
metric there exists a compact space form. We give the following:
\begin{conj}[{\cite[Conjecture 2.6]{koba01}}]\label{conj:spherical}
 There exists a compact complete pseudo-Riemannian manifold
of signature $(p,q)$ with constant sectional curvature $\kappa$
if and only if $(p,q,\kappa)$ satisfies one of the following conditions.
\begin{enumerate}
\ritem[1)] $\kappa>0$ and the pair $(p,q)$ is in the list below:
\begin{center}
\begin{tabular}{c||c|c|c|c|c}
$p$ & $\N$ & $0$ & $1$ & $3$ & $7$  \\ \hline
$q$ & $0$  & $\N$ & $2\N$ & $4\N$ & $8$
\end{tabular}
\end{center}
\ritem[2)] $\kappa=0$ and the pair $(p,q)$ is arbitrary.
\ritem[3)] $\kappa<0$ and the pair $(p,q)$ is in the list below:
\begin{center}
\begin{tabular}{c||c|c|c|c|c}
$p$ & $0$  & $\N$ & $2\N$ & $4\N$ & $8$  \\ \hline
$q$ & $\N$ & $0$ & $1$ & $3$ & $7$
\end{tabular}
\end{center}
\end{enumerate}
\end{conj}
We say $M$ is {\it complete} if the geodesic is defined for
all the time interval $(-\infty,\infty)$.
A complete space form of signature $(p,q)$ with constant
sectional curvature
$\kappa$ is a Clifford-Klein form of the following symmetric spaces:
\begin{enumerate}
\ritem[1)] $(\kappa>0)\quad O(p+1,q)/O(p,q)\quad(p\geq 2)$\quad or \quad
$\widetilde{O(2,q)}/\widetilde{O(1,q)}\quad(p=1)$.
\ritem[2)] $(\kappa=0)\quad O(p,q)\ltimes\R^{p+q}/O(p,q)$.
\ritem[3)] $(\kappa<0)\quad O(p,q+1)/O(p,q)\quad(q\geq 2)$\quad or \quad
$\widetilde{O(p,2)}/\widetilde{O(p,1)}\quad(q=1)$.
\end{enumerate}
Thus, Conjecture~\ref{conj:spherical} gives a conjectural
answer to Problem A for the above symmetric spaces.

The `if' part of Conjecture~\ref{conj:spherical} has been
proved by Borel, Kulkarni and Kobayashi
(Theorem~\ref{thm:borel} and Corollary~\ref{cor:admit_cck}).
Obvious cases for this are the case $q=0$ in (1) because the sphere $S^p$
itself is compact, and the case (2) because
$\R^{p+q}/\Z^{p+q}$ is a compact flat space form.
The `only if' part remains unsolved though
in spite of some recent
progress that provides evidence
supporting the conjecture
(Corollaries \ref{cor:spherical}, 
\ref{cor:benoist} and \ref{cor:grassman:parity}).

Conjecture~\ref{conj:spherical} is known to be true, particularly 
for the Riemannian and Lorentzian cases (i.e. $p\leq 1$ or $q\leq 1$).
To be more precise:
\begin{fact}[\cite{calabi},\cite{kulkarni},\cite{klingler}]\label{fact:lorentz2}
Concerning $n$-dimensional Lorentz space forms of sectional curvature $\kappa$:
\begin{enumerate}
\ritem[1)] $(\kappa>0)$ Compact forms never exist.
\ritem[2)] $(\kappa=0)$ Compact forms always exist.
\ritem[3)] $(\kappa<0)$ Compact forms exist if and only if $n$ is odd.
\end{enumerate}
\end{fact}
This fact can be proved as special cases of
Corollary~\ref{cor:spherical} $(\kappa>0)$, and Corollaries~\ref{cor:admit_cck}
and \ref{cor:grassman:parity} $(\kappa<0)$.
The case $\kappa=0$ is clear, as we have already seen that
$\Z^{k+1}\backslash\R^{k+1}$ is a compact flat space form.
Note that compactness implies completeness for the Lorentz space form
by a theorem of Klingler
 (\cite{klingler}).

Section~\ref{sec:spherical} will summarize the current status
of Conjecture~\ref{conj:spherical} for general
signature $(p,q)$, and will 
explain its `if' part including a complete
proof of Theorem~\ref{thm:spherical78}, namely, the existence of
a compact space form of signature $(p,q)=(7,8)$.

\subsection{Tangent space version of space form problem}
Associated to a reductive Lie group $G$ with
a Cartan decomposition $\g=\fk+\p$,
one defines the Cartan motion group $G_\theta:=K\ltimes \p$
as a semidirect group.
Likewise, 
one can associate a (non-reductive) symmetric space $G_\theta/H_\theta$
to a reductive symmetric space $G/H$,
which we say the {\it tangential symmetric space}.

Problem A for the tangential symmetric space $G_\theta/H_\theta$ 
turns out to be much simpler
than that for the reductive symmetric space $G/H$. In fact,
we can solve a `tangential' version of `space form conjecture'
(Conjecture~\ref{conj:spherical}).
To formulate the result, we recall the 
Hurwitz-Radon number $\rho(n)$ (see \cite{hurwitz}, \cite{radon})
defined as follows:
A positive integer $n$ is uniquely written as
$n=u\cdot 2^{4\alpha+\beta}$,
where $u, \alpha$ and $\beta$ are non-negative integers,
$u$ is odd and $\beta\leq 3$. Then, the Hurwitz-Radon number is given by
\begin{equation}
\rho(n):=8\alpha+2^\beta.\label{eq:df:rho}
\end{equation}
It is useful to set $\rho(0):=\infty$.

Let $G=O(p+1,q)$ and $H=O(p,q)$.
Then Theorem~\ref{thm:tangential} asserts that
the tangential symmetric space $G_\theta/H_\theta$ has a compact Clifford-Klein form if and only if $p<\rho(q)$.
Then, we obtain the
 list of the pair $(p,q)$ such that
$G_\theta/H_\theta$ has a compact Clifford-Klein form.
\begin{center}
\begin{tabular}{c||c|c|c|c|c|c|c|c|c|c}
$p$ & $\N$ & $0$  & $1$   & $2$   & $3$   & $4$   & $5$   & $6$   & $7$ & $\cdots$\\ \hline
$q$ & $0$ & $\N$ & $2\N$ & $4\N$ & $4\N$ & $8\N$ & $8\N$ & $8\N$ & $8\N$ & $\cdots$
\end{tabular}
\end{center}
We may observe that the above list is different from, but somewhat close to
the original
space form conjecture (see Conjecture~\ref{conj:spherical}
for $\kappa>0$).

\section{Methods}\label{sec:methods}
This section provides a brief account of the methods
to investigate Problem A, namely,
the existence problem of compact Clifford-Klein
forms of symmetric spaces.
\subsection{Generalized concept of properly discontinuous
actions}
Putting aside compactness for the moment,
we concentrate on the condition that the double coset space
$\Gamma\backslash G/H$ becomes a Clifford-Klein form of the
symmetric space $G/H$ with emphasis on non-compact $H$:
\begin{enumerate}
\ritem[(a)] The $\Gamma$-action on $G/H$ is fixed point free.
\ritem[(b)] The $\Gamma$-action on $G/H$ is properly discontinuous.
\end{enumerate}
Under the condition (b), the condition (a) is fulfilled if $\Gamma$
is torsion free. Furthermore, if $\Gamma$
is finitely generated and linear then $\Gamma$ contains
a torsion free subgroup of finite index (\cite{selberg}).
Thus, the condition (a) is not very serious, provided the
condition (b) is satisfied.
In the case $G/H$ is a Riemannian symmetric space, the condition (b)
is automatically satisfied for
any discrete subgroup $\Gamma$ of $G$.
However, for a non-Riemannian symmetric space $G/H$ where $H$ is
not compact,
the condition (b) is not always satisfied, and the resulting
double coset space
$\Gamma\backslash G/H$ may not be Hausdorff with
respect to the quotient
topology.

Thus, in the study of Problem A for non-Riemannian symmetric spaces,
it is crucial to
understand why the condition (b) may fail.

It is in general a hard problem to find an `efficient' criterion
for the condition (b). The first breakthrough was done in
\cite{koba88}, where the strategy is summarized as:
\begin{enumerate}
\itemindent=25pt
\ritem[Point 1)] Forget the original assumption
that $\Gamma$ is a discrete subgroup.
\ritem[Point 2)] Reformulate the problem inside the
group $G$ so that $\Gamma$ and $H$ play a symmetric role in $G$,
instead of the traditional approach that treats directly
the $\Gamma$-action on the homogeneous space $G/H$.
\ritem[Point 3)] Reformulate the properly discontinuous
condition in terms of finite dimensional representation theory.
\end{enumerate}
In order to emphasize Point (1), we shall use the notation $L$
instead of $\Gamma$.

Based on this strategy, the criterion of proper actions
(see Lemma~\ref{la:reductive_subgroup}) was proved under the assumption
that both $L$ and $H$ are reductive subgroups of a reductive Lie
group $G$, where Point (3) was pursued by using the standard theory
of parabolic subalgebras of a reductive Lie algebra $\g$.

The proof in \cite{koba89} shows that even the group structure
of $H$ and $L$ is not important. During his stay at the Institute
for Advanced Study at Princeton invited by Armand Borel in
1991--1992, the first author studied the deformation
of crystallographic groups for $G/H$ with non-compact $H$.
This study led us to
the following simple concept:

\begin{dfn}[{\cite[Definition 2.1.1]{koba96}}]\label{dfn:proper}
Let $L$ and $H$ be subsets of a locally compact group $G$.
\begin{enumerate}
\ritem[1)]
We say the pair $(L,H)$ is {\it proper} in $G$, denoted by $L\proper H$,
if $L\cap SHS$ is relatively compact for any compact set $S$ in $G$.
\ritem[2)]
We say the pair $(L,H)$ is {\it similar} in $G$, denoted by $L\similar H$,
if there exists a compact set $S$ in $G$ such that
$L\subset SHS$ and $H\subset SLS$.
\end{enumerate}
\end{dfn}
The point here is that we do not require $L$ and $H$ to be subgroups.
However, our motivation of introducing $\proper$ goes back to the following
observation:
\begin{la}\label{la:pd}
Suppose $H$ is a closed subgroup of a locally compact group $G$,
and $\Gamma$ is a discrete subgroup of $G$. Then
the natural $\Gamma$-action on $G/H$ is properly discontinuous
if and only if $\Gamma\proper H$.
\end{la}
\begin{rem}
Lemma~\ref{la:pd} can be generalized as follows:
{\it Suppose $L$ is a closed subgroup of $G$. Then
the natural $L$-action on $G/H$ is proper in the sense of Palais
\cite{palais} if and only if $L\proper H$.}
We note that a proper action of a discrete subgroup is properly
discontinuous, and vice versa.
\end{rem}
We recall from {\cite[Section 2]{koba96}}
some basic properties of the relations $\similar$ and
$\proper$:
\begin{la}\label{la:proper}
Let $L,L'$ and $H$ be subsets of a locally compact group $G$.
\begin{enumerate}
\ritem[1)] The relation $\similar$ is an equivalence relation.
\ritem[2)] If $L\similar L'$, then $L\proper H \Leftrightarrow L'\proper H$.
\ritem[3)] $L\proper H \iff H\proper L$.
\end{enumerate}
\end{la}
\begin{rem}\label{rem:proper}
1) If $H\similar G$, then
$L\proper H$
if and only if $L$ is compact.

2) If $H$ is compact, then $L\proper H$ for
any subset $L$ in $G$.
\end{rem}
Remark~\ref{rem:proper} has the following three interesting consequences:
\begin{prop}\label{prop:proper}
 Let $H$ be a closed subgroup of a locally compact group $G$.
\begin{enumerate}
\ritem[1)] Only a finite
subgroup $\Gamma\subset G$ can act on $G/H$
properly discontinuously if $H\similar G$
(the Calabi-Markus phenomenon; see Theorem~\ref{thm:calabi2}).
\ritem[2)]Any discrete subgroup $\Gamma\subset
G$ acts on the homogeneous space $G/H$ properly discontinuously
if $H$ is a compact subgroup of $G$.
\ritem[3)] Suppose $\Gamma$ is a cocompact discrete subgroup of $G$. Then the
$\Gamma$-action on $G/H$ is properly discontinuous if and only if
$H$ is compact.
\end{enumerate}
\end{prop}

The rest of this subsection collects some further basic
properties of the relations $\proper$ and $\similar$ together with another
relation $\ssimilar$, that
will be used in Section~\ref{sec:tangential}.
\begin{la}\label{la:cone}
Let $G$ be a Euclidean space, and let $L$ and $H$ be closed cones
of $G$. Then the following equivalences hold.
\begin{enumerate}
\ritem[1)] $L\proper H \iff L\cap H=\z$.
\ritem[2)] $L\similar H \iff L=H$.
\end{enumerate}
\end{la}
\begin{proof}
1) The implication $L\proper H \Rightarrow L\cap H=\z$ is obvious
from Definition~\ref{dfn:proper}.
Conversely, let us prove $L\not\proper H \Rightarrow L\cap H\neq\z$.
We find sequences 
$\{l_n\}\subset L,\ \{h_n\}\subset H$
and a positive number $M>0$ such that
$$ \sup_{n}\| l_n-h_n\|\leq M,\quad \lim_{n\rightarrow\infty}
\|l_n\|=\infty.$$
Taking a subsequence, if necessary, we may and do assume
the limit $l:=\lim\frac{l_n}{\|l_n\|}$ exists.
Since
$$ \left\| \frac{l_n}{\|l_n\|} -\frac{h_n}{\|l_n\|} \right\|=
\frac{\|l_n-h_n\|}{\|l_n\|}\leq \frac{M}{\|l_n\|} \rightarrow 0
\qquad (n\rightarrow \infty),$$
we have $\lim\frac{h_n}{\|l_n\|}=l$. As both $H$ and $L$ are
closed cones, we conclude $l\in L\cap H$. Hence, $L\cap H \neq \z$.

2) Obviously, $L=H\Rightarrow L\similar H$.
Let us prove $L\similar H\Rightarrow L=H$. We take a compact set
$S$ such that $L\subset H+S$. For any $l\in L$ and any
$n\in\N$, we have $nl\in L$, so that we find
$h_n\in H$ and $s_n\in S$ satisfying
$nl=h_n+s_n$. Thus we have
$$ \frac{nl-s_n}{n}=\frac{h_n}{n}\in H.$$
Taking the limit as $n$ tends to $\infty$, we have
$l\in H$ because $H$ is closed.
This means $L\subset H$. Similarly, $H\subset L$.
Thus, we have proved (2).
\end{proof}
Now we introduce another relation $\ssimilar$, which
is far stronger than $\similar$. For actual
applications, we shall need the case where $L$ and $H$ are subgroups.

\begin{dfn}
Let $L$ and $H$ be subsets of a locally compact group $G$.
We say the pair $(L,H)$ is {\it (right) strongly similar} in $G$,
denoted by $L\ssimilar H$, if there exist $g\in G$ and
a compact set $S$ in $G$ such that $L\subset gHS$ and $H\subset gLS$.
\end{dfn}
Clearly, $\ssimilar$ is an equivalence relation.
The following lemma is an obvious consequence of the definition:
\begin{la}\label{la:quotient_cpt}
Let $G$ be a locally compact group, and let $L\supset H$
be two closed subgroups. Then, $L/H$ is compact if and only if
$L\ssimilar H$.
\end{la}
However, the relation $\ssimilar$ will be useful later, when
there is no inclusive relation between $L$ and $H$.
We pin down some basic properties of $\ssimilar$ (and $\similar, \proper$):
\begin{la}\label{la:ssimilar}
Let $L,L'$ and $H$ be closed subgroups of a locally compact group $G$.
Assume that $L\ssimilar L'$, then the following equivalence holds:
$$ L\backslash G/H \text{ is compact} \iff L'\backslash G/H
\text{ is compact.}$$
\end{la}

\begin{la}\label{la:inclusion} Let $G_1$ be a locally compact group, and $G_2$
its closed subgroup.
Suppose $L$ and $H$ are subgroups of $G_2$. Then,
\begin{enumerate}
\ritem[1)] $\properin{L}{H}{G_1}\Rightarrow \properin{L}{H}{G_2}$.
\ritem[2)] $\similarin{L}{H}{G_2}\Rightarrow \similarin{L}{H}{G_1}$.
\ritem[3)] $\ssimilarin{L}{H}{G_2}\Rightarrow \ssimilarin{L}{H}{G_1}$.
\end{enumerate}
\end{la}

\newcommand{\oL}{\overline{L}}
\newcommand{\oH}{\overline{H}}
\newcommand{\oS}{\overline{S}}

\begin{la}\label{la:bieber_normal}
Let $G$ be a locally compact group and $N$ a normal subgroup.
We write the quotient map as
$\varpi:G\longrightarrow G/N$.
For two subsets $\oL$ and $\oH$ in $G/N$, we set
$L:=\varpi^{-1}(\oL)$ and $H:=\varpi^{-1}(\oH)$. Then,
\begin{equation*}
\ssimilarin{L}{H}{G}\iff\ssimilarin{\oL}{\oH}{G/N}.
\end{equation*}
\end{la}
\begin{rem}
An analogous statement of Lemma~\ref{la:bieber_normal}
is true if we replace $\ssimilar$ by $\similar$. However,
an analogous statement for $\proper$ is not true. In fact,
$\oL\proper \oH$ does not imply $L\proper H$.
\end{rem}

The proof of Lemmas~\ref{la:quotient_cpt}, \ref{la:ssimilar}
and \ref{la:inclusion}
are straightforward from the definition. Let us show
Lemma~\ref{la:bieber_normal}.
\begin{proof}[Proof of Lemma~\ref{la:bieber_normal}]
$\Rightarrow$)
Take $g\in G$ and a compact set $S$ in $G$ such that
$L\subset g H S$ and $H\subset g L S$.
Applying $\varpi$, we have
$\oL\subset\varpi(g)\oH \varpi(S)$ and $\oH\subset \varpi(g)\oL\varpi(S)$.
Since $\varpi(S)$ is compact, $\ssimilarin{\oL}{\oH}{G/N}$.

$\Leftarrow$) Take $\overline{g}\in G/N$ and a compact
set $\oS$ in $G/N$ such that
$\oL\subset \overline{g}\,\oH\,\oS$ and $\oH\subset \overline{g}\,\oL\,\oS$.
We take $g\in \varpi^{-1}(\overline{g})$ and a compact set $S$ in $G$ such that $\oS\subset \varpi(S)$. 
Then,
$$L=\varpi^{-1}(\oL)\subset\varpi^{-1}(\varpi(g)\oH\varpi(S))=gHS.$$
Similarly, we have
$H\subset gLS$.
Thus, the implication `$\Leftarrow$' has been proved.
\end{proof}

\subsection{Criterion of properness for reductive groups}
This subsection explicates a criterion of $\proper$ and $\similar$
for subsets of a linear reductive Lie group $G$.
Our criterion is `computable' by means of the Cartan projection,
and will play
a key role in the study of Clifford-Klein forms of
non-Riemannian symmetric spaces $G/H$.
An analogous result for the Cartan motion group will be given
in Section~\ref{sec:tangential}.

Let $G$ be a linear reductive Lie group, and 
$K$ a maximal compact subgroup. We write $\theta$ for the corresponding
Cartan involution, and
write $\g=\fk+\p$ for the
Cartan decomposition of the Lie algebra $\g$.
Choose a maximal
abelian subspace $\fa$ in $\p$.
The dimension of $\fa$ is 
independent of the choice, and is
called the {\it real rank} of $G$, denoted by $\rrank G$.
 Then the Cartan decomposition
$$G=K\exp(\p)=K\exp(\fa)K$$
has the following property:
any element $g\in G$ can be written as $g=k_1e^Xk_2$
for some $k_1,k_2\in K$,
and the element $X$ of $\fa$
is uniquely determined
up to the conjugation by the Weyl group $W_G:=N_G(\fa)/Z_G(\fa)$.
The correspondence $g\mapsto X$ induces the {\it Cartan projection}
$$\nu:G\rightarrow \fa/W_G.$$
For a subset $L$ of $G$, the
image of the Cartan projection is described as a $W_G$-invariant
subset as follows:
\begin{align}\label{eq:df:fl}
\fa(L)&:=\set{X\in\fa}{k_1e^Xk_2\in L\quad \text{for some $k_1,k_2\in K$}}\\
&=\set{X\in\fa}{e^X\in KLK}.\notag
\end{align}

Here is a criterion of properness $\proper$ and similarity $\similar$
for reductive Lie groups.
\begin{la}\label{la:reductive}
For any subsets $L$ and $H$ in a linear reductive Lie group $G$, we have
\begin{enumerate}
\ritem[1)] $\properin{L}{H}{G} \iff \properin{\fa(L)}{\fa(H)}{\fa}$.
\ritem[2)] $\similarin{L}{H}{G} \iff \similarin{\fa(L)}{\fa(H)}{\fa}$.
\end{enumerate}
Here, we regard $\fa\simeq\R^n$ as an abelian Lie group in the
right-hand side.
\end{la}
\begin{proof} See \cite[Theorem 1.1]{koba96} for a proof in the same circle
of ideas of \cite{koba89}. See also
\cite[Th\'eor\`eme 5.2]{benoist} for a proof where $G$ is a
reductive group over a local field.
\end{proof}
\begin{rem}
There is mysterious similarity between the criterion (1) of
Lemma~\ref{la:reductive} for proper actions and the criterion
for the non-existence of continuous spectrum in the irreducible
decomposition of unitary representations (\cite{koba98d}).
See \cite[Section~4D]{koba05} for details.
\end{rem}

In this paper, by a {\it reductive subgroup}, we shall mean a closed
subgroup $H$ which has
the following two properties:
\par\noindent
(a) $H$ has at most finitely many connected components.
\par\noindent
(b) $H$ is stable under some Cartan involution of $G$.

Then
$H$ itself is a reductive Lie group.
We say the resulting homogeneous space
$G/H$ is of {\it reductive type}. A reductive (or semisimple) symmetric
space is a typical example.

Suppose $H$ is a reductive subgroup of $G$,
and we take a maximal split abelian 
subspace $\fa'$ of the Lie algebra $\h$ of $H$. Then, there is $g\in G$ such that $\Ad(g)\fa'$
is contained in the fixed maximal split abelian subspace $\fa$ for $G$.
We set
\begin{equation}
\fa_H:=\Ad(g)\fa'.\label{eq:df:fah}
\end{equation}
Then $\fa_H$ is well-defined up to the
action of the Weyl group $W_G$. If $H$ is stable under the
fixed Cartan involution $\theta$, it is not hard to see
\begin{equation}
\fa(H)=W_G\cdot \fa_H\label{eq:ah}.
\end{equation}
For such $\theta$-stable subgroups $L$ and $H$, the above criterion of
properness (Lemma~\ref{la:reductive}) is of a much simpler
form.
\begin{la}[{\cite[Theorem~4.1]{koba89}}]\label{la:reductive_subgroup}
Suppose $L$ and $H$ are $\theta$-stable subgroups, with finitely
many connected components. Then
\begin{equation}
\properin{L}{H}{G} \iff \fa(L)\cap \fa(H)=\z.\label{eq:reductive_subgroup}
\end{equation}
\end{la}
\begin{ex}\label{ex:opq}
Let $G$ be an indefinite orthogonal group $O(p,q)\ (p\geq q)$.
We set $H_i:=E_{p+i,i}+E_{i,p+i} \in \g\ (1\leq i\leq q)$.
Then, the $q$-dimensional vector space
$$\fa:=\bigoplus_{i=1}^q \R H_i$$
is a maximal split abelian subspace.
For a subgroup $H=O(r,s)$, where $r\leq p$ and $s\leq q$,
the subset $\fa(H)$ is of the form:
$$ \fa(H)=\bigcup_{1\leq i_1<\dots<i_{\min(r,s)}\leq q}
\left(\bigoplus_{k=1}^{\min(r,s)}\R H_{i_k}\right).$$
\end{ex} 

\medskip
Calabi and Markus proved in \cite{calabi} that
there does not exist an infinite discontinuous group for
the Lorentz symmetric space $O(n+1,1)/O(n,1)$.
Named after their result, the non-existence of infinite
discontinuous group is sometimes referred to as the
{\it Calabi-Markus phenomenon}.
The Calabi-Markus phenomenon in the reductive case has been
studied by Calabi and Markus \cite{calabi}, Wolf \cite[64]{wolf},
Kulkarni \cite{kulkarni} and Kobayashi \cite{koba89},
and it has been settled completely in terms of the
rank condition as follows. This criterion
is obtained as an application of Lemma~\ref{la:reductive}.
\begin{thm}[\cite{koba89}]
\label{thm:calabi2}
Let $G/H$ be a homogeneous space of reductive type. Then the
following two conditions are equivalent:
\begin{enumerate}
\ritem[(i)] Any discontinuous group for $G/H$ is finite.
\ritem[(ii)] $\rrank G=\rrank H$.
\end{enumerate}
\end{thm}
\begin{proof}[Sketch of the proof]
The condition $\rrank H=\rrank G$ is equivalent to 
$\fa(H)=\fa$, and then equivalent to $H\similar G$ by
Lemma~\ref{la:reductive}.
Then, only a finite
group acts properly discontinuously on $G/H$ (Proposition~\ref{prop:proper}).
Conversely, if $\rrank G>\rrank H$, we can take $X\in \fa\setminus \fa(H)$.
Then, $\exp(\Z X)$ acts on $G/H$ properly discontinuously and freely
because of Lemma~\ref{la:reductive_subgroup} and $\R X\cap\fa(H)=\{0\}$.
\end{proof}
\begin{rem}
If $M$ is a compact topological space, then clearly only
a finite group  can act properly discontinuously on $M$.
The homogeneous space $G/H$ satisfying the rank condition (ii)
behaves as if $G/H$ were a compact space with respect to
the transformation group $G$.
\end{rem}

Geometrically, a homogeneous space of reductive type is a good
example of pseudo-Riemannian manifold.
\begin{prop}\label{prop:pseudo-R}
A homogeneous space $G/H$ of reductive type has a $G$-invariant
pseudo-Riemannian structure of signature
\begin{equation}
(d(G)-d(H),\dim G-\dim H-d(G)+d(H)).\label{eq:sig}
\end{equation}
\end{prop}
For a Lie group $G$, we define
\begin{equation}
d(G):=\dim G - \dim K \label{eq:dfdg}\\
\end{equation}
where $K$ is a  maximal compact subgroup of $G$. 
Then, $d(G)$ is well-defined
because any two maximal compact subgroups are conjugate to each other,
and may be regarded as the
`non-compact dimension' of $G$. For a reductive Lie group
$G$ with a Cartan decomposition $\g=\fk+\p$, we note
$$d(G)=\dim\p.$$
Here are some examples.
\begin{center}
\begin{tabular}{c|c|c|c|c|c|c}
$G$ & $GL(n,\R)$ & $SO(p,q)$ & $SU(p,q)$ & $Sp(p,q)$ & $Sp(n,\R)$
& $GL(n,\C)$ \\ \hline
$d(G)$ & $\tfrac{n(n+1)}{2}$ & $pq$ & $2pq$ & $4pq$ & $n^2+n$ & $n^2$
\end{tabular}
\end{center}
\smallskip
\begin{proof}[Proof of Proposition~\ref{prop:pseudo-R}]
Without loss of generality, we shall assume that $H$ is $\theta$-stable.

The adjoint action $\Ad: G\rightarrow GL(\g)$ is completely reducible
when restricted to $H$. Then, we take
$\mathfrak{q}$ to be the $H$-invariant complementary subspace of $\h$
in $\g$.
Since the decomposition $\g=\h+\mathfrak{q}$ is $\theta$-stable,
we have 
\begin{equation}
\g=\fk+\p=(\fk\cap\h)+(\fk\cap\mathfrak{q})+(\mathfrak{p}\cap
\h)+(\p\cap\mathfrak{q}). \label{eq:four:decomp}
\end{equation}

We fix a $G$-invariant symmetric bilinear form $B$
on $\g$ such that $B|_{\fk\times\fk}$ is negative definite,
$B|_{\p\times\p}$ is positive definite,
and $\fk$ is orthogonal to
$\p$ with respect to $B$ (e.g.\ the Killing form if $\g$ is
semisimple). 
Then the restriction of $B$ to $\mathfrak{q}\times\mathfrak{q}$
is an $H$-invariant non-degenerate symmetric bilinear form of signature
$(\dim(\mathfrak{q}\cap\p),\dim(\mathfrak{q}\cap\fk))$,
which equals \eqref{eq:sig}. Then, the left $G$-translation
of this bilinear form
gives rise to a $G$-invariant pseudo-Riemannian structure
on $G/H$.
\end{proof}

\subsection{Construction of compact Clifford-Klein forms}
\label{subsec:construction}
So far, we have discussed properly discontinuous actions
by putting aside the compactness condition. Building on the
criterion of properness (Lemma~\ref{la:reductive} or
Lemma~\ref{la:reductive_subgroup}),
we shall study the existence problem of compact Clifford-Klein
forms (Problem A).

A cocompact discrete subgroup $\Gamma$ of $G$ cannot act properly
discontinuously on the homogeneous space $G/H$ unless $H$ is compact
(Proposition~\ref{prop:proper}(3)).
This means that a crystallographic group for $G/H$ must be `smaller' than
a uniform lattice of $G$.
Instead of $G$ itself, a uniform lattice of a subgroup $L$
may be a candidate of a crystallographic group for $G/H$.
This idea leads us to the following
simple method that constructs a compact Clifford-Klein form of
$G/H$ with $H$ non-compact.
\begin{thm}[{\cite[Theorem 4.7]{koba89}}]\label{thm:sufficient}
Let $G$ be a Lie group, and $H$ its closed subgroup.
Then the homogeneous space $G/H$ admits a compact Clifford-Klein
form if there is a subgroup $L$ with the following three properties:

\smallskip
\noindent
\begin{tabular}{ll}
{\rm \eqref{thm:sufficient}(a)}& $\properin{L}{H}{G}$ 
(see Definition~\ref{dfn:proper}),\\
{\rm \eqref{thm:sufficient}(b)}& $L\backslash G/H$ is compact in the quotient topology,\\
{\rm \eqref{thm:sufficient}(c)}& $L$ has a cocompact, torsion free discrete subgroup $\Gamma$.
\end{tabular}
\end{thm}
\stepcounter{equation}
\begin{proof}[Sketch of the proof]
From (a) and (c), $\Gamma$ acts on $G/H$ properly discontinuously and freely.
Then, $\Gamma\backslash G/H$ is a compact Clifford-Klein
form of $G/H$ by (b) and (c).
\end{proof}
\begin{rem}
If $L$ is a linear group, then one can drop the 
torsion free assumption in (c). In fact, a finitely generated
linear group $\Gamma$ contains a torsion free subgroup
$\Gamma'$ of finite index by a result of Selberg \cite{selberg}.
\end{rem}
In the reductive case, we have a much explicit formulation of
Theorem~\ref{thm:sufficient}.
In what follows, the role of $L$ and $H$ will be symmetric.
This symmetric role goes back to Lemma~\ref{la:proper}(3).

\begin{thm}[\cite{koba89}]\label{thm:red:suf}
Let $G$ be a reductive linear Lie group, and let $L$ and $H$ 
be $\theta$-stable subgroups with finitely many connected
components. If $L$ and $H$ satisfy the following two properties,
then both $G/H$ and $G/L$ admit compact Clifford-Klein forms.

\smallskip
\noindent
\begin{tabular}{ll}
{\rm \eqref{thm:red:suf}(a)} & $\fa(L)\cap\fa(H)=0$,\\
{\rm \eqref{thm:red:suf}(b)} & $d(L)+d(H)=d(G)$.
\end{tabular}
\end{thm}
\stepcounter{equation}
\stepcounter{equation}
\begin{proof}[Sketch of the proof]
The condition \eqref{thm:sufficient}(a) follows from the
criterion of $L\proper H$ (see Lemma~\ref{la:reductive_subgroup}),
while the condition \eqref{thm:sufficient}(c) always holds for reductive
linear Lie groups by Borel's result (see Theorem~\ref{thm:borel}).
The remaining condition \eqref{thm:sufficient}(b) follows from the
numerical criterion given in Lemma~\ref{la:dg}.
\end{proof}

\begin{la}[{\cite[Theorem 4.7]{koba89}, see also \cite[Theorem 9.1]{iozzi}}]
\label{la:dg}
Let $G$ be a linear reductive Lie group. Suppose
that connected closed subgroups $L$ and $H$ satisfy
$\properin{L}{H}{G}$. Then we have:
\begin{enumerate}
\ritem[1)] $d(L)+d(H)\leq d(G)$.
\ritem[2)] $d(L)+d(H)=d(G)$ if and only if $L\backslash G/H$
is compact.
\end{enumerate}
\end{la}
\begin{rem}\label{rem:la:dg:fmcc}
Clearly, Lemma~\ref{la:dg} still holds if
$L$ and $H$ have finitely many connected
components.
\end{rem}
\begin{rem}
It might be illuminative to see that 
Lemma~\ref{la:dg} applied to $L=\{e\}$ means:
\begin{equation}
G/H \text{ is compact} \iff d(G)=d(H). \label{eq:ghcpt}
\end{equation}
\end{rem}
\begin{cor}\label{cor:admit_cck}
The following semisimple symmetric spaces $G/H$ have
compact Clifford-Klein forms, as the triple $(G, H, L)$
fulfills the conditions \eqref{thm:red:suf}{\rm (a)} and {\rm (b)}.
Here, $n=1,2,\dots$.
\setcounter{symcnt}{1}
\begin{center}
\begin{tabu}{c|c|c}
& $G/H$ & $L$ \\ \Hline
\symml{SU(2,2n)/Sp(1,n)}{U(1,2n)}
\symml{SU(2,2n)/U(1,2n)}{Sp(1,n)}
\symml{SO(2,2n)/U(1,n)}{SO(1,2n)}
\symml{SO(2,2n)/SO(1,2n)}{U(1,n)}
\symml{SO(4,4n)/SO(3,4n)}{Sp(1,n)}
\symml{SO(4,4)/SO(4,1)\times SO(3)}{Spin(4,3)}
\symml{SO(4,3)/SO(4,1)\times SO(2)}{G_{2(2)}}
\symml{SO(8,8)/SO(7,8)}{Spin(1,8)}
\symml{SO(8,\C)/SO(7,\C)}{Spin(1,7)}
\symml{SO(8,\C)/SO(7,1)}{Spin(7,\C)}
\symml{SO^*(8)/U(3,1)}{Spin(1,6)}
\symml{SO^*(8)/SO^*(6)\times SO^*(2)}{Spin(1,6)}
\end{tabu}
\end{center}
\end{cor}
\medskip
In Table~\thesubsection, the role of $L$ and $H$ is symmetric in
the cases 1 and 2, and in the cases 3 and 4 (see Theorem~\ref{thm:red:suf}).

Kulkarni \cite{kulkarni} proved the existence of a compact
Clifford-Klein form of $G/H$ in the
cases 4 and 5. The cases
1, 2 and 3 were proved
in \cite{koba89}. The cases 6 and 7 are
given in \cite[Corollary~3]{koba92}.
The cases 11 and 12 were given in
\cite[Examples 5.15 and~5.16]{koba96acad}.
See also \cite[Corollary 4.7]{koba96acad}
for some few non-symmetric examples 
such as $G/H=SO(4,4n)/Sp(1,n)$ ($n=1,2,\dots$) that admit compact
Clifford-Klein forms.
We shall give a detailed proof of the cases $8$, $9$ and $10$
in Section~\ref{sec:spherical}.

\begin{rem}
By using local isomorphisms of Lie
groups $SO^*(8)\simeq SO(2,6)$ and $SO^*(2)\times SO^*(6)\simeq U(1,3)$,
we see that $SO^*(8)/U(3,1)$ (the case 11), $SO^*(8)/SO^*(6)\times SO^*(2)$
(the case 12) are infinitesimally isomorphic to $SO(2,6)/U(1,3)$ (the case 2).
\end{rem}
\begin{rem}
Suppose $L$ and $H$ are $\theta$-stable connected subgroups of
a connected reductive Lie group $G$. If $\fl+\h=\g$, or
equivalently, if $\fl_\C+\h_\C=\g_\C$, then the inclusion
$L\subset G$ induces a natural diffeomorphism:
$L/L\cap H\simeq G/H$ (\cite[Lemma 5.1]{koba94}).
If $L\cap H$ is furthermore compact and if $L$ is linear,
then all the assumptions of Theorem~\ref{thm:red:suf} are
satisfied, and consequently, $G/H$ admits a compact Clifford-Klein form.
See \cite[Example 5.2]{koba94} for some list
of the triple $(G,H,L)$ of real reductive linear
Lie groups satisfying $\fl+\h=\g$.
This gives an alternative proof of Corollary~\ref{cor:admit_cck}.
\end{rem}

To the best of our knowledge,
all of the semisimple symmetric spaces
that have been proved
to admit compact Clifford-Klein forms so far are infinitesimally isomorphic
to a direct product of those $G/H$ in Table \thesubsection,
Riemannian symmetric spaces and reductive group manifolds.
We venture to make the following:
\begin{conj}\label{conj:inv:suff}
If the homogeneous space $G/H$ of reductive type
admits a compact Clifford-Klein form, then there
exists a reductive subgroup $L$ satisfying the assumptions
of Theorem~\ref{thm:sufficient}.
\end{conj}

\begin{rem}
Conjecture~\ref{conj:inv:suff} does not
assert the following (false) statement:
{\it A crystallographic group for $G/H$ is contained in
a reductive subgroup $L$ satisfying the assumptions of
Theorem~\ref{thm:sufficient}.}
The deformation argument in
\cite{goldman}, \cite{koba98} and \cite{salein} gives counterexamples
of this (too strong) statement.
\end{rem}
\begin{rem}
A `tangential' analog of Conjecture~\ref{conj:inv:suff}
will be formulated and proved in Section~\ref{sec:tangential}
(see Theorem~\ref{thm:carmot:criterion}).
\end{rem}

\subsection{Obstruction I --
Calabi-Markus phenomenon}\label{subsec:obstruction1}

In Section~\ref{subsec:construction}, we have given a sufficient condition
for symmetric spaces to admit compact Clifford-Klein forms.
Conversely, Sections~\ref{subsec:obstruction1} to \ref{subsec:obstruction3}
will explain some obstructions for 
symmetric spaces to have compact forms.
By using these obstructions, we shall find a number of symmetric
spaces without cocompact discontinuous groups.
These examples support Conjecture~\ref{conj:inv:suff}.

Loosely, one might think of the reason why some $G/H$ cannot have a
compact Clifford-Klein form as follows:
\begin{enumerate}
\ritem[I)] $G/H$ does not admit a `large' discontinuous group.
\ritem[II)] $H$ does not attain its maximum of `non-compact dimension'.
\end{enumerate}

These ideas will be formulated into theorems in
Sections~\ref{subsec:obstruction1} and \ref{subsec:obstruction2},
respectively.
Such formulations are made possible because
we know an explicit criterion for $\proper$
(see Lemma~\ref{la:reductive}).

\medskip
We start with a first obstruction obtained
by
the criterion (see Theorem~\ref{thm:calabi2}) of the Calabi-Markus phenomenon
(\cite{koba89}):
\begin{thm}\label{thm:cr:calabi}
Suppose $G$ is a  reductive Lie group, and $H$ is its
reductive subgroup with equal real rank. Then,
$G/H$ does not admit a compact
Clifford-Klein form, unless $G/H$ itself is compact,
or equivalently, if $d(G)>d(H)$.
\end{thm}
Applying Theorem~\ref{thm:cr:calabi} to the `space form problem' (see
Section~\ref{subsec:space_form}), we obtain immediately:
\begin{cor}[\cite{wolf}]\label{cor:spherical}
The symmetric space $G/H=O(p+1,q)/O(p,q)$ does not admit compact
Clifford-Klein forms if $p\geq q>0$.
\end{cor}
\begin{proof}
As $\rrank G=\min(p+1,q)$ and $\rrank H=\min(p,q)$, they
coincide if $p\geq q$. Furthermore,
$G/H$ is not compact because $d(G)-d(H)=(p+1)q-pq=q>0$.
Thus, Corollary~\ref{cor:spherical} follows from Theorem~\ref{thm:cr:calabi}.
\end{proof}
\begin{rem}
The case $q=1$ in Corollary~\ref{cor:spherical} was proved
by Calabi-Markus \cite{calabi}. See Section~\ref{subsec:space_form2}
for the condition on $(p,q)$ for which $G/H$ is known to have
compact Clifford-Klein forms.
\end{rem}

The criterion for $\proper$
(Lemma~\ref{la:reductive}) strengthens Theorem~\ref{thm:cr:calabi}
if we consider
non-abelian discontinuous groups instead of abelian discontinuous groups
as follows.

Let $\Sigma\equiv\Sigma(\g,\fa)$ be the restricted root system,
and let $W_G$ be its Weyl group. We fix
a positive system $\Sigma^+\equiv\Sigma^+(\g,\fa)$, and denote
by $\fa_+$ the corresponding dominant chamber,
and write $w_0$ for the longest element of $W_G$. We set
\newcommand{\bp}{\mathfrak{b}_+}
$$ \bp:=\set{X\in\fa_+}{w_0X=-X}.$$

Suppose $H$ is a reductive subgroup of $G$,
and we write $\fa_H\subset \fa$ as in \eqref{eq:df:fah}.
\begin{thm}[{\cite{benoist}}]\label{thm:benoist}
 Let $G/H$ be a
homogeneous space of reductive type. Then the following
two conditions are equivalent:
\begin{enumerate}
\ritem[(i)] Any discontinuous group for $G/H$ is virtually abelian.
\ritem[(ii)] $\bp\subset w\fa_H$ for some $w\in W_G$.
\end{enumerate}
\end{thm}
Here, an abstract group is said to be {\it virtually abelian}
if it contains an abelian subgroup of finite index.
\begin{thm}[{\cite[Corollary 7.6]{benoist}}]\label{thm:benoist2}
Suppose $G$ is a reductive Lie group, and $H$ its reductive subgroup.
If $\bp\subset w\fa_H$ for some $w\in W_G$, then $G/H$ does not
admit a compact Clifford-Klein form, unless $G/H$ itself is compact.
\end{thm}

Theorem~\ref{thm:calabi2} concerns with infinite discontinuous
groups, whereas Theorem~\ref{thm:benoist} concerns with
non-abelian discontinuous groups. Here is a comparison:
\begin{center}
\begin{tabular}{ccc}
Theorem~\ref{thm:calabi2}(i) & $\Rightarrow$ & Theorem~\ref{thm:benoist}(i)\\
$\Updownarrow$& &$\Updownarrow$\\
Theorem~\ref{thm:calabi2}(ii) & $\Rightarrow$ & Theorem~\ref{thm:benoist}(ii)
\end{tabular}
\end{center}
We note that if $\bp=\fa_+$, then all of the above four conditions
are equivalent. On the other hand, if $\bp\subsetneq\fa_+$, then
Theorem~\ref{thm:benoist2} is stronger than Theorem~\ref{thm:cr:calabi}.
We remark that $\bp\subsetneq\fa_+$ occurs when
$\Sigma$ is of type $A_n (n\geq 2)$, $D_{2m+1}$ or $E_6$.
\begin{ex}
1) Suppose $\Sigma(\g,\fa)$ is of type $A_{n-1}$.
We choose a positive system $\Sigma^+$ to be
$\set{e_i-e_j}{1\leq i<j\leq n}$ so that
 $\fa_+=\set{(x_1,\dots,x_n)\in \R^n}{x_1\geq\dots\geq x_n,
\sum x_j=0}$. Then $w_0\neq -1$ if and only if $n\geq 3$, and we have
$$
\bp=
\begin{cases}
\set{(x_1,\dots,x_m,-x_m,\dots,-x_1)}{x_1\geq\dots\geq x_m\geq 0}
& (n=2m),\\
\set{(x_1,\dots,x_m,0,-x_m,\dots,-x_1)}{x_1\geq\dots\geq x_m\geq 0}
& (n=2m+1).
\end{cases}
$$
2) Suppose $\Sigma(\g,\fa)$ is of type $D_n$.
We choose $\Sigma^+(\g,\fa)=\set{e_i\pm e_j, e_k}
{1\leq i<j\leq n, 1\leq k\leq n}$ so that
$\fa_+=\set{(x_1,\dots,x_n)\in\R^n}{x_1\geq\dots\geq x_{n-1}\geq |x_n|}$.
Then $w_0\neq -1$ if and only if $n$ is odd, and
$$
\bp=
\begin{cases}
\set{(x_1,\dots,x_{n-1},0)}{x_1\geq\dots\geq x_{n-1}\geq 0}
& (n \text{ is odd}),\\
\fa_+ & (n \text{ is even}).
\end{cases}
$$
\end{ex}
\begin{cor}\label{cor:benoist}
The following symmetric spaces $G/H$ do not admit
compact Clifford-Klein forms.
\setcounter{symcnt}{1}
\begin{center}
\begin{tabu}{c|c}
& $G/H$ \\ \Hline
\symm{SL(2n,\R)/Sp(n,\R)}
\symm{SL(2n+1,\R)/SO(n,n+1)}
\symm{SO(2n+1,2n+1)/SO(2n,2n+1)}
\symm{SL(2n,\C)/Sp(n,\C)}
\symm{SO(4n+2,\C)/SO(4n+1,\C)}
\symm{E_{6,\C}/F_{4,\C}}
\end{tabu}
\end{center}
\end{cor}
\subsection{Obstruction II -- maximality of non-compactness}
\label{subsec:obstruction2}
If two subgroups $H$ and $L$ satisfy $H\similar L$,
then any discontinuous group for $G/H$
is also a discontinuous group for $G/L$. Then, the comparison of 
two homogeneous spaces $G/H$ and $G/L$
gives rise to a second necessary condition for
the existence of compact Clifford-Klein forms for
$G/H$:
\begin{thm}[{\cite[Theorem 1.5]{koba92n}}]\label{thm:non-cpt-dim}
Let $G/H$ be a homogeneous space of reductive type.
If there exists a reductive subgroup $L$ in $G$ such that
$L\similar H$ and $d(L)>d(H)$, then
$G/H$ does not admit compact Clifford-Klein forms.
\end{thm}
\begin{rem}
It is readily seen that
the subgroup $L$ in this Theorem
becomes
an obstruction for the existence of the subgroup
$L$ in Theorem~\ref{thm:sufficient}.
In this sense, Theorem~\ref{thm:non-cpt-dim}
affords evidence for Conjecture~\ref{conj:inv:suff}.
\end{rem}

A key step of proving Theorem~\ref{thm:non-cpt-dim}
is to formulate a discrete analog of Lemma~\ref{la:dg}
in terms of the virtual cohomological dimension
$vcd(\Gamma)$ of an abstract group $\Gamma$
as follows:

\begin{la}[\cite{koba89}]\label{la:vcd}
Let $G/H$ be a homogeneous space of reductive type.
If a discrete subgroup $\Gamma\subset G$ satisfies $\Gamma\proper H$,
then we have:
\begin{enumerate}
\ritem[1)] $vcd(\Gamma)+d(H)\leq d(G)$.
\ritem[2)] $vcd(\Gamma)+d(H)= d(G)$ if and only if $\Gamma\backslash G/H$
is compact.
\end{enumerate}
\end{la}
\begin{rem}
A special case of Lemma~\ref{la:vcd}(2) with $H=\{e\}$ asserts that
$$ vcd(\Gamma)=d(G)$$
if $\Gamma$ is a cocompact discrete subgroup of $G$.
This result is due to Serre \cite{serre}.
\end{rem}
\begin{proof}[Proof of Theorem~\ref{thm:non-cpt-dim} (sketch)]
Let $\Gamma$ be a uniform lattice for $G/H$, then we have
$vcd(\Gamma)+d(H)=d(G).$ On the other hand, we have $\Gamma\proper L$,
thus $vcd(\Gamma)+d(L)\leq d(G)$. This contradicts to $d(L)>d(H)$.
\end{proof}
\begin{rem}\label{rem:smallL}
As one can see from the proof,
Theorem~\ref{thm:non-cpt-dim}
can be strengthened by replacing
the assumption $L\similar H$ by
$\fa(L)\subset \fa(H)$.
\end{rem}
\begin{rem}
Taking $L$ to be $G$ in Theorem~\ref{thm:non-cpt-dim},
we obtain an alternative proof of
Theorem~\ref{thm:cr:calabi}.
\end{rem}
The following symmetric space $G/H$ is a good test case, as it
contains many parameters. We note that $G/H$ becomes
a pseudo-Riemannian space form if $i=0$ and $k=1$.
\begin{cor}[{\cite[Example 1.7]{koba92n}}]\label{cor:grassman:ineq}
Let $G/H$ be a pseudo-Riemannian Grassmannian manifold
$$ O(i+j,k+l)/O(i,k)\times O(j,l).$$
Without loss of generality, we may and do assume $i\leq j,k,l$.
We assume that neither $H$ nor $G/H$ is compact, that is, $j,k,l>0$.
If $G/H$ admits compact Clifford-Klein forms, then $i=0$ and
$0<l\leq j-k$.
\end{cor}
\begin{proof}
Assume $G/H$ admits compact Clifford-Klein forms.
From Theorem~\ref{thm:cr:calabi}, we have $\rrank H <\rrank G$, that is,
$$i+\min(j,l) < \min(i+j,k+l).$$
In particular, $l<j$.
On the other hand, we define two subgroups $L_1$ and $L_2$ by
\begin{align*}
L_1&:=O(i+j,i+l) \subset O(i+j,k+l),\\
L_2&:=O(i+l,k+l) \subset O(i+j,k+l).
\end{align*}
Then, we have $\fa(H)=\fa(L_1)=\fa(L_2)$ (see Example~\ref{ex:opq}).
Thus $H\similar L_1\similar L_2$ in $G$ by Lemma~\ref{la:reductive}.
Now, applying Theorem~\ref{thm:non-cpt-dim}, we have
the following two inequalities
\begin{align*}
d(H)=ik+jl &\geq d(L_1)=(i+l)(k+l),\\
d(H)=ik+jl &\geq d(L_2)=(i+j)(i+l).
\end{align*}
Combining with $l<j$, we conclude that
$i=0$ and $l\leq j-k$.
\end{proof}
\begin{rem}
The same results as Corollary~\ref{cor:grassman:ineq} hold
if we consider indefinite unitary groups over $\C$ or $\bbH$
instead of $\R$, namely, if we replace $O(p,q)$ by
$U(p,q)$ or $Sp(p,q)$.
\end{rem}

\setbox0=\hbox{\scriptsize $r=\min(n,2p+1,2q+1)$}
\setbox1=\vbox to 10pt {\hbox to\wd0{\hfil$\mathfrak{so}^*(2r)$\hfil}\vskip-4pt\box0}
\setbox2=\hbox to 0pt {\hskip-60pt$\Bigl(\ $\lower 5pt \vbox to 16pt
{\footnotesize \hbox{$3p\leq 2n\leq 6p$}\hbox{\hskip15pt$n\geq 3$}}$\ \Bigr)$}

\begin{cor}[{\cite[Table 4.4]{koba92n} and \cite[Table 5.18]{koba96}}]
\label{cor:not_exist}
The following symmetric spaces $G/H$ do not admit
compact Clifford-Klein forms
as the pair $(G/H,L)$ satisfies
the assumptions of Theorem~\ref{thm:non-cpt-dim}
(or Remark~\ref{rem:smallL}).
In the table, we always assume $n=p+q$. We denote by $\floor{m}$ the
maximal integer that does not exceed $n$.
\setcounter{symcnt}{1}
\newcommand{\symmlc}[3]{\symm{#1$ & #3 & $#2}}
\newcommand{\symmlcv}[4]{$\thesymcnt$ \stepcounter{symcnt} & $#1$ &
#3 & $#2$\\#4}
\begin{center}
\begin{tabu2}{c|cr|c}{1.5}
& $G/H$ && $Lie(L)$\\ \Hline
\symmlc{SL(2n,\R)/SO(n,n)}{\mathfrak{sp}(n,\R)}{}
\symmlc{SU^*(2n)/SO^*(2n)}{\mathfrak{sp}(\floor{\tfrac{n}{2}},n-
\floor{\tfrac{n}{2}})}{}
\symmlc{SU(2n,2n)/SO^*(4n)}{\mathfrak{sp}(n,n)}{}
\symmlc{Sp(2n,\R)/U(n,n)}{\mathfrak{sp}(n,\C)}{}
\symmlc{SO(2n,2n)/SO(2n,\C)}{\mathfrak{u}(n,n)}{}
\symmlc{SO^*(2n)/SO^*(2p)\times SO^*(2q)}
         {\mathfrak{so}^*(2)+\mathfrak{so}^*(2n-2)}{$(p,q>1)$}
\symmlc{SL(n,\C)/SO(n,\C)}{\mathfrak{u}(\floor{\tfrac{n}{2}},
           n-\floor{\tfrac{n}{2}})}{}
\symmlcv{SO(n,\C)/SO(p,\C)\times SO(q,\C)}
        {\mathfrak{so}(n-1,\C)}{$(p,q>1)$}{[3pt]}
\symmlcv{SO^*(2n)/U(p,q)}{\raise 1pt \box1}{\box2}{[3pt]}
\symmlc{SL(2n,\C)/SU(n,n)}{\mathfrak{sp}(n,\C)}{}
\symmlc{Sp(n,\R)/Sp(p,\R)\times Sp(q,\R)}{\mathfrak{sp}(n,\R)}{}
\symmlc{E_{6(6)}/Sp(4,\R)}{\mathfrak{f}_{4(4)}}{}
\symmlc{E_{6(2)}/SU(4,2)\times SU(2)}{\mathfrak{so}^*(10)+\sqrt{-1}\R}{}
\symmlc{E_{7(7)}/SU(4,4)}{\mathfrak{e}_{6(2)}+\sqrt{-1}\R}{}
\symmlc{E_{7(7)}/SU^*(8)}{\mathfrak{so}^*(12)+\mathfrak{su}(2)}{}
\symmlc{E_{7(-5)}/SU(6,2)}{\mathfrak{e}_{6(-14)}+\sqrt{-1}\R}{}
\symmlc{E_{7(-25)}/SU(6,2)}{\mathfrak{e}_{6(-14)}+\sqrt{-1}\R}{}
\symmlc{E_{8(8)}/SO^*(16)}{\mathfrak{e}_{7(-5)}+\mathfrak{su}(2)}{}
\symmlc{E_{6,\C}/Sp(4,\C)}{\mathfrak{f}_{4,\C}}{}
\end{tabu2}
\end{center}
\end{cor}
Some remarks on Table~\thesubsection\ are in order.
\begin{enumerate}
\ritem[1)] In the case 6, there {\it exists} a compact Clifford-Klein form
for $(p,q)=(1,3)$ and $(3,1)$.
\ritem[2)] In the cases 6 and 8 with $pq$ even, and in the case 11,
there is an alternative proof by using Theorem~\ref{thm:cr:calabi}.
\ritem[3)] Taking this opportunity, we correct
typographical errors in \cite[Table 5.18]{koba96};
the pairs $(SL(2n,\C), SU(p,q))$ and $(SL(2n,\R), SO(p,q))$ there should
read $(SL(n,\C), SU(p,q))$ and $(SL(n,\R), SO(p,q))$, respectively.
\end{enumerate}

In what follows, we omit an explanation of
the technical terms
``associated pair'' and ``basic in $\epsilon$-family'',
but we just indicate how to find
$L$ systematically in most of the cases in Table~\thesubsection.
\begin{la}[{\cite[Theorem 1.4]{koba92n}}]
A reductive symmetric space $G/H$ admits a compact Clifford-Klein
form only if its associated symmetric pair $(G,H^a)$ is basic in
the $\epsilon$-family.
\end{la}

\subsection{Obstruction III
-- compatibility of three compact spaces}
\label{subsec:obstruction3}

A Clifford-Klein form $\Gamma\backslash G/H$ inherits local geometric
structure from $G/H$ through the covering map $G/H\rightarrow \Gamma\backslash
G/H$.
Such geometric structures may become some constraints
on the topology of Clifford-Klein forms.

Another constraint is derived from the fiber bundle structure
$G/H\rightarrow K/H\cap K$. The compatibility of 
these constraints yields a necessary condition for
the existence of compact Clifford-Klein form, formulated
as follows:
\begin{thm}[{\cite[Proposition 4.10]{koba89}}]\label{thm:rank}
Suppose $G$ is a reductive Lie group, and $H$ is its
reductive subgroup.
If $\rank G=\rank H$ and $\rank K> \rank H\cap K$,
then $G/H$ does not admit a compact Clifford-Klein form.
\end{thm}

Let $G_\C$ be a complexification of $G$, and
let $G_U$ be a compact real form of $G_\C$
(for example, $G_U=SU(n)$ if $G=SL(n,\R)$).
Likewise $H_U$ for $H$.
Thus, we have the following setting:
\begin{center}
\begin{tabular}{ccccc}
$G$ & $\subset$ & $G_\C$ & $\supset$ & $G_U$\\
$\cup$ & &$\cup$ & &$\cup$\\
$H$ & $\subset$ & $H_\C$ & $\supset$ & $H_U$
\end{tabular}
\end{center}
The compact homogeneous space $G_U/H_U$ plays a role of
a `real form' of $G_\C/H_\C$. 
The idea of analytic continuation leads us to:
\begin{la}[\cite{kobaono}]\label{la:homology}
There is a canonical $\C$-algebra homomorphism between
de Rham cohomology groups:
$$ \eta: H^*(G_U/H_U;\C)\rightarrow H^*(\Gamma\backslash G/H;\C).$$
If $\Gamma\backslash G/H$ is compact, then $\eta$
is injective. Furthermore, $\eta$ sends the characteristic classes for
$G_U/H_U$ to the corresponding ones for $\Gamma\backslash G/H$.
\end{la}

\begin{rem}
If $H=K$, a maximal compact subgroup of $G$,
then there exists a compact Clifford-Klein form
$\Gamma\backslash G/K$ (see Theorem~\ref{thm:borel}).
In this case, the image of
$\eta$ is isomorphic to the $(\g,K)$-cohomology group
$H^*(\g,K;\C)$ via the Matsushima-Murakami isomorphism (\cite{borelwallach}):
$$H^*(\Gamma\backslash G/H;\C)\simeq \bigoplus_{\pi\in\hat{G}}
\Hom_G(\pi,L^2(\Gamma\backslash G))\otimes H^*(\g,K;\pi_K),$$
where $\pi_K$ is the $(\g,K)$-module of an irreducible
unitary representation $\pi$ of $G$.
\end{rem}

\begin{proof}[Proof of Theorem~\ref{thm:rank} (sketch)]
The idea is to compare three compact spaces $G_U/H_U$, $K/H\cap K$
and $\Gamma\backslash G/H$.
The Euler number $\chi(G_U/H_U)\neq 0$ if $\rank G=\rank H$
by a theorem of Hopf and Samelson \cite{hopf}.
On the other hand, we can show $\chi(\Gamma\backslash G/H)=0$ by using
the fiber bundle structure $G/H \rightarrow K/H\cap K$, provided
$\rank K > \rank K\cap H$.
By Lemma~\ref{la:homology} $\chi(\Gamma\backslash G/H)=\eta(\chi(G_U/H_U))$.
Therefore, $\eta$ cannot be injective. Hence, $\Gamma\backslash G/H$
cannot be compact.
\end{proof}

Here are some examples of Theorem~\ref{thm:rank}:
\begin{cor}
A semisimple symmetric space $Sp(2n,\R)/Sp(n,\C)$ does not admit
a compact Clifford-Klein form.
\end{cor}
\begin{proof} Both the ranks of $Sp(2n,\R)$ and $Sp(n,\C)$
are $2n$, while the ranks of their
maximal compact subgroups
$U(2n)$ and $Sp(n)$ are $2n$ and $n$, respectively.
Therefore, Corollary follows from Theorem~\ref{thm:rank}.
\end{proof}
\begin{cor}\label{cor:grassman:parity}
If $jkl$ is odd, an indefinite Grassmannian manifold
$$ G/H = O(j,k+l)/O(k)\times O(j,l)$$
does not admit a compact Clifford-Klein form.
\end{cor}
\begin{proof}
The condition
$\rank G=\rank H$ amounts to
\begin{equation}
\floor{\tfrac{j+k+l}{2}} = \floor{\tfrac{k}{2}}+\floor{\tfrac{j+l}{2}}.
\label{eq:rank}
\end{equation}
On the other hand, the condition $\rank K > \rank H\cap K  $ amounts to
$$\floor{\tfrac{j}{2}}+\floor{\tfrac{k}{2}}+\floor{\tfrac{l}{2}} >
\floor{\tfrac{j}{2}}+\floor{\tfrac{k+l}{2}}.
$$
This holds if and only if $kl$ is odd. Under the assumption $kl$ is odd,
\eqref{eq:rank} is equivalent to that
$j$ is odd.
Thus, Corollary follows from Theorem~\ref{thm:rank}.
\end{proof}

\begin{rem}
Theorem~\ref{thm:tangential} implies particularly that
there do not exist compact Clifford-Klein forms
of the tangential symmetric space $G_\theta/H_\theta$
for $G/H=O(p+1,q)/O(p,q)$
if $p>0$ and $q$ is odd. This parity condition may be regarded
as a `tangential' analog of
Corollary~\ref{cor:grassman:parity} with $k=1,l=p$ and $j=q$.
\end{rem}

\def\mck{\mathcal K}
\def\mcl{\mathcal L}
\def\mcm{\mathcal M}
\def\E{E}
\newcommand{\mon}{_{\operatorname{mon}}}
\newcommand{\V}{^\vee}

\section{Space form problem}\label{sec:spherical}
By the ``space form problem'' we mean an (unsolved)
problem whether or not there exists
a compact complete pseudo-Riemannian manifold of general
signature $(p,q)$ with constant sectional curvature.
This section provides a proof of the `if' part of 
Conjecture~\ref{conj:spherical} to the space form
problem, together with a proof of the `if' part of
Conjecture~\ref{conj:complex_sphere} for complex version of
space form problem (namely the cases $8$ and~$9$ of
Theorem~\ref{thm:admit_cck}).
An algebraic machinery for this is the Clifford algebra
associated to indefinite quadratic form
that we develop in Sections \ref{subsec:Clifford}--\ref{subsec:auto_iso}.
This leads us to another family of finitely many homogeneous spaces
that admit compact Clifford-Klein forms arisen from
the spin representation (see Theorem~\ref{thm:spin_ck_form}).
We also complete the proof of the remaining case of
Theorem~\ref{thm:admit_cck}, namely, the case 10 in Table~\ref{subsec:admit_cck}.

\subsection{Space form conjecture}\label{subsec:space_form2}
This subsection summarizes the current status of
the `space form conjecture' (Conjecture~\ref{conj:spherical}).
First, we equip $\R^{p+q}$ with the pseudo-Riemannian metric
$$ds^2=dx_1^2+\dots+dx_p^2-dx_{p+1}^2-\dots-dx_{p+q}^2.$$ The
resulting pseudo-Riemannian manifold, denoted by $\R^{p,q}$,
is a flat space form with signature $(p,q)$.

Next, let $Q_{p,q}$ be a quadratic form on $\R^{p+q}$ defined by
$$ Q_{p,q}(x):=x_1^2+\dots+x_p^2-x_{p+1}^2-\dots-x_{p+q}^2.$$
We induce a pseudo-Riemannian structure on the hypersurface
$$ X(p,q):=\set{x\in\R^{p+1,q}}{Q_{p+1,q}(x)=1}$$
from $\R^{p+1,q}$. Then $X(p,q)$ is a space form
of signature $(p,q)$ with
sectional curvature $+1$. It is simply connected if $p\geq 2$.
\begin{rem}\label{rem:replace_k}
Note that $X(p,q)$ has signature $(q,p)$ with sectional curvature
$-1$, if we replace the pseudo-Riemannian metric $g$ by $-g$.
\end{rem}

Here are some well-known cases.
\begin{enumerate}
\ritem[1)] $\R^{3,1}$ is the Minkowski space.
\ritem[2)] $X(p,0)$ is the $p$-sphere $S^p$.
\ritem[3)] $X(0,q)$ is a hyperbolic space.
\ritem[4)] $X(p,1)$ is a de Sitter space.
\ritem[5)] $X(1,q)$ is an anti-de Sitter space.
\end{enumerate}

The indefinite orthogonal group
\begin{equation}
O(p,q)=\set{g\in GL(p+q,\R)}{Q_{p,q}(gx)=Q_{p,q}(x)\text{ for any
$x\in\R^{p+q}$}}\label{eq:df:opq}
\end{equation}
acts transitively on $X(p,q)$, so that
$X(p,q)$ is represented as the homogeneous space $O(p+1,q)/O(p,q)$.
We notice that $O(p+1,q)/O(p,q)$ is a semisimple symmetric space
of rank one. Then, it is plausible that
Problem A for this symmetric space has
the following solution:
\begin{conj}\label{conj:spherical:group}
The symmetric space $M=O(p+1,q)/O(p,q)$ admits compact Clifford-Klein forms
if and only if $(p,q)$ is in the list below:
\begin{center}
\begin{tabular}{c||c|c|c|c|c}
$p$ & $\N$ & $0$ & $1$ & $3$ & $7$  \\ \hline
$q$ & $0$  & $\N$ & $2\N$ & $4\N$ & $8$
\end{tabular}
\end{center}
\end{conj}
We note:
\vskip8pt
\noindent
\begin{tabular}{rll}
$\cdot$ & $p=0$ or $q=0$ & $M$ is a Riemannian symmetric space.\\
$\cdot$ & $p=1,3,7$ & $M$ is a pseudo-Riemannian symmetric space.
\end{tabular}
\vskip8pt
\begin{rem}
Conjecture~\ref{conj:spherical:group} is equivalent to
Conjecture~\ref{conj:spherical}. In fact, in the
case $\kappa=0$, Conjecture~\ref{conj:spherical} is true because
$\R^{p,q}$ has a uniform lattice $\Z^{p+q}$.
In the case $\kappa<0$, Conjecture~\ref{conj:spherical} is
reduced to the case $\kappa>0$ as we saw in Remark~\ref{rem:replace_k}.
Conjecture~\ref{conj:spherical} and Conjecture~\ref{conj:spherical:group}
are also true in the case
$p=1$ and $\kappa>0$: compact forms exist by Theorem~\ref{thm:admit_cck}
if $q$ is even, while compact forms do not exist by
Corollary~\ref{cor:grassman:parity} if $q$ is odd.

For the remaining case where $p\geq 2$ and $\kappa>0$,
the universal covering manifold of any complete space form of signature
$(p,q)$ with sectional curvature $\kappa$ is isometrically
(up to a positive scalar multiplication of metric) diffeomorphic
to the pseudo-Riemannian symmetric space $X(p,q)\simeq O(p+1,q)/O(p,q)$.
Thus, Conjecture~\ref{conj:spherical} is equivalent to
Conjecture~\ref{conj:spherical:group}.
\end{rem}

As for the `only if' part of Conjecture~\ref{conj:spherical:group},
to the best of our knowledge,
the following is all the cases that
the non-existence of compact space forms
has been proved so far.
\vskip8pt
\noindent
\begin{tabular}{cll}
(a)& $p\geq q>0$  &(see Corollary~\ref{cor:spherical}).\\
(b)& $p+1=q$ is odd\qquad & (see Corollary~\ref{cor:benoist}).\\
(c)& $pq$ is odd & (see Corollary~\ref{cor:grassman:parity}).
\end{tabular}
\vskip8pt

Conversely, the `if' part of the Conjecture~\ref{conj:spherical:group}
is true. Here is a brief account of each case.
\par\medskip\noindent
{\bf 1) $q=0$}; \quad $M=S^p$ itself is compact.
\par\medskip\noindent
{\bf 2) $p=0$};\quad $M$ is a Riemannian symmetric space and has a compact
Clifford-Klein form by Theorem~\ref{thm:borel}.

For the remaining three pseudo-Riemannian cases, let us see
$G/H$ admits compact Clifford-Klein forms by applying
Theorem~\ref{thm:red:suf}.
\par\medskip\noindent
{\bf 3) $(p,q)=(1,2n)$};\quad $G=O(2,2n),\quad H=O(1,2n)$.\\
Take a subgroup $L$ to be $U(1,n)$.
By choosing a suitable basis, we identify
$\fa\simeq \R^2$, on which
the Weyl group $W_G\simeq \mathfrak{S}_2\ltimes (\Z/2\Z)^2$
acts by permutation and change of
signatures of coordinates. Then one finds
$$\fa(H)=W_G\cdot \R
\begin{pmatrix} 1 \\ 0\end{pmatrix},\quad \fa(L)=W_G\cdot\R
\begin{pmatrix} 1 \\ 1\end{pmatrix},$$
and therefore $\fa(L)\cap\fa(H)=\z$.
On the other hand, we have $d(L)+d(H)=2n+2n=4n=d(G)$.
Now, it follows from Theorem~\ref{thm:red:suf} that a uniform lattice
of $L$ becomes a cocompact discontinuous group for $G/H$.
\par\medskip\noindent
{\bf 4) $(p,q)=(3,4n)$};\quad $G=O(4,4n),\quad H=O(3,4n)$.\\
Take $L=Sp(1,n)$.
Similarly to the previous case, one finds
$$\fa=\R^4,\quad
W_G\simeq \mathfrak{S}_4\ltimes(\Z/2\Z)^4,\quad
 \fa(H)=W_G\cdot \R
\begin{pmatrix} 1 \\ 1\\ 1\\ 0\end{pmatrix},\quad
\fa(L)=W_G\cdot\R\begin{pmatrix} 1 \\ 1\\ 1\\ 1\end{pmatrix},
$$
and so $\fa(L)\cap\fa(H)=\z$. Furthermore,
$d(L)+d(H)=12n+4n=16n=d(G)$.
Hence, $G/H$ admits a compact Clifford-Klein form by
Theorem~\ref{thm:red:suf}.

\medskip
The rest of this section is devoted to proving
Theorem~\ref{thm:spherical78},
namely, the existence of a compact Clifford-Klein form for the last
case:

\medskip
\noindent
{\bf 5) $(p,q)=(7,8)$};\quad $G=O(8,8),\quad H=O(7,8)$.\\
We shall take a subgroup $L$ to be $Spin(1,8)$. We need to embed
$Spin(1,8)$ into $G$ and then to find $\fa(L)\subset \fa$. For this
purpose, we study the Clifford algebra $C(p,q)$ for the indefinite
quadratic form $Q_{p,q}$ of signature $(p,q)$.
This will be carried out in Sections~\ref{subsec:Clifford} to
\ref{subsec:auto_iso}  based on an unpublished note
\cite{koba94c}, in the way that we need
for the proof of Theorems~\ref{thm:intro:complex_sphere} and
\ref{thm:spherical78}.

\subsection{Compact Clifford-Klein forms}
\label{subsec:spin_cck}
As was stated in Theorem~\ref{thm:admit_cck}, there exist compact Clifford-Klein
forms for the symmetric spaces $O(8,\C)/O(7,\C)$ and
$O(8,8)/O(7,8)$. We shall provide a unified proof of
these results. That is, the following theorem exhibits
a family of finitely many homogeneous
spaces $G/H$ that admit uniform lattices arisen from the
spin representation.
Then the homogeneous space $O(8,\C)/O(7,\C)$ and $O(8,8)/O(7,8)$
appear corresponding to
 $q=7$ and $8$, respectively
in Theorem~\ref{thm:spin_ck_form} below.

This family of homogeneous spaces will be treated
systematically
based on Clifford algebras of indefinite quadratic forms.
\begin{thm}\label{thm:spin_ck_form}
The  triples of Lie groups $(G,H,L) = (G(1,q), H(1,q),
 Spin(1,q))$ $(1\leq q\leq 8)$ satisfy the conditions \eqref{thm:red:suf}(a) and (b),
and therefore, both of the homogeneous spaces $G/H$ and $G/L$ admit compact
Clifford-Klein forms.
$$
\begin{array}{|c|ccc|ccccccc|} \hline
q& G & H & L & W && \hbox to 0pt {\kern-10pt $d(G)$}&=&\hbox to 0pt
{\kern-9pt $d(L)$}&+&
\hbox to 0pt{\kern-13pt $d(H)$}\\ \hline
1&GL(1,\R) & \{ e\} & Spin(1,1) & \R &&1&=&1&+&0 \\
2&Sp(1,\R) & \{ e\} & Spin(1,2) & \R^2&&2&=&2&+&0 \\
3&Sp(1,\C) & \{e \} & Spin(1,3) & \C^2&&3&=&3&+&0 \\
4&Sp(1,1) & Sp(1) & Spin(1,4) & \bbH^2&& 4&=&4&+&0\\
5&GL(2,\bbH) & GL(1,\bbH) & Spin(1,5) & \bbH^2 && 6&=&5&+&1\\
6&O^*(8) & O^*(6) & Spin(1,6) & \bbH^4&& 12&=&6&+&6\\
7&O(8,\C) & O(7,\C) &Spin(1,7) & \C^8&& 28&=&7&+&21\\
8&O(8,8) & O(7,8) &Spin(1,8) & \R^{16}&& 64&=&8&+&56\\ \hline
\end{array}
$$
\begin{center}\bf Table~\thesubsection\end{center}
\end{thm}
Here, $G(p,q)$ is a Lie group contained in the Clifford algebra
$C(p,q)$ associated to the quadratic form of
indefinite signature $(p,q)$ (see Definition~\ref{df:gpq}).
The general structural theorem of $G(p,q)$ will be
stated in Proposition~\ref{prop:auto_iso}.
In particular, $G(p,q)$ acts on the vector space $W$ as a
natural representation, and $H$ is defined as the isotropy subgroup
at $e_1$, where $\{e_i\}$ is the standard basis of $W$.
\begin{rem}
1) For $q\leq 4$, Theorem~\ref{thm:spin_ck_form}
is just a special case of Borel's theorem (see Theorem~\ref{thm:borel}),
because $H$
is compact then.
This is also the case if $q=5$ because $H$ is locally
isomorphic to $\R\times SO(3)$.
\par\medskip\noindent
2) The cases $6\leq q\leq 8$ are remarkable as both of $H$ and
$L$ are non-compact. Theorems \ref{thm:intro:complex_sphere}
and~\ref{thm:spherical78} correspond to the cases
$q=7$ and $q=8$ of Theorem~\ref{thm:spin_ck_form}, respectively.
Corollary~\ref{cor:admit_cck} (12) (up to a compact factor) corresponds
to the case $q=6$ of Theorem~\ref{thm:spin_ck_form}.
\par\medskip\noindent
3) In light of the case $q=7$, we can prove that the triple
$(O(8,\C), O(1,7), Spin(7,\C))$ satisfies the conditions~\eqref{thm:red:suf}(a)
and (b) so that the symmetric space $O(8,\C)/O(7,\C)$ admits
a compact Clifford-Klein form (see Theorem~\ref{thm:o8c}).
This can be justified by an outer automorphism of $\fo(8,\C)$.
This can be also proved by using \cite{koba94c}. However, we shall
give another proof based on Clifford algebras in Section~\ref{subsec:o8c}.
\par\medskip\noindent
4) For $q\geq 9$, the condition \eqref{thm:red:suf}(a) holds, but
the condition \eqref{thm:red:suf}(b) fails. 
It is likely that $G/H$ does not admit a compact Clifford-Klein form if
$q\geq 9$. For example, $G/H=GL(16,\R)/GL(15,\R)$ for $q=9$.
We note that 
Conjecture~\ref{conj:inv:suff} implies that $GL(16,\R)/GL(15,\R)$ does
not admit compact Clifford-Klein forms.
But even in this
case, this is still an unsolved problem (see [Ko02] for the list of 
various results in the literature about
the (non-)existence of compact Clifford-Klein forms of non-symmetric
spaces $GL(n,\F)/GL(m,\F)\ (n>m;\ \F=\R,\C$ or $\bbH$). 
\end{rem}
\begin{proof}[Outline of proof]
The condition \eqref{thm:red:suf}(b) is clearly satisfied from the
Table~\thesubsection.
We shall prove the condition \eqref{thm:red:suf}(a) in a more
general setting in Section~\ref{subsec:clifford_proper}
(see Proposition~\ref{prop:clifford_proper}), 
after a general theory of Clifford algebras.
The last statement of Theorem~\ref{thm:spin_ck_form} is a direct
consequence of Theorem~\ref{thm:red:suf}.
\end{proof}

\subsection{Clifford algebra associated to an
indefinite quadratic form}\label{subsec:Clifford}

This subsection elucidates the structure of the Clifford algebra
associated to the indefinite quadratic form
$Q(a)=a_1^2+\dots +a_p^2
-a_{p+1}^2-\dots-a_{p+q}^2$ of signature $(p,q)$.
In contrast to the case that $Q$ is negative (or positive) definite
(e.g. \cite{abs}), there are more periodicity theorems and
abundant structures, in particular, through the shift
$(p,q)\mapsto (p+1,q+1)$
and the flip $(p,q)\mapsto (q+1,p-1)$.

 We have tried to make the exposition
here self-contained as far as possible.

Let $\E=\R^{p,q}$ be the vector space $\R^{p+q}$ equipped with 
the quadratic form $Q(a)$ for $a=\sum_{j=1}^{p+q}a_je_j\in \R^{p+q}$,
where $\{e_j\}$
is the standard basis of $\R^{p+q}$. Let $T(\E)$ be the tensor
algebra over $\R$, and let $I(Q)$ be the two sided ideal generated by the
element $a\otimes a-Q(a)\cdot 1$ in $T(\E)$.
The quotient algebra $T(\E)/I(Q)$ is called the {\it Clifford algebra},
and will be denoted by $C(Q)$ or $C(p,q)$.
We define $i_Q : \E\rightarrow C(p,q)$ to be the canonical map by the
composition $\E\rightarrow T(\E)\rightarrow
C(p,q)$. Then $i_Q$ is injective, and we shall regard 
$\E$ as a
subspace of $C(p,q)$.

The Clifford algebra $C(p,q)$ is the
$\R$-algebra generated by the unit $1$ and the symbols $e_j$
subject to the relations:
\begin{equation}
e_j^2=1\ (1\leq j\leq p),\ e_j^2=-1\ (p+1\leq j\leq p+q),\ 
e_ie_j+e_je_i=0\ (i\neq j).\label{eq:clifford}
\end{equation}
We state the universality of the $\R$-algebra $C(p,q)$ 
in the form that will be frequently used later:
Let $\phi:\E\mapsto A$ be a linear map of $\E$ into an $\R$-algebra $A$
with unit 1 satisfying the following conditions \eqref{eq:universal}:
\begin{equation}
\begin{split}
\phi(e_j)^2=1\ (1\leq j\leq p),\quad \phi(e_j)^2=-1\ (p+1\leq j\leq p+q),\\
\phi(e_i)\phi(e_j)+\phi(e_j)\phi(e_i)=0\ (i\neq j).\qquad\qquad\qquad
\end{split}
\label{eq:universal}
\end{equation}
Then there exists a unique $\R$-algebra homomorphism
$\tilde{\phi}:C(p,q)\rightarrow A$,
such that $\tilde{\phi}\circ i_Q=\phi$.
$\tilde{\phi}$ automatically satisfies $\tilde{\phi}(1)=1$.
We refer to $\tilde{\phi}$ as the
{\it extension} of $\phi$, and shall use the same notation $\phi$ instead of
$\tilde{\phi}$. In particular, if $\phi$ is a linear isomorphism
$\phi:\E\rightarrow\E (\subset C(p,q))$ that preserves the
quadratic form $Q$, then $\phi$ is extended to an algebra isomorphism
$\phi:C(p,q)\simarrow C(p,q)$.

The Clifford algebra is naturally a $\Z_2$-graded algebra
$$C(p,q)=\codd(p,q)\oplus \ceven(p,q),$$
as is defined to be the image in $C(p,q)$ of
$T(\E)=\bigoplus_{j=0}^\infty T^{2j+1}(\E)+\bigoplus_{j=0}^\infty
T^{2j}(\E).$

Further, we define an anti-ring homomorphism 
\begin{equation}
C(p,q)\rightarrow C(p,q),\quad
a\mapsto \tr{a}\label{eq:tr}
\end{equation}
 as the extension of the identity map $E\rightarrow E$.
If $a$ is of the form $a=a_1\cdots a_s\in C(p,q)$ for some $a_i\in E$,
then $\tr{a} = a_s\cdots a_1$.

\begin{dfn}\label{df:gpq}
We define the subset $C(p,q)\mon$ and the groups $G(p,q)$
 and $Spin(p,q)$
in $C(p,q)$ as follows:
\begin{align*}
G(p,q)&:=\set{a\in \ceven(p,q)}{\tr{a}a=1},\\
C(p,q)\mon&:=\bigcup_{s=0}^{\infty}\set{\pm a_1\cdots a_s}{a_i\in E},\\
Spin(p,q)&:=G(p,q)\cap C(p,q)\mon.
\end{align*}
\end{dfn}
We pin down an obvious inclusive relation:
\newcommand{\rsub}[2]{\lower #1 \hbox{\rotatebox{#2}{$\subset$}}}
$$
\begin{array}{cccccc}
& \rsub{10pt}{20} & C(p,q)\mon & \rsub{7pt}{-20} &\\
Spin(p,q)&& && C(p,q)\\
& \rsub{-7pt}{-20} & G(p,q)& \rsub{-4pt}{20}&
\end{array}
$$

\begin{rem}
The signature $\pm$ in the definition of $C(p,q)\mon$ is not very
important. In fact, an easy computation shows 
$C(p,q)\mon=\bigcup_{s=0}^{\infty}
\set{a_1\cdots a_s}{a_i\in E}$ if $(p,q)\neq (0,0), (1,0)$.
\end{rem}

The quadratic form $Q$ on $E$ defines an
inner product $\langle\ ,\ \rangle$ of signature $(p,q)$ on $E$ by
$$
\langle a,b\rangle := \tfrac{1}{2}(Q(a+b)-Q(a)-Q(b))
= \tfrac{1}{2}(ab+ba)
$$
for $a,b\in E$. We note $\langle a,a\rangle=Q(a)$.

Take $v\in E$ such that $Q(v)\neq 0$ and $x\in Spin(p,q)$. We set
\begin{align}
\tau_v(a)&:=\frac{-1}{Q(v)}v a v\label{eq:df:tau}\\
\rho(x)(a)&:=xa\tr{x}
\end{align}
for $a\in E$. It turns out that both $\tau_v(a)$ and $\rho(x)(a)$ belong
to $E$. In fact, both $\tau_v$ and $\rho(x)$ belong to 
the indefinite orthogonal group $O(p,q)$
as is seen in the following:

\begin{prop}
\begin{enumerate}
\ritem[1)] $\tau_v(a) \in E$ for $a\in E$. More precisely,
$\tau_v$ is a reflection in $E$
with respect to the inner product  $\langle \cdot, \cdot\rangle$
given by
\begin{equation}
\tau_v(a)= a - 2\frac{\langle a,v\rangle}{Q(v)}v.
\label{la:pq-reflection}
\end{equation}
In particular, $\tau_v\in O(p,q)$ and $\det(\tau_v)=-1$.
\ritem[2)] $\rho(x)= \tau_{x_1} \cdots\tau_{x_{2n}}$
for $x=x_1\cdots x_{2n}\in Spin(p,q)$.
\ritem[3)] If $x\in Spin(p,q)$, then $\rho(x)(a)\in E$ for any $a\in E$
and $\rho(x)\in SO(p,q)$.
\ritem[4)] $\rho$ defines a surjective group homomorphism $Spin(p,q)\rightarrow
SO(p,q)$ with $\Ker\rho=\{\pm 1\}$.
\end{enumerate}
\end{prop}
\begin{proof}
Use $\langle a,v\rangle =\tfrac{1}{2}(av+va)$ for (1). 
(2) and (3) are obvious.
The proof of the statement (4) parallels to that of the well-known case $q=0$.
\end{proof}

Later, we shall use frequently the following lemma
for specific choice of $v$.
\begin{la}\label{la:tau}
 The reflection $\tau_v: E\rightarrow E$ extends to
the $\R$-algebra involutive isomorphism $\tau_v: C(p,q)\simarrow C(p,q)$.
\end{la}

\subsection{Clifford algebra}\label{subsec:clifford_algebra}
This subsection summarizes structural results of the
Clifford algebra $C(p,q)$ and the group $G(p,q)$
(see Definition~\ref{df:gpq}) associated to
the quadratic form $Q$ of arbitrary signature $(p,q)$.
The proof will be given
in Sections~\ref{subsec:clifford_iso}
and \ref{subsec:auto_iso}, respectively.

For $\F=\R, \C$ and $\bbH$, we denote by $M(n,\F)$ the $\R$-algebra
($\C$-algebra if $\F=\C$) of $n$ by $n$ matrices whose
entries are in the field $\F$.
\begin{prop}\label{prop:clifiso}
For $p,q\geq 0$, we define $\alpha\in \{0,1,2,3\}$ as in the table below and
set $n=2^{(p+q-\alpha)/2}$. Then, the Clifford algebra
$C(p,q)$ is isomorphic to one of the following algebras:
$$
\arraycolsep 20pt
\renewcommand{\arraystretch}{1.3}
\begin{array}{c|c||c}
C(p,q) &  \alpha & p-q-1\mod 8\\ \hline
M(n,\R)\oplus M(n,\R) & 1 & 0\\
M(n,\R) & 0 & \pm 1\\
M(n,\C) & 1 & \pm 2\\
M(n,\bbH) & 2 & \pm 3\\
M(n,\bbH) \oplus M(n,\bbH) & 3 & 4
\end{array}
$$
\end{prop}
We note that $\alpha$ depends only on $p-q\mod 8$. The Clifford algebra
$C(p,q)$ carries a $\C$-algebra structure if $p-q-1\equiv \pm 2\mod 8$
(see Remark~\ref{rem:complex}).

In an obvious sense, $C(p,q)$ is of period $(8,0)$ and $(0,8)$, and
of $(1,1)$. By Proposition~\ref{prop:clifiso}, we have the following
 table of $C(p,q)$ for $p,q\leq 8$.
\unitlength 1pt
\setbox0=\hbox{
\begin{picture}(15,12)
\put(-6,11){\line(3,-2){24}}
\put(9,5.5){$q$}
\put(-2,-1){$p$}
\end{picture}}

\tabcolsep=3pt
\begin{center}
\begin{tabud}{|c|c|c|c|c|c|c|c|c|c|}{1}{1.3}\hline
\box0
& $0$ & $1$ & $2$ & $3$ & $4$ & $5$ & $6$ & $7$ & $8$ \\ \hline $0$ &
$\R$ & $\C$ & $\bbH$ & $\bbH^2$ & $\bbH(2)$ & $\C(4)$ & $\R(8)$ &
$\R(8)^2$ & $\R(16)$ \\ \hline $1$ & $\R^2$ & $\R(2)$ & $\C(2)$ &
$\bbH(2)$ & $\bbH(2)^2$ & $\bbH(4)$ & $\C(8)$ & $\R(16)$ & $\R(16)^2$
\\ \hline $2$ & $\R(2)$ & $\R(2)^2$ & $\R(4)$ & $\C(4)$ & $\bbH(4)$ &
$\bbH(4)^2$ & $\bbH(8)$ & $\C(16)$ & $\R(32)$ \\ \hline $3$ & $\C(2)$
& $\R(4)$ & $\R(4)^2$ & $\R(8)$ & $\C(8)$ & $\bbH(8)$ & $\bbH(8)^2$ &
$\bbH(16)$ & $\C(32)$ \\ \hline $4$ & $\bbH(2)$ & $\C(4)$ & $\R(8)$ &
$\R(8)^2$ & $\R(16)$ & $\C(16)$ & $\bbH(16)$ & $\bbH(16)^2$ &
$\bbH(32)$ \\ \hline $5$ & $\bbH(2)^2$ & $\bbH(4)$ & $\C(8)$ &
$\R(16)$ & $\R(16)^2$ & $\R(32)$ & $\C(32)$ & $\bbH(32)$ &
$\bbH(32)^2$ \\ \hline $6$ & $\bbH(4)$ & $\bbH(4)^2$ & $\bbH(8)$ &
$\C(16)$ & $\R(32)$ & $\R(32)^2$ & $\R(64)$ & $\C(64)$ & $\bbH(64)$ \\
\hline $7$ & $\C(8)$ & $\bbH(8)$ & $\bbH(8)^2$ & $\bbH(16)$ & $\C(32)$
& $\R(64)$ & $\R(64)^2$ & $\R(128)$ & $\C(128)$ \\ \hline $8$ &
$\R(16)$ & $\C(16)$ & $\bbH(16)$ & $\bbH(16)^2$ & $\bbH(32)$ &
$\C(64)$ & $\R(128)$ & $\R(128)^2$ & $\R(256)$ \\ \hline
\end{tabud}
\end{center}
Here, we have abbreviated $M(n,\F)$ to $\F(n)$, and $M(n,\F)\oplus M(n,\F)$
to $\F(n)^2$ for $\mathbb{F}=\R,\C$ or $\bbH$.
\begin{rem}
In the case $p=0$ or $q=0$, Table~\thesubsection.1 is the
same with \cite[Table 1]{abs}, where their $C_k$ corresponds to our $C(0,k)$
and likewise $C_k'$ to $C(k,0)$.
\end{rem}

A key lemma to Proposition~\ref{prop:clifiso} is the following
(see Section~\ref{subsec:clifford_iso} for the proof):
\begin{la}\label{la:iso_family}
We have two families of isomorphisms of Clifford algebra $C(p,q)$:
\begin{gather}
C(p,q)\simeq C(q+r,p-r)\qquad\text{for $r\equiv 1\mod 4$,} 
\label{eq:iso1}\\
C(p,q)\simeq C(p+r,q-r)\qquad\text{for $r\equiv 0\mod 4$.}
\label{eq:iso2}
\end{gather}
\end{la}

Next, we determine the group structure of $G(p,q)$
(see Definition~\ref{df:gpq}):
\newcommand{\pmid}[1]{\vbox to 0pt{\hbox{$#1$}}}
\begin{prop}\label{prop:auto_iso}
For $p,q\geq 1$, the group $G(p,q)$ is isomorphic to one of
the following classical Lie groups.
Here, $n=2^{(p+q-\alpha)/2}$ and $\alpha\in\{3,4,5,6\}$
are given according to $p-q\mod 8$ as in the table.
$$
\begin{array}{c|c|c||c|c}
G(p,q) & \ceven(p,q) & \alpha & p-q \mod 8& p+q \mod 8\\ \hline
O(n,n)^2 &             &  &    & 0 \\
GL(2n,\R) & M(2n,\R)^2 & 4& 0 & \pm 2 \\ 
Sp(n,\R)^2 &           &  &   & 4\\ \hline
O(n,n)  & \pmid{M(2n,\R)} &\pmid{3} & \pmid{\pm 1} & \pm 1\\
Sp(n,\R) &        &   &       & \pm 3\\ \hline
O(2n,\C) &        &   &       & 0\\
U(n,n)  & M(2n,\C)& 4 & \pm 2 &\pm 2\\
Sp(n,\C) &       &    &       & 4\\ \hline
O^*(4n) & \pmid{M(2n,\bbH)} & \pmid{5} &\pmid{\pm 3} &\pm 1\\
Sp(n,n) &                &          &             &\pm 3\\ \hline
O^*(4n)^2 &              &          &             & 0 \\ 
GL(2n,\bbH) & M(2n,\bbH)^2 & 6      & 4 & \pm 2\\
Sp(n,n)^2 & & & & 4 
\end{array}
$$
\end{prop}
In the above table, $n=\tfrac{1}{2}$ if $(p,q)=(1,1)$.

As is observed in Proposition~\ref{prop:auto_iso},
the type of $G(p,q)$ depends only on $p-q\mod 8$ and
$p+q\mod 8$. The following table
exhibits the ``basic pattern'' of
the group $G(p,q)$.
\begin{gather*}
G(p,q)\quad\text{for $p\pm q \mod 8$}\\
\begin{array} {c|ccccc}
\setbox0=\hbox{
\begin{picture}(20,12)
\put(-9,10){\line(5,-2){35}}
\put(5,7){\footnotesize $p+q$}
\put(-7,-2){\footnotesize $p-q$}
\end{picture}}
\box0
& 0 & \pm 1 & \pm 2 & \pm 3 & 4\\ \hline
0 & O(n,n)^2 & & GL(2n,\R) & & Sp(n,\R)^2\\
\pm 1 & & O(n,n) & & Sp(n,\R) \\
\pm 2 & O(2n,\C) & & U(n,n) & & Sp(n,\C)\\
\pm 3 & & O^*(4n) & & Sp(n,n)\\
4 & O^*(4n)^2 & & GL(2n,\bbH) & & Sp(n,n)^2
\end{array}\\[4pt]
\text{\bf Table~\thesubsection.2}
\end{gather*}
For the convenience of the reader, we present
the following table of the group $G(p,q)$ for $p,q\leq 8$ by adding
the easier cases $p=0$ and $q=0$.
\tabcolsep=0pt
\begin{center}
\unitlength 1pt
\setbox0=\hbox{
\begin{picture}(13,10)
\put(-3.3,15){\line(3,-4){16}}
\put(7,5.5){$q$}
\put(-1,-1){$p$}
\end{picture}}
$G(p,q)\quad (p,q\leq 8)$\nopagebreak\\
\vbox{\kern10pt
\begin{tabud}{|c|c|c|c|c|c|c|c|c|c|}{3}{1.8}\hline
\box0
& $0$ & $1$ & $2$ & $3$ & $4$ & $5$ & $6$ & $7$ & $8$ \\ \hline $0$ &
$O(1)$ & $O(1)$ & $U(1)$ & $Sp(1)$ & $Sp(1)^2$ & $Sp(2)$ &
$U(4)$ & $O(8)$ & $O(8)^2$ \\ \hline $1$ &
$O(1)$ & $GL(1,\R)$ & $Sp(1,\R)$ & $Sp(1,\C)$ & $Sp(1,1)$ & $GL(2,\bbH)$
& $O^*(8)$ & $O(8,\C)$ & $O(8,8)$ \\ \hline $2$ &
$U(1)$ & $Sp(1,\R)$ & $Sp(1,\R)^2$ & $Sp(2,\R)$ & $U(2,2)$ & $O^*(8)$ &
$O^*(8)^2$ & $O^*(16)$ & $U(8,8)$\\ \hline $3$ &
$Sp(1)$ & $Sp(1,\C)$ & $Sp(2,\R)$ & $GL(4,\R)$ & $O(4,4)$ & $O(8,\C)$ &
$O^*(16)$ & $GL(8,\bbH)$ & $Sp(8,8)$ \\ \hline $4$ &
$Sp(1)^2$ & $Sp(1,1)$ & $U(2,2)$ & $O(4,4)$ & $O(4,4)^2$ & $O(8,8)$ &
$U(8,8)$ & $Sp(8,8)$ & $Sp(8,8)^2$ \\ \hline $5$ & $Sp(2)$ & $GL(2,\bbH)$
& $O^*(8)$ & $O(8,\C)$ & $O(8,8)$ & $GL(16,\R)$ & $Sp(16,\R)$ &
$Sp(16,\C)$ & $Sp(16,16)$ \\ \hline $6$ &
$U(4)$ & $O^*(8)$ & $O^*(8)^2$ & $O^*(16)$ & $U(8,8)$ & $Sp(16,\R)$ &
$Sp(16,\R)^2$ & $Sp(32,\R)$ & $U(32,32)$ \\ \hline $7$ &
$O(8)$ & $O(8,\C)$ & $O^*(16)$ & $GL(8,\bbH)$ & $Sp(8,8)$ & $Sp(16,\C)$ &
$Sp(32,\R)$ & $GL(64,\R)$ & $O(64,64)$ \\ \hline $8$ &
$O(8)^2$ & $O(8,8)$ & $U(8,8)$ & $Sp(8,8)$ & $Sp(8,8)^2$ & $Sp(16,16)$ &
$U(32,32)$ & $O(64,64)$ & $O(64,64)^2$\\ \hline
\end{tabud}}
\end{center}
\begin{rem}
As one may observe from Table~\thesubsection.3, the pair
($G(p+1,q), G(p,q)$) forms a
symmetric pair, and likewise $(G(p,q+1), G(p,q))$.
In fact, $\tau_1^+$ is the involutive automorphism of 
the group $G(p+1,q)$ such that the group of the fixed points
$G(p+1,q)^{\tau_1^+}=G(p,q)$, and likewise, $G(p,q+1)^{\tau_1^-}=G(p,q)$,
where $\tau_1^\pm$ will be defined in \eqref{eq:tau}.
\end{rem}

\newcommand{\bk}{{\mathcal K}}
\newcommand{\bl}{{\mathcal L}}
\newcommand{\bm}{{\mathcal M}}
\newcommand{\bx}{{\mathcal X}}
\newcommand{\by}{{\mathcal Y}}
\subsection{Definition and Lemmas}\label{subsec:spin_def_la}
For later argument, it will be convenient to use the following 
notation for the basis $\{e_j\}$ of $E$:
$$v^+_j:=e_j\ (1\leq j\leq p),\quad v^-_j:=e_{p+j}\ (1\leq j\leq q).$$
Furthermore, we write
$\tau_i^\pm$ for the involutive automorphism induced from
reflection with respect to the
normal direction $v_i^\pm$, namely,
\begin{equation}
\tau_i^\pm:= \tau_{v_i^\pm}.\label{eq:tau}
\end{equation}
Note that $\tau_i^\pm$ mutually commute because $e_i^\pm$ are orthogonal
to each other.

For pairs of natural numbers $\bk=(k_+,k_-),\ \bl=(l_+,l_-)$,
we define:
\begin{align*}
\Delta \bk &:=k_+ - k_-,\\
\bk\V &:= (k_-, k_+),\\
\bk\pm \bl &:= (k_+ \pm l_+, k_- \pm l_-).
\end{align*}

\begin{dfn}\label{df:vk}
Let $\bm=(p, q), \bk=(k_+, k_-)$ and $\bl=(l_+, l_-)$ be pairs of natural
numbers
such that $\bm=\bk+\bl$.
We define elements $V_\bk, J_\pm, J$ in $C(p,q)$, and
involutive automorphisms $T_\bl$ of $C(p,q)$ by:
\begin{align*}
V_\bk&:=v_1^+\cdots v_{k_+}^+v_1^-\cdots v_{k_-}^-,\\
J_+&:=V_{p,0},\\
J_-&:=V_{0,q},\\
J&:=V_{p,q}=J_+J_-,\\
T_\bl&:=\tau_{p-l_++1}^+\cdots\tau_{p}^+\tau_{q-l_-+1}^-
\cdots\tau_{q}^-.
\end{align*}
\end{dfn}
We also use the notation $m_+:=p$ and $m_-:=q$ so that $\bm=(m_+, m_-)$.
In particular, $\Delta\bm=p-q$.
Then, we have from definition:
\begin{equation}
T_\bl(v_i^\pm)=
\begin{cases}
v_i^\pm & \IF 1\leq i\leq m_\pm - l_\pm,\\
-v_i^\pm & \IF m_\pm-l_\pm< i \leq m_\pm.
\end{cases}\label{eq:tl}
\end{equation}

A straightforward computation shows
\begin{equation}
V_\bk^2=
\begin{cases}
1 & \IF \Delta \bk\equiv 0,1 \mod 4,\\
-1 & \IF \Delta \bk \equiv 2,3 \mod 4.
\end{cases}\label{eq:vk2}
\end{equation}
In particular, we have:
\begin{equation}
J^2=
\begin{cases}
1 & \IF\Delta \bm\equiv 0,1\mod 4,\\
-1 & \IF\Delta \bm\equiv 2,3\mod4.
\end{cases}\label{eq:j2}
\end{equation}
\begin{rem} \label{rem:complex}
We have $J a= a J$ for any $a\in C(p,q)$ if $\Delta \bm\equiv 1\mod 2$.
It then follows from \eqref{eq:j2} that
$J$ becomes a central element satisfying $J^2=-1$ if
and only if $\Delta \bm\equiv 3 \mod 4$. In this case,
$C(p,q)$ is endowed with the $\C$-algebra structure with
complex structure $J$. We recall Proposition~\ref{prop:clifiso}
asserts particularly that $C(p,q)$ is isomorphic to $M(n,\C)$
 if $m_+-m_--1\equiv \pm 2\mod 8$, namely, if $\Delta \bm\equiv 3\mod 4$.
\end{rem}

Next, given an automorphism $\tau$ of $C(p,q)$, we denote by
$C(p,q)^{\tau}$
the subalgebra consisting of $\tau$-fixed elements:
$$ C(p,q)^\tau:=\set{a\in C(p,q)}{\tau(a)=a}.$$

\begin{ex}\label{ex:tau} Retain the notation as in \eqref{eq:tau}
and Definition~\ref{df:vk}.
\begin{enumerate}
\ritem[1)] $C(p,q)^{T_{p,q}} = \ceven(p,q)$.
\ritem[2)] $C(p,q)^{\tau_i^+} \simeq C(p-1,q) \qquad (1\leq i\leq p)$.
\ritem[3)] $C(p,q)^{\tau_i^-} \simeq C(p,q-1) \qquad (1\leq i\leq q)$.
\end{enumerate}
\end{ex}

The composition of $a\mapsto \tr{a}$ (see \eqref{eq:tr}) and $T_{0,q}$
gives rise to an anti-automorphism, denoted by $a\mapsto \st{a}$,
of the $\R$-algebra $C(p,q)$:
\begin{align}
\st{a}&:=T_{0,q}(\tr{a})&&(a\in C(p,q)). \label{eq:st}
\intertext{We shall see in Lemma~\ref{la:gpq_aut} that}
\st{a}&=T_{p,0}(\tr{a})&&(a\in \ceven(p,q)). \label{eq:st2}
\end{align}
All of the isomorphisms, say $\phi$, that we shall need
later satisfy
 $\phi(\st{a})=\st{\phi(a)}$ (see
Lemma~\ref{la:st_compat}), but
$\phi(\tr{a})$ is not always
equal to $\tr{\phi(a)}$. Therefore, in order to
find the structure of the group $G(p,q)$, it is more convenient to
rewrite the definition of $G(p,q)$ (see Definition~\ref{df:gpq})
by means of $\st{a}$ instead of $\tr{a}$.
For this purpose, we prepare the following definition and
lemma:
\begin{dfn}\label{df:aut} Retain the notation as in
Definition~\ref{df:vk} and \eqref{eq:st}.
We define the groups $\Aut(J_+)$ and $\Aut(J_-)$ as
$$\Aut(J_\pm):=\set{a\in C(p,q)}{\st{a}J_\pm a=J_\pm}.$$
\end{dfn}
Then, we have:
\begin{la}\label{la:gpq_aut}
\begin{enumerate}
\ritem[1)] $J_+^{-1}aJ_+=J_-^{-1}aJ_-=T_{p,0}(a)=T_{0,q}(a) \qquad 
\text{for any $a\in \ceven(p,q)$}$.
\ritem[2)] $G(p,q)=\ceven(p,q)\cap\Aut(J_+)=\ceven(p,q)\cap\Aut(J_-)$.
\end{enumerate}
\end{la}
\begin{proof}
\noindent
1)
By the definition of $J_+$ and $J_-$, we have:
\begin{align*}
J_+^{-1}v_i^+J_+&=(-1)^{p-1}v_i^+,\qquad
&J_+^{-1}v_i^-J_+&=(-1)^{p}v_i^-,\\
J_-^{-1}v_i^+J_-&=(-1)^{q}v_i^+,\qquad
&J_-^{-1}v_i^-J_-&=(-1)^{q-1}v_i^-.
\end{align*}
In light of \eqref{eq:tl}, we
obtain the following equation for any $a\in C(p,q)$:
$$
J_+^{-1}aJ_+=
\begin{cases}
T_{p,0}(a)&\IF p \text{ is even},\\
T_{0,q}(a)&\IF p \text{ is odd},
\end{cases}\qquad
J_-^{-1}aJ_-=
\begin{cases}
T_{0,q}(a)&\IF q \text{ is even},\\
T_{p,0}(a)&\IF q \text{ is odd}.
\end{cases}
$$
On the other hand, because $T_{p,q}(a)=a$ for any $a\in\ceven(p,q)$
by Example~\ref{ex:tau}(1), we obtain:
$$ T_{p,0}(a)=T_{p,0}\circ T_{p,q}(a)=T_{0,q}(a).$$
Thus, the statement (1) has been proved.
\par\medskip\noindent
2) 
Note that $\st{a}\in\ceven(p,q)$ if $a\in\ceven(p,q)$.
Then we have:
$$ J_+^{-1}\st{a}J_+=J_-^{-1}\st{a}J_-= T_{0,q}(\st{a})=\tr{a}.$$
Here, the last equality follows from \eqref{eq:st}.
Thus, (2) follows. 

\end{proof}
\begin{rem}\label{rem:pq_odd}
If $p-q\equiv 1\mod 2$, then 
the above proof shows that $J_+^{-1} a J_+ = J_-^{-1} a J_-$.
Therefore, we have
$C(p,q)\cap\Aut(J_+)=C(p,q)\cap\Aut(J_-)$.
\end{rem}
\subsection{Proof of Proposition~\ref{prop:clifiso}}\label{subsec:clifford_iso}
This subsection establishes various isomorphisms among Clifford
algebras (Lemma~\ref{la:phi}) and completes the proof of the
structural result (Lemma~\ref{la:iso_family} and
Proposition~\ref{prop:clifiso}) of the Clifford algebra
$C(p,q)$.
First, for $(p,q)$ such that $p+q\leq 2$,
Proposition~\ref{prop:clifiso} is directly verified
by the following lemma:
\begin{la}\label{la:clifford_hom2}
We have the following isomorphisms as $\R$-algebras:
\begin{align*}
\psi_{0,1}:C(0,1)\simarrow\ &\C,\\
\psi_{1,0}:C(1,0)\simarrow\ &\R\oplus\R,\\
\psi_{0,2}:C(0,2)\simarrow\ &\bbH,\\
\psi_{1,1}:C(1,1)\simarrow\ &M(2,\R),\\
\psi_{2,0}:C(2,0)\simarrow\ &M(2,\R).
\end{align*}
\end{la}
\begin{proof}
The linear map $\psi_{1,1}: \R^{1,1}\rightarrow
M(2,\R)$ given by
\begin{equation}
v_1^+\mapsto
\begin{pmatrix}
1 &  \\
 & -1
\end{pmatrix},\quad
v_1^-\mapsto
\begin{pmatrix}
& -1 \\
1 &
\end{pmatrix}\label{eq:clif11}
\end{equation}
induces the algebra isomorphism $\psi_{1,1}: C(1,1)\simarrow M(2,\R)$
by the universality of the $\R$-algebra $C(1,1)$.
Other isomorphisms are well-known, but
we write down these maps explicitly on the generators
for later purpose.
\begin{align*}
\psi_{0,1}&: C(0,1)\simarrow \C, &v_1^- &\mapsto i,\\
\psi_{1,0}&: C(1,0)\simarrow \R\oplus\R, &v_1^+ &\mapsto (1,-1),\\
\psi_{0,2}&: C(0,2)\simarrow \bbH, &v_1^-&\mapsto i, &v_2^-&\mapsto j,\\
\psi_{2,0}&: C(2,0)\simarrow M(2,\R), &v_1^+&\mapsto
\begin{pmatrix}
1& \\
 & -1
\end{pmatrix}, &v_2^+&\mapsto
\begin{pmatrix}
 & 1\\
1 &
\end{pmatrix}.
\end{align*}
\end{proof}
For general $(p,q)$ such that $p+q\geq 3$, one shall find
the structure of $C(p,q)$
inductively based on the following:
\begin{proof}[Proof of Proposition~\ref{prop:clifiso}]
\eqref{eq:iso2} reduces the proof of Proposition~\ref{prop:clifiso}
to the case $|p-q-1|\leq 4$. Then \eqref{eq:iso1} with $r=1$
reduces it to the case $-4\leq p-q-1\leq 0$.
\par\noindent
Case 1) $p-q-1=0$.\\
An iterative use of $\phi_{1,1}$ yields:
$$ C(p,q)=C(q+1,q)\simeq \left(\bigotimes^q C(1,1)\right)\otimes C(1,0)
\simeq M(2^q,\R)\otimes (\R\oplus \R).$$
Here, the last isomorphism is induced by $\psi_{1,1}: C(1,1)\simarrow M(2,\R)$
and $\psi_{1,0}: C(1,0)\simarrow \R\oplus\R$.
\par\noindent
Case 2) $p-q-1=-1$, $-2$, $-3$, or~$-4$.\\
We write $q=p+s$ ($s=0$, $1$, $2$, or~$3$). Then we have
$$ C(p,q)=C(p,p+s)\simeq \left(\bigotimes^p C(1,1)\right)\otimes C(0,s)
\simeq M(2^p,\R)\otimes C(0,s). $$
Thus, Proposition~\ref{prop:clifiso} is reduced to the isomorphism
$C(0,0)\simeq\R$, $C(0,1)\simeq\C$, $C(0,2)\simeq\bbH$ (see
Lemma~\ref{la:clifford_hom2}) or $C(0,3)\simeq \bbH\oplus\bbH$,
where the last isomorphism is derived from the composition:
$$C(0,3)\maparrow{\phi_{0,2}} C(0,2)\otimes C(1,0)
\maparrow{\psi_{0,2}\otimes \psi_{1,0}} \bbH\otimes (\R\oplus\R)
\simeq \bbH\oplus\bbH.$$
Hence, we have proved Proposition~\ref{prop:clifiso}.
\end{proof}

\begin{la}\label{la:phi}
For $\bm=(m_+,m_-)\in\N^2$, we write simply $C(\bm)$ for the
Clifford algebra $C(m_+, m_-)$.
\begin{enumerate}
\ritem[1)]
Suppose $\bk,\bl\in \N^2$ such that $\bm=\bk+\bl$.
Then, we have the following $\R$-algebra isomorphism $\phi_\bk$
as an extension of the linear map below:
\begin{equation*}
\phi_\bk: C(\bk+\bl)\rightarrow
\begin{cases}
C(\bk)\otimes C(\bl) & \IF \Delta \bk\equiv 0 \mod 4,\\
C(\bk+\bl\V) & \IF \Delta \bk\equiv 1 \mod 4,\\
C(\bk)\otimes C(\bl\V) & \IF \Delta \bk\equiv 2 \mod 4,\\
C(\bk+\bl) & \IF \Delta \bk\equiv 3 \mod 4.
\end{cases}
\end{equation*}
\begin{enumerate}
\ritem[Case i)]  $\Delta \bk\equiv 0 \mod 4$. $\phi_\bk:\E\rightarrow
C(\bk)\otimes C(\bl)$ is given by
\begin{equation*}
v_i^\pm \mapsto v_i^\pm\otimes 1\quad (1\leq i \leq k_\pm),\qquad
v_{k_\pm+\,i}^\pm \mapsto V_\bk\otimes v_i^\pm
\quad (1\leq i \leq l_\pm).
\end{equation*}
\ritem[Case ii)]  $\Delta \bk\equiv 1 \mod 4$. $\phi_\bk:\E\rightarrow
C(\bk+\bl\V)$ is given by
\begin{equation*}
v_i^\pm \mapsto v_i^\pm\quad (1\leq i \leq k_\pm),\qquad
v_{k_\pm+\,i}^\pm \mapsto V_\bk v_{k_\mp+i}^\mp
\quad (1\leq i \leq l_\pm).
\end{equation*}
\ritem[Case iii)]  $\Delta \bk\equiv 2 \mod 4$. $\phi_\bk:\E\rightarrow
C(\bk)\otimes C(\bl\V)$ is given by
\begin{equation*}
v_i^\pm \mapsto v_i^\pm\otimes 1\quad (1\leq i \leq k_\pm),\qquad
v_{k_\pm+\,i}^\pm \mapsto V_\bk\otimes v_i^\mp\quad (1\leq i \leq l_\pm).
\end{equation*}
\ritem[Case iv)] $\Delta \bk\equiv 3 \mod 4$. $\phi_\bk:\E\rightarrow
C(\bk+\bl)$ is given by
\begin{equation*}
v_i^\pm \mapsto v_i^\pm\quad (1\leq i \leq k_\pm),\qquad
v_{k_\pm+\,i}^\pm \mapsto V_\bk v_{k_\pm+i}^\pm\quad (1\leq i \leq l_\pm).
\end{equation*}
\end{enumerate}
\ritem[2)] 
Furthermore, $\phi_\bk$ transfers the automorphism
$\tau_i^{\pm}$ of the Clifford algebra $C(p,q)$ to the following
automorphism:
\medskip
\halign{ \hfil\qquad {\rm (#)}&\ $\tau_i^\pm \mapsto #\quad (1\leq i\leq k_\pm),
\hfil$\quad
& $\tau_{k_\pm+\,i}^\pm\mapsto #$\hfil\cr
i&\tau_i^\pm\otimes T_{\bl}&\id\otimes \tau_i^\pm\quad (1\leq i\leq l_\pm),\cr
ii&\tau_i^\pm T_{\bl\V}&\tau_{\mp+i}^\mp\quad (1\leq i\leq l_\pm),\cr
iii&\tau_i^\pm\otimes T_{\bl\V}&
\id\otimes \tau_i^\mp\quad (1\leq i\leq l_\pm),\cr
iv&\tau_i^\pm T_{\bl}&\tau_{\pm+i}^\pm\quad (1\leq i\leq l_\pm).\cr
}

\end{enumerate}

\end{la}
\begin{proof} 1)
Let us give a proof in Case (i). Other cases are similar.
To see $\phi_\bk$ extends to an $\R$-algebra homomorphism,
it is enough to verify the condition \eqref{eq:universal}.
This is so because $V_\bk^2=1$ by \eqref{eq:vk2} and
$V_\bk v_i^\pm=-v_i^\pm V_\bk\  (1\leq i \leq k_\pm)$. Thus,
the linear map $\phi_\bk:E\rightarrow C(p,q)$ extends to
an $\R$-algebra homomorphism $C(p,q)\rightarrow C(p,q)$.

Next, to see $\phi_\bk$ is surjective,
it is sufficient to check that the image of $\phi_\bk$ contains
generators of the algebra $C(\bk)\otimes C(\bl)$. This is so because
$$ \phi_\bk(v_i^\pm)=v_i^\pm\otimes 1\ (1\leq i\leq k_\pm),\qquad
\phi_\bk(V_\bk v_{k_\pm+i}^\pm)=1\otimes v_i^\pm\ (1\leq i\leq l_\pm).$$
Finally, the injectivity of $\phi_\bk$ follows from
$$\dim C(p,q)=2^{p+q}=2^{k_++k_-} 2^{l_++l_-}=\dim C(\bk)\otimes C(\bl).$$
Thus, the statement (1) is proved.

2)
Let us verify that $\phi_\bk$ transfers $\tau_1^+$ to $\tau_1^+\otimes
T_\bl$, or equivalently, $ \phi_\bk\circ \tau_1^+ =
 (\tau_1^+\otimes T_\bl)\circ\phi_\bk$ in Case (i) ($k_+\geq 1$).
Other cases are similar.
This equation can be checked by the evaluation at the generators
($v_1^+,\dots,v_p^+,v_1^-,\dots,v_q^-$) as follows:
\begin{align*}
\phi_\bk\circ\tau_1^+(v_1^+)&=-v_1^+&&=(\tau_1^+\otimes T_\bl)\circ \phi_\bk(v_1^+),\\
\phi_\bk\circ\tau_1^+(v_i^+)&=v_i^+&&=(\tau_1^+\otimes T_\bl)\circ\phi_\bk(
v_i^+) && (2\leq i\leq k_+),\\
\phi_\bk\circ\tau_1^+(v_i^-)&=v_i^-&&=(\tau_1^+\otimes T_\bl)\circ\phi_\bk(
v_i^-) && (1\leq i\leq k_-),\\
\phi_\bk\circ\tau_1^+(v_{k_\pm+i}^\pm)&=V_\bk\otimes v_i^\pm
&&=(\tau_1^+\otimes T_\bl)\circ\phi_\bk(v_{k_\pm+i}^\pm) &&
(1\leq i\leq k_\pm).
\end{align*}
\end{proof}

We pin down $\phi_\bk$ explicitly in the cases
$\bk=(0,1), (1,0), (0,2), (1,1)$ and $(2,0)$, $\bl=(p,q)$.
\makeatletter
\def\arraytag{\@eqnswtrue\make@display@tag}
\makeatother
$$
\begin{array}{rccclr}
\phi_{0,1}:& C(p,q+1) &\simarrow& C(p,q+1)&&\kern42.2pt\arraytag\\[5pt]
& v_i^+ &\mapsto & v_1^-v_i^+ & (1\leq i\leq p)\\
& v_1^- &\mapsto& v_1^- \\
& v_{i+1}^- &\mapsto & v_1^-v_i^- & (1\leq i\leq q),\\[10pt]
\phi_{1,0}:& C(p+1,q) &\simarrow& C(q+1,p)&&\arraytag\\[5pt]
& v_1^+ &\mapsto & v_1^+\\
& v_{i+1}^+ &\mapsto & v_1^+v_i^- & (1\leq i \leq p)\\
& v_i^- &\mapsto & v_1^+v_{i+1}^+ & (1\leq i \leq q),\\[10pt]
\phi_{0,2}: & C(p,q+2) &\simarrow& C(0,2)\otimes C(q,p)&&\arraytag\\[5pt]
& v_i^+&\mapsto & v_1^- v_2^-\otimes v_i^-& (1\leq i\leq p)\\
& v_i^-&\mapsto & v_i^-\otimes 1 & (i=1, 2)\\
& v_{i+2}^-&\mapsto & v_1^- v_2^-\otimes v_i^+ &(1\leq i\leq q),\\[10pt]
\phi_{1,1}: & C(p+1,q+1) &\simarrow& C(1,1)\otimes C(p,q)&&\arraytag\\[5pt]
& v_1^\pm &\mapsto & v_1^\pm\otimes 1\\
& v_{i+1}^+&\mapsto & v_1^+v_1^-\otimes v_i^+ & (1\leq i\leq p)\\
& v_{i+1}^-&\mapsto & v_1^+v_1^-\otimes v_i^- & (1\leq i\leq q),\\[10pt]
\phi_{2,0}: & C(p+2,q) &\simarrow& C(2,0)\otimes C(q,p)&&\arraytag\\[5pt]
& v_i^+&\mapsto & v_i^+\otimes 1 & (i=1, 2)\\
& v_{i+2}^+&\mapsto & v_1^+ v_2^+\otimes v_i^-& (1\leq i\leq p)\\
& v_i^-&\mapsto & v_1^+ v_2^+\otimes v_i^+ &(1\leq i\leq q).
\end{array}
$$
\begin{rem}
Among the above isomorphisms, positive (or negative) definite
cases such as the isomorphisms $\phi_{0,2}: C(0,q+2)\simarrow C(0,2)\otimes
C(q,0)$ and $\phi_{2,0}: C(p+2,0)\simarrow C(2,0)\otimes C(0,p)$
are found in \cite{abs}.
\end{rem}
\subsection{Isomorphism $\lambda$}

In this subsection, we introduce an $\R$-algebra isomorphism $\lambda_\bk$
(see Lemma~\ref{la:df:lambda}). 
The effect of $\lambda_\bk$, as we shall show in
Lemmas \ref{la:eta} and~\ref{la:xi}, is to make the `degree' of
$J_\pm$, $T_{p,0}$ or $T_{0,q}$ down to one.
These two lemmas will play a key role in finding the
group structure of $G(p,q)$ in Section~\ref{subsec:auto_iso}.

Suppose that $a\in C(p,q)$ is of the form
\begin{equation}
a=a_1\cdots a_s, \label{eq:monomial}
\end{equation}
where $a_1,\dots,a_s$ are distinct elements among the set of the
generators $v_1^+,\cdots,v_p^+,v_1^-,\cdots,v_q^-$.
Then
$$\{a_1,\dots,a_s\}=\{v_{i_1}^+,\dots,v_{i_t}^+,
v_{j_1}^-,\dots,v_{j_r}^-\}\qquad (s=t+r)$$
is independent of the expression $a=a_1\cdots a_s$. We
set
\begin{align*}
\ind(a)&:= t-r,\\
\deg(a)&:= s\quad (=t+r).
\end{align*}
Suppose another element $b\in C(p,q)$ is also expressed as a product
$b=b_1\cdots b_{s'}$ of the generators.
We denote by $\deg(a,b)$ the cardinality 
of the generators occurred in
both $a$ and $b$. For example:
\begin{align*}
\ind(v_1^+v_2^-v_2^+)&=\#\{v_1^+,v_2^+\}-\#\{v_2^-\}=1,\\
\deg(v_1^+v_2^-v_2^+)&=\#\{v_1^+,v_2^-,v_2^+\}=3,\\
\deg(v_1^+v_2^-v_2^+, v_1^- v_2^+ v_1^+)&=
\#\{v_1^+, v_2^+\}=2.
\end{align*}
With this notation, we have:
\begin{la} \label{la:sign}
For $a$ and $b$ $\in C(p,q)$ of the form
\eqref{eq:monomial},
we have:
\begin{enumerate}
\ritem[1)] $a^2=(-1)^{\tfrac{1}{2} \ind(a)(\ind(a)-1)}=
\begin{cases}
1 & \IF \ind(a)\equiv 0,1 \mod 4,\\
-1 &\IF \ind(a)\equiv 2,3 \mod 4.
\end{cases}$
\ritem[2)]
$ab=(-1)^{\deg(a)\cdot\deg(b)-\deg(a, b)} ba.$
\end{enumerate}
\end{la}
\begin{proof} Follows from straightforward computation.
\end{proof}

\begin{la}\label{la:df:lambda}
Let $\bm=(m_+,m_-):=(p,q)$ and $\bk=(k_+,k_-)$ satisfy
$1\leq k_\pm\le m_\pm$. We set $V:=V_{k_+-1,k_--1}\ 
(=v_1^+\cdots v_{k_+-1}^+v_1^-\cdots v_{k_--1}^-)$.
Then the linear map $\lambda_\bk: E\rightarrow C(p,q)$ given below
extends to an $\R$-algebra automorphism $\lambda_\bk: C(p,q)\simarrow C(p,q)$.
\begin{enumerate}
\ritem[Case 1)] $\Delta \bk\equiv 0 \mod 2$.
\begin{align*}
v_i^\pm&\mapsto v_{k_+}^+v_{k_-}^-v_i^\pm && (1\leq i\leq k_\pm-1),\\
v_{k_\pm}^\pm &\mapsto  
\begin{cases}
V  v_{k_\pm}^\pm & \IF\Delta \bk\equiv 0 \mod 4,\\
V  v_{k_\mp}^\mp & \IF \Delta \bk\equiv 2 \mod 4,
\end{cases}\\
v_i^\pm&\mapsto v_{i}^\pm&& (k_\pm+1\leq i\leq m_\pm).
\end{align*}
\ritem[Case 2)] $\Delta \bk\equiv 1 \mod 2$.
\begin{align*}
v_i^\pm&\mapsto v_i^\pm&& (1\leq i\leq k_\pm-1),\\
v_{k_\pm}^\pm &\mapsto
\begin{cases}
V v_{k_\mp}^\mp & \IF \Delta \bk\equiv 1 \mod 4,\\
V v_{k_\pm}^\pm & \IF \Delta \bk\equiv 3 \mod 4,
\end{cases}\\
v_i^\pm&\mapsto v_{k_+}^+v_{k_-}^-v_i^\pm&& (k_\pm+1\leq i\leq m_\pm).
\end{align*}
\end{enumerate}
\end{la}
\begin{proof}
It is easy to verify the condition \eqref{eq:universal}.
Then, $\lambda_\bk$ extends
to an $\R$-algebra automorphism $C(p,q)\rightarrow C(p,q)$
by the universality.
Furthermore,
a simple computation on the generators shows:
$$
(\lambda_\bk)^2=
\begin{cases}
\id & \IF \Delta \bk\equiv 0,1 \mod 4,\\
\tau_{k_+}^+\tau_{k_-}^- & \IF \Delta \bk\equiv 2,3 \mod 4.
\end{cases}
$$
Hence, $\lambda_\bk$ is bijective.
\end{proof}

\begin{la}\label{la:tau:lambda}
 Retain the setting as in Lemma~\ref{la:df:lambda}.
Then, the $\R$-algebra isomorphism
$\lambda_\bk$ transfers the automorphism $\tau_i^\pm$
as follows:

\begin{enumerate}
\ritem[Case 1)] $\Delta \bk\equiv 0 \mod 2$.
\begin{align*}
\tau_i^\pm&\mapsto \tau_i^\pm && (1\leq i\leq k_\pm-1),\\
\tau_{k_\pm}^\pm &\mapsto  
\begin{cases}
\tau_{k_\pm}^\pm T_\bm T_{\bm-\bk} & \IF\Delta \bk\equiv 0 \mod 4,\\
\tau_{k_\mp}^\mp T_\bm T_{\bm-\bk} & \IF \Delta \bk\equiv 2 \mod 4,
\end{cases}\\
\tau_i^\pm&\mapsto \tau_{i}^\pm\tau_{k_+}^+\tau_{k_-}^- && (k_\pm+1\leq i\leq m_\pm).
\end{align*}
\ritem[Case 2)] $\Delta \bk\equiv 1 \mod 2$.
\begin{align*}
\tau_i^\pm&\mapsto \tau_i^\pm\tau_{k_+}^+\tau_{k_-}^- && (1\leq i\leq k_\pm-1),\\
\tau_{k_\pm}^\pm &\mapsto  
\begin{cases}
\tau_{k_\pm}^\pm T_{\bm-\bk} & \IF\Delta \bk\equiv 1 \mod 4,\\
\tau_{k_\mp}^\mp T_{\bm-\bk} & \IF \Delta \bk\equiv 3 \mod 4,
\end{cases}\\
\tau_i^\pm&\mapsto \tau_{i}^\pm && (k_\pm+1\leq i\leq m_\pm).
\end{align*}

\end{enumerate}
\end{la}
\begin{proof}
The proof goes similarly as in Lemma~\ref{la:phi}(2).
\end{proof}
\begin{rem}
The isomorphism $\lambda_\bk$ is essentially
the composition $\phi_{1,1}^{-1}\circ\phi_{\bk'}^{-1}\circ\phi_{1,1}\circ
\phi_{\bk'}$ up to an automorphism induced by switching of
generators, where $\bk'=(k_+-1,k_--1)$.

For example, in the case $\Delta \bk\equiv 0\mod 4$, 
$\lambda_\bk$ may be defined as the composition of the
following $\R$-algebra isomorphisms:
\begin{gather*}
C(p,q)=C(\bk'+\bl')\simeq C(\bk')\otimes C(\bl')\simeq C(\bk')\otimes
C(1,1)\otimes
C(l'_+-1,l'_--1)\\
\simeq C(1,1)\otimes C(\bk')\otimes C(l'_+-1,l'_--1)
\simeq C(1,1)\otimes C(p-1,q-1)\simeq C(p,q).
\end{gather*}
In Lemma~\ref{la:df:lambda}, we have adopted an alternative
and direct definition of $\lambda_\bk$.
\end{rem}

As we shall see in the following lemma,
the point of introducing the isomorphism $\lambda_\bk$
is that it sends certain
elements of ``higher degree'' such as $V_*$ or $T_*$ into
those of degree one.

\begin{la}\label{la:lambda}
Suppose that the pairs $\bm=(m_+,m_-)$ and $\bk=(k_+,k_-)$ satisfy
$1\leq k_\pm\leq m_\pm$. Then we have:
\begin{enumerate}
\ritem[1)]$\lambda_\bk(V_{k_+-1,k_--1} v_{k_\pm}^\pm)=
\begin{cases}
\varepsilon v_{k_\pm}^\pm & \IF \Delta \bk\equiv 0,3 \mod 4,\\
\varepsilon v_{k_\mp}^\mp & \IF \Delta \bk\equiv 1,2 \mod 4.
\end{cases}$
\ritem[2)]$\lambda_\bk(\tau_{k_\pm}^\pm T_{\bm-\bk})=
\begin{cases}
\tau_{k_\mp}^\mp & \IF \Delta \bk\equiv 1 \mod 4,\\
\tau_{k_\pm}^\pm & \IF \Delta \bk\equiv 3 \mod 4.
\end{cases}$
\end{enumerate}
Here, $\varepsilon\in\{1,-1\}$.
\end{la}
\begin{proof}
1) By a simple computation, one finds
$\lambda_\bk(V_{k_+-1,k_--1})=V_{k_+-1,k_--1}$.
Then, the statement (1) follows readily from the definition of
$\lambda_\bk(v_{k_\pm}^\pm)$ (see Lemma~\ref{la:df:lambda}).
\par\medskip\noindent
2) Similarly, the statement (2) is obtained by combining
the formula $\lambda_\bk(T_{\bm-\bk})=T_{\bm-\bk}$ if $\Delta \bk\equiv 1\mod 2$
and Lemma~\ref{la:tau:lambda}(1).
\end{proof}
\begin{rem}
It is not hard to specify $\varepsilon=\pm 1$
in Lemma~\ref{la:lambda} (1), but we do not need this information
in this article.
\end{rem}

We have already seen in Lemma~\ref{la:phi} and its examples
that $\phi_{1,1}$ gives an isomorphism of $\R$-algebras between $C(p,q)$
and $C(1,1)\otimes C(p-1,q-1)$. More than this, we need an isomorphism,
denoted by $\eta_\bk$,
having an additional property \eqref{eq:eta}.
\begin{la}\label{la:eta}
Suppose that $\bm=(p,q)$ and $\bk=(k_+,k_-)$ are
pairs of natural numbers satisfying
$k_+\leq p$, $k_-\leq q$, $1\leq p,q$, $\bk\neq (0,0)$, and $\bk\neq \bm$.
Then, there exists an $\R$-algebra
isomorphism $\eta_\bk: C(p,q)\simarrow C(1,1)\otimes C(p-1,q-1)$ satisfying:
\begin{equation}
\eta_\bk(V_\bk)=
\begin{cases}
v_1^+\otimes 1 & \IF \Delta \bk\equiv 0,1\mod 4,\\
v_1^-\otimes 1 & \IF \Delta \bk\equiv 2,3\mod 4.
\end{cases}\label{eq:eta}
\end{equation}

\end{la}
\begin{proof}
Assume first $k_+< m_+$ and $k_->0$. We set
$\bk':=(k_++1,k_-)$. (If this assumption is not satisfied, then
it is easy to see $k_+>0$ and $k_-<m_-$. The proof goes similarly
if we set $\bk':=(k_+,k_-+1)$ in this case.)

Clearly, we have $k'_{\pm}\geq 1$ and 
$V_{k'_+-1,k'_--1} v_{k_-}^-=V_\bk$, and therefore it follows from
Lemma~\ref{la:lambda} that

\begin{equation*}
\lambda_{\bk'}(V_\bk)=
\begin{cases}
\varepsilon v_{{k'}_+}^+ &\IF \Delta \bk\equiv 0,1 \mod 4,\\
\varepsilon v_{{k'}_-}^- &\IF \Delta \bk\equiv 2,3 \mod 4.
\end{cases}
\end{equation*}
\noindent
We set
\begin{align*}
v&:=
\begin{cases}
\varepsilon v_{{k'}_+}^+ - v_1^+ &\IF \Delta \bk\equiv 0,1\mod 4,\\
\varepsilon v_{{k'}_-}^- - v_1^- &\IF \Delta \bk\equiv 2,3 \mod 4.
\end{cases}\\
\eta_\bk&:=\phi_{1,1}\circ\tau_v\circ\lambda_{\bk'}.
\end{align*}
Here, we interpret $\tau_v$ as the identity if $v=0$, which
happens when $k_+'=1$ and $\varepsilon=1$.
Let us show
that the isomorphism
$\eta_\bk: C(p,q)\simarrow C(1,1)\otimes C(p-1,q-1)$
satisfies desired property \eqref{eq:eta}.
\par\noindent
Case i) $\Delta \bk\equiv 0,1\mod 4$. \\
If $k'_+\neq 1$, then $\tau_v$ exchanges $\varepsilon v_{{k'}_+}^+$ and
$v_1^+$. If $k'_+=1$ and $\varepsilon=-1$, then $\tau_v$ changes
the signature of $v_1^+$. If $k'_+=1$ and $\varepsilon=1$,
then $\tau_v=\id$. Thus, $\tau_v\circ\lambda_{\bk'}(V_\bk)=v_1^+$
in all the cases. Applying $\phi_{1,1}$, we have:
$$\phi_{1,1}\circ \tau_v\circ\lambda_{\bk'}(V_\bk)= v_1^+\otimes 1\quad
\in C(1,1)\otimes C(p-1,q-1).$$
Hence, $\eta_\bk$ satisfies \eqref{eq:eta}.

\par\noindent
Case ii) $\Delta \bk\equiv 2,3\mod 4$.\\
The proof goes similarly.
\end{proof}
\begin{rem}
If $\eta_\bk$ is an isomorphism such that $\eta_\bk(V_\bk)=v_1^\delta\otimes 1$,
then we can tell a priori the signature $\delta=\pm$ by \eqref{eq:vk2}.
\end{rem}
\begin{la}\label{la:xi}
Suppose $\bm=(p,q)$, $\bk=(k_+,k_-)$ and~$\bl=(l_+, l_-)$ satisfy
$\bm=\bk+\bl$, $\bk\neq (0,0), (p,q)$, $1\leq p$ and
$\Delta \bm\equiv \Delta \bl\equiv 1\mod 2$.
Then there exists an $\R$-algebra isomorphism
$\xi_\bl: C(p,q)\simarrow C(\bx)\otimes C(\by)$ such that
\begin{equation*}
\xi_\bl(T_\bl)=
\begin{cases}
\id\otimes \tau_1^+ &\IF \Delta \bm\equiv 1\mod 4,\\
\id\otimes \tau_1^- &\IF \Delta \bm\equiv 3\mod 4.
\end{cases}
\end{equation*}
Here, $\bx$ and $\by$ are given by:
\begin{equation*}
\bx=\left\{
\begin{aligned}
(q,p-1) && \Delta \bk\equiv 0 \mod 4,\\
(p-1,q) && \Delta \bk\equiv 2 \mod 4.
\end{aligned}\right. \quad
\by=\left\{
\begin{aligned}
(1,0) && \Delta \bm\equiv 1 \mod 4,\\
(0,1) && \Delta \bm\equiv 3 \mod 4.
\end{aligned}
\right.
\end{equation*}
\end{la}
\begin{proof}
It follows from $\bk\neq (0,0)$, that at least one of $k_+$ and 
$k_-$ is non-zero. We shall give a proof in the case
$k_-\geq 1$. The case $k_-\geq 1$ goes similarly.
We set $\bk':=(k_++1,k_-)$. In light of
$k'_{\pm}\geq 1$ and $\tau^+_{k'_+} T_{\bm-\bk'}=T_\bl$,
we have from Lemma~\ref{la:lambda}:
\begin{equation*}
\lambda_{\bk'}(T_\bl)=
\begin{cases}
\tau^-_{k'_-} &\IF \Delta \bk\equiv 0 \mod 4,\\
\tau^+_{k'_+} &\IF \Delta \bk\equiv 2 \mod 4.
\end{cases}
\end{equation*}
\noindent
(i) Case $\Delta \bk\equiv 0\mod 4$.\\
We set $v=v_{k'_-}^--v_q^-$ and $\xi_\bl:=\phi_{q,p-1}\circ\phi_{1,0}\circ
\tau_v\circ\lambda_{\bk'}$. Here, we regard $\tau_v=\id$ if $v=0$ as
before. In light that $\tau_v$ is an $\R$-algebra isomorphism
of $C(p,q)$ that exchanges $v_{k'_-}^-$ and $v_q^-$,
we have from Lemma~\ref{la:phi}:
$$
\begin{array}{ccccccc}
C(p,q)&\maparrow{\lambda_{\bk'}}&C(p,q)&\maparrow{\tau_v}&C(p,q)&
\maparrow{\phi_{1,0}}&C(q+1,p-1)\\
T_\bl&\mapsto&\tau_{k'_-}^-&\mapsto&\tau_q^-&\mapsto&
\tau_{q+1}^+.
\end{array}
$$
Applying furthermore the isomorphism $\phi_{q,p-1}$, we obtain
\begin{equation*}
\phi_{q,p-1}(\tau_{q+1}^+)=
\begin{cases}
\id\otimes \tau_1^+\in C(q,p-1)\otimes C(1,0) &\IF \Delta \bm\equiv 1 \mod 4,\\
\id\otimes \tau_1^-\in C(q,p-1)\otimes C(0,1) &\IF \Delta \bm\equiv 3 \mod 4.
\end{cases}
\end{equation*}

\par\noindent
(ii) Case $\Delta \bk\equiv 2\mod 4$.\\
We set $v=v_{k'_+}^+-v_p^+$ and
$\xi_\bl:=\phi_{p-1,q}\circ \tau_v\circ\lambda_{\bk'}$.
Then the proof goes similarly.
\end{proof}
\begin{la}\label{la:st_compat}
Suppose $\pi$ is any of the isomorphisms
$\tau_v, \psi_{p,q}, \phi_\bk, \lambda_\bk, \eta_\bk$
and $\xi_\bk$ (see Lemmas~\ref{la:tau}, \ref{la:clifford_hom2}, \ref{la:phi},
\ref{la:df:lambda}, \ref{la:eta} and \ref{la:xi}, respectively for definition).
Then
\begin{equation}
\pi(\st{a})=\st{\pi(a)} \text{ for any $a\in C(p,q)$}.\label{eq:st_compat}
\end{equation}
Here, for a matrix $A\in M(n,\F)\ (\F=\R,\C,\bbH)$, $\st{A}$ denotes
the standard conjugation $\overline{\tr{A}}$.
\end{la}
\begin{proof}
It is sufficient to show \eqref{eq:st_compat} for
the generators $\{v_1^+,\dots, v_p^+, v_1^-,\dots, v_q^-\}$.
\par\medskip\noindent
1) For $\pi=\psi_{p,q}$, \eqref{eq:st_compat} follows straightforward
from definition.
\par\medskip\noindent
2) For $\pi=\phi_\bk$, we first note  $\st{A}A=A\st{A}=1$ if
$A\in C(p,q)\mon$ (see Definition~\ref{df:gpq}).
Furthermore, $\phi_\bk(v_i^\pm)\in C(p,q)\mon$ by definition.
Therefore, we have
$\phi_\bk(v_i^\pm)\st{\phi_\bk(v_i^\pm)}=1$.
On the other hand, it follows from $\st{v_i^\pm}v_i^\pm=1$
that $\phi_\bk(\st{v_i^\pm})\phi_\bk(v_i^\pm)=1$ because
$\phi_\bk$ is an algebra isomorphism.
Hence, we have proved $\phi_\bk(\st{v_i^\pm})=\st{\phi_\bk(v_i^\pm)}$.
The cases $\pi=\tau_v$ and $\lambda_\bk$ are similar.
\par\medskip\noindent
3) For $\pi=\eta_\bk$ and $\xi_\bk$,
\eqref{eq:st_compat} holds because $\eta_\bk$ and $\xi_\bk$
are defined as the composition of $\tau, \phi$ and $\lambda$.
\end{proof}

\subsection{Proof of Proposition~\ref{prop:auto_iso}}\label{subsec:auto_iso}
In this subsection, we shall find the group structure
of $C(p,q)\cap\Aut(J_\pm)$ (see Definition~\ref{df:aut}) and
$G(p,q)$ (see Definition~\ref{df:gpq}). Let us prove the next lemma:
\begin{la}\label{la:even}
We have the following isomorphisms of $\R$-algebras:
\begin{enumerate}
\ritem[1)] $\ceven(p,q)\simeq C(p,q)^{\tau_1^-}\simeq C(p,q-1)
\qquad (q>0).$
\ritem[2)] $\ceven(p,q)\simeq \ceven(q,p)$.
\ritem[3)] $G(p,q)\simeq G(q,p)$.
\ritem[4)] $G(p,q)\simeq C(p,q-1)\cap\Aut(J_+)$ 
if $p$ is even and $q>0$.
\ritem[5)] $G(p,q)\simeq C(p,q-1)\cap\Aut(J_-)$
if $q$ is even and $q>0$.
\end{enumerate}
\end{la}
\begin{proof}
1)
We recall $C(p,q)^{\tau_1^-}\simeq C(p,q-1)$
and $\ceven(p,q)=C(p,q)^{T_{p,q}}$
by Example~\ref{ex:tau}.
Applying $\phi_{0,1}$, we obtain
$$\ceven(p,q)=C(p,q)^{T_{p,q}}\maparrow{\phi_{0,1}}
C(p,q)^{\phi_{0,1}(T_{p,q})}.$$
Since $\phi_{0,1}(T_{p,q})=\tau_1^-$ by Lemma~\ref{la:phi},
the statement (1) follows.
\par\noindent
2) The statement is obvious in the case $(p,q)=(0,0)$.
Suppose now $(p,q)\neq (0,0)$. Then, without loss of generality,
we may assume $q\geq 1$. In view of
$\ceven(p,q)=C(p,q)^{T_{p,q}}$,
the isomorphism (2) is obtained as the composition:
$$C(p,q)^{T_{p,q}}\maparrow{\phi_{0,1}} C(p,q)^{\tau_1^-}\simeq C(p,q-1)
\simeq C(p+1,q-1)^{\tau_1^+}\maparrow{\phi_{1,0}} C(q,p)^{T_{q,p}}.$$
Here, the last isomorphism follows from $\phi_{1,0}(\tau_1^+)=T_{q,p}$
by Lemma~\ref{la:phi}.
\par\noindent
3) We denote by $\pi:\ceven(p,q)\simarrow \ceven(q,p)$
the isomorphism constructed in (2). Then
$\pi(\st{a})=\st{\pi(a)}$ for any $a\in\ceven(p,q)$ by
Lemma~\ref{la:st_compat}.
Therefore, to see $\pi(\tr{a})=\tr{\pi(a)}$, it is sufficient to
show that $\pi$ transfers $T_{p,0}$ to $T_{0,p}$ because
$\tr{a}=T_{p,0}(\st{a})$ by \eqref{eq:st2}.
Chasing the automorphisms in each step:
$$
\begin{array}{ccccccc}
\ceven(p,q)&\maparrow{\phi_{0,1}}&C(p,q)^{\tau_1^-}&
\simeq& C(p+1,q-1)^{\tau_1^+}&\maparrow{\phi_{1,0}}&
\ceven(q,p)\\
\tau_1^+\cdots\tau_p^+&\mapsto&\tau_1^+\cdots\tau_p^+&\mapsto
&\tau_2^+\cdots\tau_{p+1}^+&\mapsto&\tau_1^-\cdots\tau_p^-
\end{array}
$$
by using Lemma~\ref{la:phi}, we see $\pi$ transfers $T_{p,0}$ to $T_{0,p}$.
Hence (3) follows.
\par\noindent
4) We notice $J_+\in \ceven(p,q)$ if $p$ is even.
Then, a simple computation shows that the isomorphism
$\ceven(p,q)\simarrow C(p,q-1)$ constructed in (1) sends
$J_+\in\ceven(p,q)$ to $J_+\in C(p,q-1)$. Hence,
the statement (4) follows.
\par\noindent
5) The proof is similar to (4).
\end{proof}
\begin{rem}
Lemma~\ref{la:even} implies, in particular, that the groups
$C(p,q)\cap\Aut(J_+)$ and $C(p,q)\cap\Aut(J_-)$ are isomorphic
if $p$ is even and $q$ is odd. In fact, these two groups
coincide (see Remark~\ref{rem:pq_odd}).
\end{rem}
In light of Lemma~\ref{la:even}, we can tell the group structure of $G(p,q)$
in the case $pq$ is even by finding the group structure of
$C(p,q-1)\cap\Aut(J_\pm)$. This is given
in the following proposition:

\begin{prop} \label{prop:autj}
For $p,q\geq 1$, we set $\alpha$ as in the table below
and $n:=2^{(p+q-\alpha)/2}$.
Then, the group $C(p,q)\cap\Aut(J_+)$ is isomorphic to one of
the following Lie groups.
$$
\arraycolsep 15pt
\begin{array}{c|ccc}
C(p,q)\cap\Aut(J_+) & p \mod 4 & \alpha & p-q-1 \mod 8\\ \hline
O(n,n)^2 & 0,1 &\pmid{3} & \pmid{0}\\
Sp(n,\R)^2 & 2,3\\\hline
O(n,n) & 0,1 & \pmid{2} & \pmid{\pm 1}\\
Sp(n,\R) & 2,3 \\ \hline
U(n,n) &\text{any}  & 3 & \pm 2\\ \hline
Sp(n,n) & 0,1 & \pmid{4} & \pmid{\pm 3}\\
O^*(4n) & 2,3 \\ \hline
Sp(n,n)^2 & 0,1 & \pmid{5} & \pmid{4}\\
O^*(4n)^2 & 2,3
\end{array}
$$
Likewise, the group $C(p,q)\cap\Aut(J_-)$ is given,
if we replace the condition $p\mod 4$ by $q-1\mod 4$ in the
above table.
\end{prop}
\begin{proof}
Combining the isomorphism $\psi_{1,1}\simarrow M(2,\R)$
(see Lemma~\ref{la:clifford_hom2}) with
$\eta_{p,0} :C(p,q)\simarrow C(1,1)\otimes C(p-1,q-1)$
(see Lemma~\ref{la:eta}), we obtain an isomorphism
$\iota:C(p,q)\simarrow M(2,\R)\otimes C(p-1,q-1)$ satisfying
$$
\iota(J_+)=
\begin{cases}
\twomatrix{1}{ }{ }{-1}\otimes 1 &\IF p\equiv 0,1\mod 4,\\
\twomatrix{ }{-1}{1}{ }\otimes 1 &\IF p\equiv 2,3\mod 4.
\end{cases}
$$

\par\noindent
Case i) $p-q\equiv 0,2,3\mod 4$.\\
Combining the isomorphism $C(p-1,q-1)\simeq M(n,\F$)
 ($\F=\R, \C$ or $\bbH$) given in Proposition~\ref{prop:clifiso},
we find that the group $C(p,q)\cap\Aut(J_+)$ is isomorphic to
$O(n,n)$, $Sp(n,\R)$ $(\F=\R);$ or $Sp(n,n), O^*(4n)$ $(\F=\bbH)$.
It is isomorphic to $U(n,n)$ (independently of $p\mod 4$) if $\F=\C$.
\par\noindent
Case ii) $p-q\equiv 1\mod 4$.\\
We note $C(p,q)\simeq M(2n,\F)^2$ ($\F=\R$ or $\bbH$).
The proof goes similarly.

This completes the proof for $C(p,q)\cap\Aut(J_+)$. For the
group structure of $C(p,q)\cap\Aut(J_-)$, we use $\eta_{0,q}$ instead
of $\eta_{p,0}$.
\end{proof}

Owing to Lemma~\ref{la:even} (4) and (5), this Proposition
completes the proof of Proposition~\ref{prop:auto_iso}
in the case $pq$ is even.

Next, let us consider the case where both $p$ and $q$ are odd.
Owing to Lemma~\ref{la:gpq_aut}, Proposition~\ref{prop:auto_iso}
in this case is reduced to
the group structure of $\ceven(p,q)\cap\Aut(J_\pm)$. Here
is a result:
\begin{prop}\label{prop:odd}
Suppose $pq$ is odd. We set
$n=2^{(p+q-\alpha)/2}$ and $\alpha\in\{2,4\}$ is given in the table below.
Then, the group $G(p,q)\simeq \ceven(p,q)\cap\Aut(J_\pm)$ is isomorphic to
the following Lie groups:
$$
\begin{array}{c|ccc}
\ceven(p,q)\cap\Aut(J_\pm) & p-q\mod 8 & p+q\mod 8 & \alpha\\ \hline
GL(n,\R) & 0 & \pmid{\pm 2} & 2\\
GL(n,\bbH) & 4 &  & 4\\ \hline
O(n,\C) & \pmid{\pm 2} & 0 & 2\\
Sp(n,\C) & & 4 & 4
\end{array}
$$
\end{prop}
\begin{proof}
We first note that the isomorphism
$\phi_{0,1}: \ceven(p,q)\simarrow C(p,q-1)$ transfers the automorphism
$T_{p,0}$ of $\ceven(p,q)$ to the automorphism $T_{p,0}$ of
$C(p,q-1)$ (see Lemma~\ref{la:even}(1) and Lemma~\ref{la:phi}(2)(iv)).
Now, applying Lemma~\ref{la:xi} to the case where
$\bm=(p,q-1), \bl=(p,0)$ and $\bk=(0,q-1)$, we obtain an isomorphism
\begin{equation}
\xi_{p,0}: C(p,q-1)\simarrow C(\bx)\otimes C(\by)\label{eq:tp0}
\end{equation}
satisfying
$$\xi_{p,0}(T_{p,0})=
\begin{cases}
\id\otimes \tau_1^+ & \IF p-q\equiv 0 \mod 4,\\
\id\otimes \tau_1^- & \IF p-q\equiv 2 \mod 4.
\end{cases}
$$
Here, we have set
\begin{equation*}
\bx:=\left\{
\begin{aligned}
(p,q-2) &&  q\equiv 1 \mod 4,\\
(p-1,q-1) &&  q\equiv 3 \mod 4,
\end{aligned}\right. \quad
\by:=\left\{
\begin{aligned}
(1,0) &&  p-q\equiv 0 \mod 4,\\
(0,1) &&  p-q\equiv 2 \mod 4.
\end{aligned}
\right.
\end{equation*}
\par\noindent
Case i) $p-q\equiv 0\mod 4$.\\
We note
$\xi_{p,0}(T_{p,0})=\id\otimes\tau_1^+$.\\
In this case, $C(\by)=C(1,0)\simeq \R\oplus\R$ and $\tau_1^+$ acts on $C(\by)$ as
$(a,b)\mapsto (b,a)$. We now assert:
\par\noindent
Claim. $C(\bx)\simeq M(n,\F)$. \qquad Here, $\F$ is given by\\
\begin{equation}
\F=
\begin{cases}
\R &\IF p-q\equiv 0\mod 8,\\
\bbH &\IF  p-q\equiv 4\mod 8.
\end{cases}\label{eq:dff}
\end{equation}
This claim is a direct consequence of Proposition~\ref{prop:clifiso}.
Combining the isomorphism $\xi_{p,0}$ (see \eqref{eq:tp0}) with
$C(\bx)\simeq M(n,\F)$ and
$C(\by)\simeq \R\oplus\R$, we obtain an isomorphism
$\iota: \ceven(p,q)\simeq M(n,\F)\oplus M(n,\F)$.
In light of $\tr{a}=T_{p,0}(\st{a})$ for $a\in\ceven(p,q)$ (see \eqref{eq:st2}),
we have $\iota(\tr{a})=(\st{B},\st{A})$ if we write $\iota(a)=(A,B)$.
Hence, the following
conditions on $a\in\ceven(p,q)$ are equivalent:
\begin{align*}
a\in G(p,q) &\iff \tr{a}{a}=1\\
&\iff (\st{B},\st{A})(A,B)=(1,1)\\
&\iff \st{B}=A^{-1}.
\end{align*}
This shows that the group $G(p,q)$ is isomorphic to 
the general linear group $GL(n,\F)$, where
$\F=\R$ or $\bbH$ as in \eqref{eq:dff}.
\par\noindent
Case ii) $p-q\equiv 2\mod 4$.\\
 We note $\xi_{p,0}(T_{p,0})=\id\otimes
\tau_1^-$.
In this case, $C(\by)=C(0,1)\simeq \C$, and $\tau_1^-$ acts on
$C(\by)$ as complex conjugation. It follows from our assumption
$q\equiv 1\mod 2$ and $p-q\equiv 2\mod 4$ that
$p+q\equiv 0\mod 4$, and therefore $p+q\mod 8$ is either $0$ or $4$.
\par\noindent
Subcase ii-a) $p+q\equiv 0 \mod 8$.\\
If $q\equiv 1\mod 4$, $\Delta \bx=q-p=2q-(p+q)\equiv 2\mod 8$.
If $q\equiv 3\mod 4$, $\Delta \bx=p-q=(p+q)-2q\equiv 2 \mod 8$.
Hence, in either case, $C(\bx)$ is isomorphic to $M(n,\R)$ by
Proposition~\ref{prop:clifiso}, and thus we obtain an isomorphism
$\iota: \ceven(p,q)\simarrow M(n,\C)$.
Then by using $T_{p,0}(\st{a})=\tr{a}$ again and $\xi_{p,0}(T_{p,0})=\id\otimes
\tau_1^-$, we have
$$\iota(\tr{a})=\iota(T_{p,0}(\st{a}))=\overline{\st{\iota(a)}}=\tr{\iota(a)}.
$$
This shows that the group, $G(p,q)$ is isomorphic to $O(n,\C)$.
\par\noindent
Subcase ii-b) $p+q\equiv 4\mod 8$.\\
Similarly to the subcase (ii-a), we have an isomorphism
$C(\bx)\simeq M(n,\bbH)$.
Via the isomorphism $M(n,\bbH)\otimes\C\simarrow M(2n,\C)$
defined by
$$
(P+Q j)\otimes 1
\mapsto
\twomatrix{P}{-Q}{\overline{Q}}{\overline{P}},\qquad
1\otimes z \mapsto
\twomatrix{z I_n}{ }{ }{z I_n},$$
$\id\otimes\tau_1^-$ induces the involution of $M(2n,\C)$ given
by
$$
z\twomatrix{P}{-Q}{\overline{Q}}{\overline{P}}\mapsto
\overline{z}\twomatrix{P}{-Q}{\overline{Q}}{\overline{P}}
=\overline{z\twomatrix{\overline{P}}{-\overline{Q}}{Q}{P}},
$$
that is,
$$\twomatrix ABCD\mapsto 
\twomatrix{\overline{D}}{-\overline{C}}{-\overline{B}}{\overline{A}}.$$
Therefore, the resulting isomorphism
$\iota: \ceven(p,q)\simarrow M(2n,\C)$ has the following property:
$$\iota(T_{p,0}(a))=\twomatrix{\overline{D}}{-\overline{C}}{-\overline{B}}{\overline{A}},$$
if we write $\iota(a)=\twomatrix{A}{B}{C}{D}\in M(2n,\C)$
for $a\in\ceven(p,q)$.
Then, by using $\tr{a}=T_{p,0}(\st{a})$, we obtain
$$\iota(\tr{a})=\twomatrix{\tr{D}}{-\tr{B}}{-\tr{C}}{\tr{A}}
=J_n^{-1}\left(\tr{\iota(a)}\right)J_n,$$
where $J_n:=\twomatrix{ }{-I_n}{I_n}{ }\in M(2n,\C)$.
Hence, for $a\in\ceven(p,q)$,
$$a\in G(p,q) \iff \tr{a}a=1 \iff \tr{\iota(a)}J_n\iota(a)=J_n,$$
and thus $G(p,q)$ is isomorphic to $Sp(n,\C)$. This
completes the proof of Proposition in all the cases.
\end{proof}

\subsection{Proof of Theorem~\ref{thm:spin_ck_form}}
\label{subsec:clifford_proper}
This subsection gives a proof of Theorem~\ref{thm:spin_ck_form}.
We have already seen $d(G)=d(H)+d(L)$ (the condition \eqref{thm:red:suf}(b))
in Table~\ref{subsec:spin_cck}.

Let us verify the condition \eqref{thm:red:suf}(a). 
We write $\g$ for the Lie algebra of $G(p,q)$.
Then, $\g=\set{X\in\ceven(p,q)}{\tr{X}+X=0}$, and
$X\mapsto -\st{X}$ defines a Cartan involution of $\g$.
Let $\g=\fk+\p$ the corresponding Cartan decomposition, which
then induces a Cartan decomposition
$\fl=(\fl\cap\fk)+ (\fl\cap\p)$ of
 the Lie algebra $\fl$ of $Spin(p,q)$.
\begin{la} \label{la:spin_invertible}
Suppose $q\geq 1$, and $(\pi,W)$ is a
representation of the $\R$-algebra $\ceven(1,q)$. 
We note $\fl=\spin(1,q)\subset \ceven(1,q)$.
Then
$\pi(X)\in \End(W)$ is invertible for any non-zero $X\in \fl\cap\p$.
\end{la}
\begin{proof}
Retain the notation as in Section~\ref{subsec:spin_def_la}. We
set $X_0:=v_1^+v_1^-$. Then $X_0\in\p$ because $\st{X_0}=\st{(v_1^-)}
\st{(v_1^+)}=X_0$. As $X_0^2=1$, we have 
$\pi(X_0)^2=\pi(1)=\id$ because $\pi$ is a representation
of an $\R$-algebra. Hence, $\pi(X_0)$ is invertible.
For $t\in\R$, we set
$$x(t):=\exp(t X_0)=(\cosh t)+(\sinh t)v_1^+v_1^-=v_1^+((\cosh t)v_1^+
+(\sinh t)v_1^-).$$
Then, $x(t)\in Spin(1,q)$ by Definition~\ref{df:gpq}, and
$\set{x(t)}{t\in\R}$ forms a
one-parameter subgroup of $Spin(1,q)$ consisting of hyperbolic elements.
Since $\R$-$\rank( Spin(1,q))=1$, any
non-zero $X\in \fl\cap\p$ can be written as
$$ X=\Ad(l) t X_0$$
for some $l\in L\cap K$ and some non-zero $t\in\R$. Hence,
$\pi(X)^2=t^2\id$ and we have seen $\pi(X)$ is invertible.
\end{proof}

Then, the proof of Theorem~\ref{thm:spin_ck_form} will be completed
if we show the following:
\begin{prop} \label{prop:clifford_proper}
Each triple of Lie groups $(G,H,L)$ in the following
table satisfies the condition
\eqref{thm:red:suf}(a), or equivalently, $\properin{L}{H}{G}$.
Here $\alpha\in\{0,1,2,3\}$ is given in the table, and
$n=2^{4r-\alpha}$.
$$
\begin{array}{|ccc|cc|} \hline
 G & H & L & W & \alpha\\\hline
GL(n,\R) & GL(n-1,\R) &Spin(1,8r+1)  & \R^n & 0\\
Sp(n,\R) & Sp(n-1,\R) &Spin(1,8r+2)  & \R^{2n} & 0\\
Sp(n,\C) & Sp(n-1,\C) &Spin(1,8r+3)  & \C^{2n} & 0\\
Sp(n,n) & Sp(n-1,n) &Spin(1,8r+4) & \bbH^{2n} & 0\\
GL(n,\bbH) & GL(n-1,\bbH) &Spin(1,8r+5)  & \bbH^n & 1\\
O^*(2n) & O^*(2n-2) &Spin(1,8r+6) & \bbH^n & 2\\
O(n,\C) & O(n-1,\C) &Spin(1,8r+7) & \C^{n} & 3\\
O(n,n) & O(n-1,n) &Spin(1,8r+8)   & \R^{2n} & 3\\\hline
\end{array}
$$
\begin{center}
\bf Table~\thesubsection
\end{center}
\end{prop}
\begin{proof}
Let $W$ be the linear space as a given in Table~\thesubsection, and 
$\{e_i\}$ the standard basis of $W$. We note that the
subgroup $H$ in Table~\thesubsection\ is given by
$$ H=\set{g\in G}{g(e_1)=e_1},$$
and thus its Lie algebra $\h$ is given by
$\h=\set{X\in\g}{Xe_1=0}$. In particular, 
any element in $\h$ is not invertible in $\End(W)$.

Let $\fa\subset\p$ be a maximally split abelian subalgebra of $\g$.
Then $\fa_H$ consists of non-invertible elements in $\End(W)$,
while any non-zero element of $\fa_L$ is invertible by
Lemma~\ref{la:spin_invertible}. Thus we have
$$ \fa(L)\cap\fa(H)=0.$$
The condition \eqref{thm:red:suf}(a) is proved.
Thus, Proposition~\ref{prop:clifford_proper} is proved.
\end{proof}

\subsection{Compact Clifford-Klein form of $SO(8,\C)/SO(7,1)$}
\label{subsec:o8c}
This subsection proves the case $10$ of Theorem~\ref{thm:admit_cck},
that is:
\begin{thm}\label{thm:o8c}
There exists a compact Clifford-Klein
form of the symmetric space $SO(8,\C)/SO(7,1)$.
\end{thm}
The main machinery for the proof here is the Clifford algebra associated to
indefinite quadratic forms and Theorem~\ref{thm:red:suf}.
We shall show that the triple of Lie groups $(G,H,L):=(O(8,\C),
O(7,1), Spin(7,\C))$ satisfies the conditions
\eqref{thm:red:suf}(a) and (b).

First, the embedding $L\subset G$ is obtained as a complexification
of the compact real form $Spin(7)\subset G(0,7)\simeq O(8)$
(see Table~\ref{subsec:clifford_algebra}.3).

In order to verify the condition \eqref{thm:red:suf}(a), let us compute
$\fa_L$ in $\fa$. To do this, we first write explicitly
the group isomorphism $G(0,7)\simeq O(8)$, which is
obtained as the restriction of the $\R$-algebra isomorphism
$\iota:\ceven(0,7)\simarrow M(8,\R)$ to $G(0,7)$. We recall
from Proposition~\ref{prop:clifiso} that $\iota$ is defined as
the composition of the following isomorphism:

\begin{equation}
\begin{split}
\ceven(0,7)&\maparrow{\phi_{0,1}} C(0,6)\maparrow{\phi_{0,3}} C(3,3)
\maparrow{\phi_{1,1}} C(1,1)\otimes C(1,1)\otimes C(1,1)\\
\maparrow{\psi_{1,1}}& M(2,\R)\otimes M(2,\R)\otimes M(2,\R)
\simeq M(8,\R).
\end{split}
\label{eq:o8c}
\end{equation}
Let $E_{ij}$ be the matrix unit, and we choose
a Cartan subalgebra of the Lie algebra $\fo(8)$ in $M(8,\R)$,
$\ft=\bigoplus_{k=1}^4 \R f_k$ 
by taking
$$ f_k:=E_{2k,2k-1} - E_{2k-1,2k}\qquad (k=1,2,3,4).$$
Then $\fa:=\sqrt{-1}\,\ft$ becomes a maximally split abelian subalgebra
of $\fo(8,\C)=\fo(8)+\sqrt{-1}\,\fo(8)$.
We shall identify $\ft$ with $\R^4$ by means of the basis
$ f_1,\dots, f_4$. We assert:
\par\medskip\noindent
Claim: \quad
$\displaystyle \fa(L)=\sqrt{-1}\!\!\bigcup_{\oalign{\scriptsize$\varepsilon_i\in\{1,-1\}$\cr
\scriptsize$\varepsilon_1 \varepsilon_2 \varepsilon_3 \varepsilon_4=1$}}
\set{\tr{(a,b,c,d)}}{ \varepsilon_1 a +\varepsilon_2 b +
\varepsilon_3 c+ \varepsilon_4 d=0}.$

To see Claim, we set
$$x_1:=\tr{(1,1,1,1)},\quad x_2:=\tr{(1,-1,1,-1)},\quad
x_3:=\tr{(1,-1,-1,1)}.$$
Then, $\iota^{-1}(x_k)$ ($k=1,2,3$) is computed by using Lemmas
\ref{la:clifford_hom2}, \ref{la:phi} and~\ref{la:even}(1) as follows:
\def\fourdiagmatrix #1,#2,#3,#4.{\diag(#1, #2, #3,#4)}
\def\imatrix{\twomatrix{ }{-1}{1}{ }}
\begin{gather*}
\begin{array}{ccccc}
M(8,\R)&\simeq&M(2,\R)\otimes M(4,\R)&\simeq&C(1,1)\otimes C(1,1)\otimes
C(1,1)\\
x_1&\mapsto&\imatrix\otimes
\fourdiagmatrix 1,1,1,1.&\mapsto & v_1^-\otimes 1\otimes 1\\
x_2&\mapsto&\imatrix\otimes
\fourdiagmatrix 1,-1,1,-1.&\mapsto & v_1^-\otimes v_1^+\otimes 1\\
x_3&\mapsto&\imatrix\otimes
\fourdiagmatrix 1,-1,-1,1.&\mapsto & v_1^-\otimes v_1^+\otimes v_1^+
\end{array}\\[10pt]
\begin{array}{cccccc}
\simarrow& C(3,3)&\simarrow& C(0,6)&\simarrow &\ceven(0,7)\\
\mapsto& v_1^- &\mapsto & v_1^-&\mapsto & v_1^- v_2^-\\
\mapsto& v_1^+v_2^+&\mapsto & v_4^- v_5^-&\mapsto & v_5^-v_6^-\\
\mapsto& v_1^-v_2^-v_3^+&\mapsto& v_3^-v_6^-&\mapsto & v_4^- v_7^-
\end{array}
\end{gather*}
Then, $v_1^-v_2^-$, $v_5^-v_6^-$, $v_4^-v_7^-$ are
elements of the Lie algebra $\fl$ of $L=Spin(7)$, as it follows from
$$\exp(tv_1^-v_2^-)=\cos t+(\sin t)(v_1^-v_2^-)
=v_1^-(-(\cos t)v_1^- +(\sin t)v_2^-)$$
and alike, where we used $(v_1^-v_2^-)^2=-1$. Furthermore,
\def\ee#1#2{$v_#1^-v_#2^-$}
since \ee12, \ee56 and \ee47 mutually commute,
we conclude that $\R x_1\oplus \R x_2\oplus \R x_3$ is
a Cartan subalgebra of $\fl$
because $\rank L=3$. Hence, $\sqrt{-1}(\R x_1\oplus\R x_2\oplus\R x_3)$
is a maximally
split abelian subspace of $\fl_\C$.
Since this subspace is contained in $\fa$, $\fa(L)$ is given by
its $W_G$-orbit, namely, $\sqrt{-1}\,\bigcup_{w\in W_G} w\set{\tr{(a,b,c,d)}}
{a+b-c-d=0}$, where $W_G$ is the Weyl group of $O(8,\C)$.
Hence Claim is proved.

On the other hand, it is easy to see
\begin{align*}
\fa(H)&=W_G(\sqrt{-1}\R f_1)\\
&=\sqrt{-1}\R\tr{(1,0,0,0)}\cup\sqrt{-1}
\R\tr{(0,1,0,0)}\cup\sqrt{-1}\R\tr{(0,0,1,0)}\cup
\sqrt{-1}\R\tr{(0,0,0,1)}.
\end{align*}
Thus, 
$ \fa(L)\cap \fa(H)=0$, and therefore the condition \eqref{thm:red:suf}(a) is proved.
Since the condition \eqref{thm:red:suf}(b) is clear from $d(G)=28=7+21=d(L)+d(H)$,
we conclude that the symmetric space $G/H=O(8,\C)/O(7,1)$ (therefore,
$SO(8,\C)/SO(7,1)$ also) admits
a compact Clifford-Klein form by Theorem~\ref{thm:red:suf}.
Hence, Theorem~\ref{thm:o8c} is proved.

\section{Tangential space form problem}
\label{sec:tangential}
In this section, we introduce the tangential homogeneous space
$G_\theta/H_\theta$ associated to a homogeneous space $G/H$ of
reductive type (see Definition~\ref{dfn:tangential}).
This is also a symmetric space if $G/H$ is symmetric.
We shall see that most of the results on discontinuous groups
for $G/H$ can be generalized or strengthened to those for the tangential
homogeneous spaces $G_\theta/H_\theta$.
Furthermore, if Conjecture~\ref{conj:inv:suff} is true, then
the existence of compact Clifford-Klein forms of $G/H$
implies that of $G_\theta/H_\theta$. Thus, the solution
to Problem A for $G_\theta/H_\theta$ might be a good approximate
to that of the original problem $G/H$, even though
there is no general theory that bridges directly between
$G/H$ and $G_\theta/H_\theta$ for Problem A.

Because the `non-compact part' of $G_\theta$ is abelian,
a simple and detailed study becomes possible for $G_\theta/H_\theta$.
In particular, we shall 
prove Theorem~\ref{thm:tangential}, which gives an answer to
the tangential version of our space form conjecture 
(Conjecture~\ref{conj:spherical:group}).

\newcommand{\ct}{\overline{\theta}}
\subsection{Tangential symmetric space}
Let $G$ be a linear reductive Lie group with
a Cartan involution $\theta$. We write the Cartan
decomposition as $G=(\exp\p)K$, and form
a semidirect group
$$G_\theta:=K\ltimes_\Ad\p$$
with multiplication: 
$(k_1,X_1)(k_2,X_2):=(k_1k_2, X_1+\Ad(k_1)X_2)$. We say
$G_\theta$ is the {\it Cartan motion group} of $G$.

We define an involutive automorphism
$\ct$ of the Cartan motion group $G_\theta=K\ltimes\p$ by
$$\ct: G_\theta\rightarrow G_\theta,\qquad (k,X)\mapsto (k,-X).$$
By using the Cartan decomposition $G=(\exp\p)K$, we define a diffeomorphism
\begin{equation}
\Phi_\theta:G\rightarrow G_\theta,\quad  e^Xk\mapsto (k,X).
\end{equation}
Clearly, we have:
$$\ct\Phi_\theta=\Phi_\theta\theta.$$
\begin{rem}
In defining $\Phi_\theta$, we have used the Cartan decomposition
$G=(\exp\p)K$ instead of $G=K\exp\p$.
The advantage of our definition is that the following natural
formula holds:
$$
\Phi_\theta(k_1)\Phi_\theta(g)\Phi_\theta(k_2)=\Phi_\theta(k_1gk_2)
\qquad (g\in G,\ k_1, k_2\in K).
$$
\end{rem}
$\Phi_\theta$ is not a group homomorphism, and in general,
the image $\Phi_\theta(H)$ of a subgroup $H$ is not
a subgroup of $G_\theta$. However, if $H$ is a
connected $\theta$-stable subgroup, then
$\Phi_\theta(H)$ becomes a $\ct$-stable subgroup of $G_\theta$.
In this case, $H$ is a reductive Lie group with the Cartan involution $\theta|_H$,
and the Cartan motion group $H_\theta$ is naturally isomorphic to
$\Phi_\theta(H)$.
Thus, we can and do regard $H_\theta$ as a closed subgroup of $G_\theta$.
The same argument is valid if we allow $H$ to have finitely many
connected components.
\begin{dfn}\label{dfn:tangential}
We say the homogeneous space $G_\theta/H_\theta$ is the
{\it tangential homogeneous space} of $G/H$.
Furthermore, if $G/H$ is a symmetric space, so is $G_\theta/H_\theta$,
which we say the {\it tangential symmetric space} of $G/H$.
\end{dfn}

This section studies
Problem A for the tangential
symmetric space $G_\theta/H_\theta$.
In contrast to a reductive symmetric space $G/H$,
we can prove affirmatively a tangential
version of Conjecture~\ref{conj:inv:suff}
(see Theorem~\ref{thm:carmot:criterion}). 
In particular, this enables us to solve completely
the tangential analog of the space form problem,
for which the original version has not yet found a
final answer (see Section~\ref{subsec:space_form2}).

\subsection{Discontinuous groups for tangential homogeneous spaces}
This subsection reveals a strong resemblance on
discontinuous group between a reductive homogeneous space $G/H$
and its tangential homogeneous space $G_\theta/H_\theta$.
The proof for $G_\theta/H_\theta$ is easier or analogous to those for $G/H$.
\begin{enumerate}
\ritem[1)] $G_\theta/H_\theta$ is diffeomorphic to $G/H$. Furthermore,
both of them are isomorphic to $K\times_{K\cap H}(\mathfrak{q}\cap\p)$ as
$K$-equivariant fiber bundles over the compact homogeneous space $K/K\cap H$
(see \eqref{eq:four:decomp} for notation).
\ritem[2)]  
$\properin{L}{H}{G} \iff \properin{L_\theta}{H_\theta}
{G_\theta}$\quad (see Lemmas~\ref{la:reductive} and \ref{la:cartan}).
\ritem[3)] $\similarin{L}{H}{G} \iff \similarin{L_\theta}{H_\theta}
{G_\theta}$.
\ritem[4)] $d(G)=d(G_\theta)$,\quad $d(H)=d(H_\theta)$.
\ritem[5)] $L\backslash G/H$ is compact $\iff$ $L_\theta\backslash G_\theta
/H_\theta$ is compact, provided $L$ and $H$ are $\theta$-stable connected
subgroups such that $\properin{L}{H}{G}$, or equivalently,
$\properin{L_\theta}{H_\theta}{G_\theta}$
(see Lemmas~\ref{la:dg} and \ref{la:tan:dg}).
\ritem[6)] If $G/H$ is in Table~\ref{subsec:construction},
both $G/H$ and $G_\theta/H_\theta$ admit compact Clifford-Klein forms.
\ritem[7)] If $G/H$ satisfies the assumption of
Theorem~\ref{thm:non-cpt-dim}, neither $G/H$ nor $G_\theta/H_\theta$ admits
compact Clifford-Klein forms.
\ritem[8)] If $\rrank G=\rrank H$, then neither $G/H$ nor $G_\theta/H_\theta$
admits infinite discontinuous groups (cf. Theorem~\ref{thm:calabi2}).
\end{enumerate}

Though $G/H$ is diffeomorphic to $G_\theta/H_\theta$,
discontinuous groups for $G/H$ are quite different from
those for $G_\theta/H_\theta$ as abstract groups in general.
In fact, we shall deal essentially with abelian discontinuous subgroup
in the case of a Cartan motion group $G_\theta$, as is
seen in Bieberbach's theorem on discrete subgroups in $O(n)\ltimes\R^n$.

The next subsection gives a result that reduces Problem A
to the existence problem of
a certain {\bf connected} subgroup (see Theorem~\ref{thm:carmot:criterion}).
This will be possible for $G_\theta/H_\theta$,
in contrast to a reductive symmetric space $G/H$ for which such statement
is still conjectural.

\medskip
Here are the correspondences between $G/H$ and $G_\theta/H_\theta$:
\begin{center}
\renewcommand{\arraystretch}{1.3}
\begin{tabular}{cccl}
$G/H$ && $G_\theta/H_\theta$\\
Conjecture~\ref{conj:inv:suff}&\hspace{8pt} & Theorem~\ref{thm:carmot:criterion} &
: \quad Reduction to connected case.\\
Conjecture~\ref{conj:spherical}& & Theorem~\ref{thm:tangential}&
: \quad Space form problem.
\end{tabular}
\renewcommand{\arraystretch}{1}
\end{center}

\subsection{
Reduction to connected subgroups}
\label{subsec:tan:key}
Hereafter, we shall denote a Cartan motion group by $G$ instead of
$G_\theta$. The main result of this subsection is
Theorem~\ref{thm:carmot:criterion}. We also reformulate this result into
Theorem~\ref{thm:tan:simple} under the assumption that $H$ is connected
and $\ct$-stable, which
will then play a key role in our solution to the tangential
version of the space form problem 
(see Section~\ref{subsec:tangential_space}).

We start with a general case where $H$ is not necessarily
$\ct$-stable.
\begin{thm}\label{thm:carmot:criterion}
Let $G$ be a Cartan motion group and $H$ a closed subgroup of $G$.
Then the following two conditions are equivalent:
\begin{enumerate}
\ritem[(i)] A homogeneous space $G/H$ admits compact
Clifford-Klein forms.
\ritem[(ii)] There exists a $\ct$-stable connected subgroup $L$ satisfying
the following two conditions (a) and (b).

\smallskip
\noindent
\begin{tabular}{ll}
{\rm \eqref{thm:carmot:criterion}(a)}&\ $\properin{L}{H}{G}$.\\
{\rm \eqref{thm:carmot:criterion}(b)}&\ $L\backslash G/H$ is compact.
\end{tabular}
\end{enumerate}
\end{thm}
\stepcounter{equation}
Note that the implication ``(ii)$\Rightarrow$(i)''
is a `tangential' analog of
Theorem~\ref{thm:sufficient}, whereas
the implication ``(i)$\Rightarrow$(ii)'' is a `tangential' solution to
Conjecture~\ref{conj:inv:suff} for a homogeneous space
of reductive type.
The proof of Theorem~\ref{thm:carmot:criterion} will be given in
Section~\ref{subsec:pf:tan:suf}.

As a second formulation of Theorem~\ref{thm:carmot:criterion},
we shall reformulate it in a more explicit way
for $\ct$-stable $H$. To do this, we introduce
a tangential analog of $\fa(L)$ given in \eqref{eq:df:fl}
as follows:
Let $G=K\ltimes\p$ be a Cartan motion group.
Recall from Section~3.2 that $\fa$ was taken as a maximal
abelian subspace of $\p$ (where, we remember the original
reductive Lie algebra $\fk+\p$). Then, for a
subset $L$ in the Cartan motion group $G=K\ltimes\p$, we define
\begin{equation}
\fa(L):=KLK\cap\fa.\label{eq:df:tan:fa}
\end{equation}
Then Theorem~\ref{thm:carmot:criterion} amounts to:
\begin{thm}\label{thm:tan:simple}
Let $G$ be a Cartan motion group, and let $H$ be a $\ct$-stable
subgroup with at most finitely many connected components.
Then, the following two conditions are equivalent:
\begin{enumerate}
\ritem[(i)] A homogeneous space $G/H$ admits compact Clifford-Klein forms.
\ritem[(ii)] There exists a subspace $W$ in $\p$ satisfying
the following two conditions (a) and (b).

\smallskip
\noindent
\begin{tabular}{ll}
{\rm \eqref{thm:tan:simple}(a)} &\ $\fa(W)\cap\fa(H)=\z$.\\
{\rm \eqref{thm:tan:simple}(b)} &\ $\dim(W)+d(H)=d(G)$.
\end{tabular}
\end{enumerate}
\end{thm}
We prove Theorem~\ref{thm:tan:simple} by admitting
Theorem~\ref{thm:carmot:criterion}. First, we prepare
some necessary lemmas for this.

\begin{la}\label{la:rectangle}
Let $G=K\ltimes\p$ be a Cartan motion group,
and $H$ a $\ct$-stable connected subgroup.
Then:
\begin{enumerate}
\ritem[1)] $H=(K\cap H)\ltimes (\p\cap H)$.
\ritem[2)] $\p\cap H$ is a subspace of $\p$.
\ritem[3)] $H$ has always a uniform lattice.
\ritem[4)] $H \ssimilar \p\cap H$.
\end{enumerate}
\end{la}
\begin{proof}

1) We note that $\exp:Lie(\p)\rightarrow\p$ is the identity map if we
identify $Lie(\p)$ with $\p$ as usual.
Let $H'$ be
the image of $(K\cap H)\times (Lie(\p)\cap\h)$ by
the following diffeomorphism:
$$\phi:=\id\times\exp: K\times Lie(\p)\rightarrow K\ltimes\p,
\quad (k,X) \mapsto (k,X).$$
Clearly, $H'$ is closed.
We claim that $H'$ is an open submanifold of $H$.
This is because the Lie algebra $\h$ of $H$ decomposes
$\h=(\fk\cap\h)+(Lie(\p)\cap\h)$, so that $\dim H'=\dim(K\cap H)+
\dim(Lie(\p)\cap\h)=\dim H$. Thus $H'=H$ since $H$ is connected.

2) Since $\phi$ surjects $H$,
we have $\p\cap H\simeq Lie(\p)\cap\h$ via $\phi$.
Therefore, $\p\cap H$ is a subspace of $\p$. 

3) A lattice of the vector space $\p\cap H$ becomes a uniform lattice of $H$.

4) This is an immediate consequence of (1).
\end{proof}
\begin{rem}
By the same argument, one sees easily that
Lemma~\ref{la:rectangle} still holds if the number of connected components of
$H$ is finite.
\end{rem}
The following lemma is a `tangential' analog of Lemma~\ref{la:dg}.
As in the proof of \cite[Theorem 4.7]{koba89}, one could prove
it by using the cohomological dimension of
abstract groups. However, we shall supply much elementary
approach for our Cartan motion groups below:
\begin{la}\label{la:tan:dg}
Let $G$ be a Cartan motion group.
Suppose that $L$ and $H$ are $\ct$-stable subgroups with at most finitely
many connected components. Assume $L\proper H$. Then, we have:
\begin{enumerate}
\ritem[1)]$d(L)+d(H)\leq d(G)$.
\ritem[2)]$d(L)+d(H)=d(G)$ if and only if $L\backslash G/H$ is compact.
\end{enumerate}
\end{la}
\begin{proof}
Let $\p_L:=L\cap\p$ and $\p_H:=H\cap\p$, then we have
$d(L)=\dim\p_L$ and $d(H)=\dim\p_H$.

1) Since $L\similar\p_L$ and $H\similar\p_H$, we have
$\p_L\proper\p_H$. Then it follows from
Lemma~\ref{la:cone} that we have $\p_L\cap\p_H=\z$, and
thus $\dim\p_L+\dim\p_H\leq\p$.

2) First, assuming $d(L)+d(H)=d(G)$, we shall prove
$G=LKH$.
To see this, we note $\p_H\similar\Ad(k)\p_H$ for any $k\in K$.
Therefore, we have $\p_L\proper\Ad(k)\p_H$, and then
$\p_L\cap \Ad(k)\p_H=\z$. Thus
we obtain $\p_L+\Ad(k)\p_H=\p$ in light of dimension.
Now, for any $(k,X)\in G$, we find $(e,Y)\in \p_L$ and $(e,Z)\in \p_H$
such that $X=Y+\Ad(k)Z$. This implies
$$(k,X)=(e,Y)(k,0)(e,Z)\in \p_LK\p_H\subset LKH.$$
Hence, $G=\p_L K\p_H=LKH$. We have thus proved that
$L\backslash G/H$ is compact.

Conversely, let us prove $d(L)+d(H)=d(G)$ if $L\backslash G/H$ is compact.
It follows from
Lemma~\ref{la:ssimilar} that $\p_L\backslash G/\p_H$ is compact.
Now, $\p$ is a closed subgroup containing $\p_L$ and $\p_H$.
Therefore, $\p_L\backslash \p/\p_H$ is a closed subset of
the compact set $\p_L\backslash G/\p_H$.
Since $\p$ is a vector space, $\p_L\backslash \p/\p_H$ is compact
only if
$\dim\p_L+\dim\p_H\geq \dim\p$.
Hence, Lemma has been proved.
\end{proof}

Next, we state a criterion of $\proper$ and $\similar$ for
Cartan motion groups (see Lemma~\ref{la:reductive} for the counterpart
in the reductive case).
\begin{la}[{\cite{yoshinocartan}}]\label{la:cartan}
For any subsets $L$ and $H$ in $G$, we have
\begin{enumerate}
\ritem[1)] $\properin{L}{H}{G} \iff \properin{\fa(L)}{\fa(H)}{\fa}$.
\ritem[2)] $\similarin{L}{H}{G} \iff \similarin{\fa(L)}{\fa(H)}{\fa}$.
\end{enumerate}
Here we recall \eqref{eq:df:tan:fa} for the definition of $\fa(L)$.
\end{la}
\begin{rem}
Similar results hold for an arbitrary semidirect product group
$G=K\ltimes_{\tau}V$ associated to a finite dimensional representation
$(\tau,V)$ of a compact group $K$, if we replace $\fa$ by any
subspace $V'$ such that $\tau(K)V'=V$.
\end{rem}

\begin{la}\label{la:cartancone}
Let $G$ be a Cartan motion group, and $H$ its $\ct$-stable subgroup with finitely
many connected components. Then, $\fa(H)$ is a closed cone.
\end{la}
\begin{proof}
In light of $H=(K\cap H)\ltimes(\p\cap H)$,
Lemma follows from
$\fa(H)=\fa\cap(K(\p\cap H)K)=\fa\cap\bigcup_{k\in K}\Ad(k)(\p\cap H)$.
\end{proof}

Finally, we prove Theorem~\ref{thm:tan:simple}.
\begin{proof}[Proof of Theorem~\ref{thm:tan:simple}]
It is enough to prove the equivalence \eqref{thm:carmot:criterion}(ii)
$\iff$ \eqref{thm:tan:simple}(ii).
For this, we shall set $W:=L\cap\p$ when we start from
\eqref{thm:carmot:criterion}, while we shall just set
$L:=W$ when we start from \eqref{thm:tan:simple}.
Then the following conditions are equivalent:
\begin{align*}
\properin{L}{H}{G}&\iff \properin{W}{H}{G}&&\text{by Lemma~\ref{la:rectangle}(4)}\\
&\iff\properin{\fa(W)}{\fa(H)}{\fa}&&\text{by Lemma~\ref{la:cartan}}\\
&\iff\fa(W)\cap\fa(H)=\z&&\text{by Lemmas~\ref{la:cone} and \ref{la:cartancone}.}
\end{align*}
On the other hand, we also have:
\begin{align*}
L\backslash G/H\text{ is compact}&\iff W\backslash G/H\text{ is compact}
&&\text{by Lemmas~\ref{la:ssimilar} and \ref{la:rectangle}}\\
&\iff \dim(W)+d(H)=d(G)&&\text{by Lemma~\ref{la:tan:dg}}.
\end{align*}
Thus, the conditions
\eqref{thm:carmot:criterion}(ii) and \eqref{thm:tan:simple}(ii) are
equivalent.
\end{proof}
Thus we have proved Theorem~\ref{thm:tan:simple}
by admitting Theorem~\ref{thm:carmot:criterion}, which we shall
prove in the next subsection.

\newcommand{\gam}{\Gamma}
\newcommand{\gamz}{{\Gamma_0}}
\newcommand{\gamzz}{{\Gamma_{(0)}}}
\newcommand{\kgam}{{K_\Gamma}}
\newcommand{\kgamz}{{(K_\Gamma)_0}}
\newcommand{\ggam}{{G_\Gamma}}
\newcommand{\vgam}{{V_\Gamma}}
\newcommand{\sgam}{\gamma}

\subsection{
Proof of Theorem~\ref{thm:carmot:criterion}}
\label{subsec:pf:tan:suf}
This subsection gives a proof of
Theorem~\ref{thm:carmot:criterion}.
The implication (ii)$\Rightarrow$(i) of Theorem~\ref{thm:carmot:criterion}
is an immediate consequence of
Theorem~\ref{thm:sufficient} and Lemma~\ref{la:rectangle}(3).

On the other hand, the implication 
(i)$\Rightarrow$(ii) follows from the next lemma.
\begin{la}\label{la:discrete}
Suppose that $(\tau, V)$ is a finite dimensional real representation
of a compact group $K$. We form a semidirect group $G=K\ltimes_\tau V$.
For any discrete
subgroup $\Gamma$ in $G$, there exists a subspace $W\subset V$ 
such that
$ \ssimilarin{\Gamma}{W}{G}$.
\end{la}
In fact, suppose $\Gamma$ is a cocompact discontinuous group for $G/H$.
Let $W$ be the subspace of $\p$ as in Lemma~\ref{la:discrete}, and
we set $L=W$. Then $L$ is a connected $\ct$-stable subgroup satisfying
\eqref{thm:carmot:criterion}(a) and (b) by Lemma~\ref{la:quotient_cpt}.
Hence, Theorem~\ref{thm:carmot:criterion} has been reduced to
Lemma~\ref{la:discrete}.

We need a generalization of results of Auslander \cite{auslander60},
\cite{auslander61} and Wang \cite{wang}, as follows:
\begin{la}[{\cite[Theorem 8.24]{ragtan}}]\label{la:auslander}
Let $G$ be a Lie group and $V$ a closed connected
solvable normal subgroup. Let $\pi: G\rightarrow G/V$ be the natural
quotient map. Given a closed subgroup $\gam$ of $G$, we write 
$\kgam:=\overline{\pi(\gam)}$ for the closure of
$\pi(\gam)$. If the identity component $\gamz$ of $\gam$
is solvable, so is  $\kgamz$.
\end{la}
Then, Lemma~\ref{la:discrete} is a consequence of the following Lemma.
\begin{la}\label{la:abelian}
Suppose $G=K\ltimes_\tau V$ is as in Lemma~\ref{la:discrete}. We
write  $\pi: G\rightarrow
G/V\simeq K$ for the natural quotient map.
Suppose $\gam$ is a closed subgroup of $G$.
Let $\kgamz$ be the subgroup of $K$ as in Lemma~\ref{la:auslander}.
If $\kgamz$ is abelian, then
there exists a subspace $W\subset V$ such that $\ssimilarin{\gam}{W}{G}$.
\end{la}
Let us give a proof of Lemma~\ref{la:discrete} by admitting
Lemma~\ref{la:abelian}.
\begin{proof}[Proof of Lemma~\ref{la:discrete}]
Apply Lemma~\ref{la:auslander} to the case:
$\Gamma$ is a discrete subgroup of a Cartan motion group $G=K\ltimes\p$.
Since $\gamz=\{e\}$ is solvable, so is $\kgamz$.
Then $\kgamz$ must be abelian because
$\kgamz$ is  compact group.
Now, Lemma~\ref{la:discrete} follows from
Lemma~\ref{la:abelian}.
\end{proof}

The rest of this subsection is devoted to the proof of
Lemma~\ref{la:abelian}. In the course of the proof,
we need to treat also the case where $\Gamma$ is not discrete.

To see Lemma~\ref{la:abelian}, we start with
the following multiplicative version of the Jordan decomposition.
Let $K$ be a closed subgroup of $O(n)$ and $V=\R^n$.
We form a semidirect product group $G:=K\ltimes V (\subset GL(n+1,\R))$.
\newcommand{\tilk}{\widetilde{K}}
We define a subset $\tilk$ of $G$ consisting of elliptic elements as follows:
\begin{align*}
\tilk&:=\bigcup_{g\in V} gKg^{-1}\\
&=\set{s\in G}{\set{s^k}{k\in\Z} \text{ is relatively compact}}.
\end{align*}
\begin{la}\label{la:decomp}
\begin{enumerate}
\ritem[1)]
For any $g\in G$, there exist uniquely $s\in \tilk$ and $w\in V$ satisfying
$g=sw=ws$.
\ritem[2)]
Both $s$ and $w$ can be written as polynomials of $g$ in $M(n+1,\R)$.
\end{enumerate}
\end{la}
\begin{proof}
1)
We fix $g=(k,v)\in G$ once and for all.
Let $V_1\subset V$ be the eigenspace of $k$
with eigenvalue 1, and $V_2$ the orthogonal complement of $V_1$ in $V$.
According to the direct sum decomposition $V=V_1\oplus V_2$,
we write $v=v_1+v_2$. We set
\begin{equation*}
s:=(k,v_2),\qquad w:=(e,v_1).
\end{equation*}
First, we claim $s\in\tilk$. Since the linear map $(\id-k)$ preserves
$V_2$ and acts as an invertible transformation,
we can take $u\in V_2$ such that $v_2=(\id-k)u$. Then, we obtain
$s=(k,v_2)=(e,u)(k,0)(e,-u)\in\tilk$.

Next, we verify $g=sw=ws$. In view of $kv_1=v_1$, we have
\begin{align*}
ws&=(e,v_1)(k,v_2)=(k,v_1+v_2)=g,\\
sw&=(k,v_2)(e,v_1)=(k,v_2+v_1)=g.
\end{align*}
Finally, the uniqueness follows immediately
from $\tilk\cap V=\{(e,0)\}$. Hence (1) is proved.

2) Admitting the existence of a polynomial $p(x)=\sum_{n=0}^Nc_nx^n$ satisfying
\begin{equation}\label{eq:polynomial}
p(k)=0,\quad p(1)=0,\quad \text{and}\quad p'(1)=1,
\end{equation}
we first finish the proof of (2). For $g=(k,v)=(k,v_1+v_2)=(k,(\id-k)u+v_1)$,
we have inductively,
\begin{equation*}
g^n=(k^n,(\id-k^n)u+nv_1).
\end{equation*}
Therefore, we have
\begin{align*}
p(g)&=\sum_{n=0}^N(c_nk^n,(c_n\id-c_nk^n)u+nc_nv_1) \\
&=(p(k), (p(1)\id-p(k))u+p'(1)v_1)\\
&=(0,v_1).
\end{align*}
Thus, $s$ and $w$ are polynomials of $g$ given by:
\begin{equation*}
s=g-p(g),\qquad w=p(g)+\id.
\end{equation*}

Finally, let us construct a polynomial $p(x)$ satisfying \eqref{eq:polynomial}.
Suppose $q(x)$ is the minimal polynomial of the matrix
\begin{equation*}
\begin{pmatrix}
k& \\
& 1
\end{pmatrix} \in O(n+1).
\end{equation*}
Clearly,  $q(k)=0$ and $q(1)=0$.
Since any element of $O(n+1)$ is semisimple, $q(x)$ has no multiple root.
Hence $q'(1)\neq 0$. Then $p(x):=\frac{q(x)}{q'(1)}$
is the desired one.
We have thus proved the lemma.
\end{proof}

We are ready to prove Lemma~\ref{la:abelian}. Without loss
of generality, we may and do assume that $K$ is a closed subgroup
of $O(V)$. Then, we can regard $G=K\ltimes V$ as a subgroup of
$GL(n+1,\R)$, where $n=\dim V$.
\begin{proof}[Proof of Lemma~\ref{la:abelian}] 
We shall divide the proof into three steps.

\noindent
{\bf Step 1}: $\Gamma$ is abelian.

Assume that $\gam$ itself is abelian. The multiplicative Jordan
decomposition $\sgam=sw$ (see Lemma~\ref{la:decomp}) for $\sgam\in \gam$
defines the projections:
\begin{alignat*}{2}
p_1&: \gam \longrightarrow G, \qquad&  \sgam & \longmapsto s,\\
p_2&: \gam \longrightarrow V, & \sgam & \longmapsto w.
\end{alignat*}
Under the assumption that $\gam$ is abelian, both $p_1$ and $p_2$
are group homomorphisms. In fact, if $\sgam_i=s_iw_i\ (i=1,2)$, then
$\sgam_1,\sgam_2,s_1,s_2, w_1$ and $w_2$ mutually commute because $s_i$ and
$w_i$ are polynomials of $\sgam_i$.
We set
\begin{equation*}
S:=\overline{p_1(\gam)},\qquad W':=p_2(\gam).
\end{equation*}
Then we claim $\ssimilarin{\gam}{W'}{G}$. To see this, it is enough to
show $S$ is compact because
\begin{equation}
\gam\subset W'S,\qquad W'\subset \gam S. \label{eq:lws}
\end{equation}
Since any element of $S$ commutes with each other,
we can find a matrix $p\in GL(n+1,\R)$ such that $psp^{-1}$ is simultaneously 
a diagonal matrix
for all $s\in S$.
As all the eigenvalues of $s$ have modulus one, we have
$pSp^{-1}\subset \set{\diag(\lambda_1,\dots,\lambda_{n+1})}{|\lambda_i|=1}$,
and thus $S$ is compact.
We have thus proved $\ssimilarin{\gam}{W'}{G}$.

On the other hand, let $W:=\rSpan W'$, then we have
$\ssimilarin{W'}{W}{G}$ by Lemma~\ref{la:quotient_cpt}.
Thus we have proved $\ssimilarin{\gam}{W}{G}$.

\noindent
{\bf Step 2}: $\pi(\gam)$ is abelian.

\noindent
Assume that $\pi(\gam)$ is abelian.
We recall $\kgam=\overline{\pi(\gam)}$, and we set
$\ggam:=\kgam\ltimes V$, and $\vgam:=\rSpan(\gam\cap V)$.
Our idea here is to remove the
`non-commutative part' of $\Gamma$ by taking a quotient map
\begin{equation}
\varpi: \ggam\rightarrow \ggam/\vgam.\label{eq:wprj}
\end{equation}
This enables us to apply Step 1 to $\ggam/\vgam$.

First, let us verify that $\vgam$ is a normal subgroup in 
$\ggam(=\kgam\ltimes V)$.
It is enough to show that $\vgam$ is $\pi(\gam)$-invariant.
For any $v\in V$ and any $\sgam\in\Gamma$, we have
$$\pi(\sgam)v=\sgam v\sgam^{-1}.$$
Here, $\pi(\sgam)v$ in the left-hand side is defined
by the action of $K$ on $V$, while $\sgam v\sgam^{-1}$ in 
the right-hand side is
the group multiplication in $G$.
Thus, the action of $\pi(\gam)$
on $V$ preserves $\gam\cap V$, and then preserves
also $\vgam=\rSpan(\gam\cap V)$.

\newcommand{\ovp}{\overline{\varpi}}
Next, let us verify $\varpi(\gam)$ is abelian. 
The natural quotient map $\ggam/\vgam\rightarrow\ggam/V$
induces a surjective homomorphism $\mu:\varpi(\gam) \rightarrow \pi(\gam)$.
Then $\mu$ is also injective
because $V\cap\gam\subset\vgam$. Hence $\varpi(\gam)$
is abelian.

Applying Step 1 to $\ggam/\vgam (\simeq \kgam\ltimes
V/\vgam)$ and its abelian subgroup $\varpi(\gam)$,
we find a subspace $W'\subset V/\vgam$ such that
$$ \ssimilarin{\varpi(\gam)}{W'}{\ggam/\vgam}.$$
It follows from Lemma~\ref{la:bieber_normal} that
$\ssimilarin{\varpi^{-1}(\varpi(\gam))}{\varpi^{-1}(W')}{\ggam}$.
We set $W:=\varpi^{-1}(W')$. Since $\ggam\subset G$, we have
$$ \ssimilarin{\varpi^{-1}(\varpi(\gam))}{W}{G}.$$
To see $\ssimilarin{\gam}{W}{G}$, it is enough to show
$\ssimilarin{\gam}{\varpi^{-1}(\varpi(\gam))}{G}$, which
will follow if we prove $\varpi^{-1}(\varpi(\gam))/\gam$ is compact
by Lemma~\ref{la:quotient_cpt}. 
This then follows from
\begin{equation}
\varpi^{-1}(\varpi(\gam))/\gam=\gam\vgam/\gam = \vgam/\vgam\cap\gam
=\vgam/V\cap\gam \label{eq:gam}
\end{equation}
because $\vgam=\rSpan(V\cap\gam)$ contains $V\cap\gam$
as a cocompact subgroup. Here, the first two equalities
in \eqref{eq:gam} are clear, and the last one follows from
$\vgam\cap\gam=V\cap\gam$, as is seen by taking
the intersection of $\gam$ with
$V\cap\gam\subset\vgam\subset V$.
Hence, we have proved $\ssimilarin{\gam}{W}{G}$.

\smallskip
\noindent
{\bf Step 3}: $\kgamz$ is abelian.

Finally, suppose $\gam$ satisfies the assumption of Lemma,
namely, $\kgamz$ is abelian.
We set a subgroup
$\gamzz:=\pi^{-1}(\kgamz)\cap \gam$ of $\gam$.
Since $\pi(\gamzz)\subset \kgamz$,
$\pi(\gamzz)$ is abelian.
Applying Step 2 to $(G,\gamzz)$, there exists a subspace
$W\subset V$ such that
$$ \ssimilarin{\gamzz}{W}{G}.$$
To see $\ssimilarin{\gam}{W}{G}$,
it is enough to prove $\gam/\gamzz$ is compact.

Since $K$ is compact, $\kgam/\kgamz$ is a finite group. Then, the
composition
$$\overline{\pi}: \gam\rightarrow \kgam \rightarrow \kgam/\kgamz$$
is a surjective homomorphism, because $\pi(\gam)$ is dense in $\kgam$.
On the other hand, we have $\Ker(\overline{\pi})=\gamzz$. Thus,
$\overline{\pi}$ induces the isomorphism:
$$ \gam/\gamzz\simeq \kgam/\kgamz.$$
Thus, $\gam/\gamzz$ is finite.
Hence, $\ssimilarin{\gam}{W}{G}$.
This completes the proof of Lemma~\ref{la:abelian}.
\end{proof}
Thus the proof of Theorem~\ref{thm:carmot:criterion} is completed.
\begin{rem}
We did not assume $\gam$ is discrete in Lemma~\ref{la:abelian}.
But, if $\gam$ is a discrete subgroup, then $\kgamz$ is
abelian by Lemma~\ref{la:auslander}, and therefore, the semidirect
product group $(\ggam)_0:=\kgamz\ltimes V$ is solvable.
Thus, $\gam$ is virtually solvable, as
$\gamzz$ is contained in
$(\ggam)_0$ and
 $\gamzz$ is of finite index in $\gam$.
Then, it follows from Fried and Goldman \cite[Theorem 1.4]{friedgoldman}
that
there exists a closed subgroup $L$
of $(\ggam)_0$ such that $\gamzz$ is cocompact in $L$ and $L$
is of finitely many components.
\end{rem}
\subsection{Solution to tangential space form problem}
\label{subsec:tangential_space}
So far, we have discussed a general theory concerning Problem A for
tangential symmetric space associated to reductive symmetric spaces.
Building on it, we shall study specific case, namely,
the tangential space form problem.

This subsection gives an outline of the proof of Theorem~\ref{thm:tangential}.
The proof will be completed in Section~\ref{subsec:tangential_bilinear}.

Theorem~\ref{thm:carmot:criterion} has reduced the existence problem of
a cocompact discontinuous group $\Gamma$ for the tangential symmetric
space $G/H$ to that
of a certain connected group $L$. However, we shall find that
the latter problem involves surprisingly rich ingredients
even for a fairly simple pair $(G,H)$. Theorem~\ref{thm:tangential}
 treats the
case where $G/H$ (previously denoted by $G_\theta/H_\theta$) 
is the tangential symmetric space of the symmetric
space $O(p+1,q)/O(p,q)$, and classifies the pair $(p,q)$ of integers
for which such $L$ exists.  Our argument
relies on Adams' theorem on vector fields on spheres (\cite{adams}).
Our strategy of proof and related classical results are listed
in the following proposition.
\begin{prop}\label{eq:hurwitz:radon}
The following seven conditions on the pair $(p,q)$ of natural numbers
are equivalent:
\begin{enumerate}
\ritem[i)] The tangential symmetric space of $O(p+1,q)/O(p,q)$
admits a compact Clifford-Klein form.
\ritem[ii)] There exists a bilinear map
$f:\R^{p+1}\times\R^q\rightarrow \R^q$ such that
$$ f(v,w)=0  \text{ only if\/ $v=0$ or $w=0$}.$$
\ritem[iii)] There exists a bilinear map
$f:\R^{p+1}\times\R^q\rightarrow\R^q$ such that
$$ \|v\|\cdot\|w\|=\|f(v,w)\|\qquad\text{for any $v\in\R^{p+1}$ and
$w\in\R^q$}.$$

\ritem[iv)] There exists an injective linear map
$F:\R^p\rightarrow M(q,\R)$ such that\/
$\Image F\subset \fo(q)\cap\R\cdot O(q)$.
\ritem[v)] There is an identity
$(x_1^2+\dots+x_{p+1}^2)\cdot (y_1^2+\dots+y_q^2)=z_1^2+\dots+z_q^2$
such that $z_1,\dots,z_q$ are bilinear functions of
$x_1,\dots, x_{p+1}$ and $y_1,\dots,y_q$.
\ritem[vi)] There exist $p$ vector fields on
the sphere $S^{q-1}$ which are linearly independent at
every point.
\ritem[vii)] $p<\rho(q)$\quad $($see \eqref{eq:df:rho}$)$.
\end{enumerate}
\end{prop}
Theorem~\ref{thm:tangential} corresponds to the equivalence of (i) and (vii)
in Proposition.
\begin{rem}
We point out that there is a mysterious duality between $p$ and $q$
in (i) and (vi). That is, in (i), the tangential symmetric
space of $O(p+1,q)/O(p,q)$ is diffeomorphic to
\begin{equation}
\text{the trivial $\R^q$-bundle over $S^p$}
\end{equation}
as $K$-equivariant fiber bundles, whereas (vi) defines
\begin{equation}
\text{the trivial $\R^p$-bundle over $S^{q-1}$}
\end{equation}
as a subbundle of the tangent bundle $T(S^{q-1})$ of $S^{q-1}$.
\end{rem}

 The proof of Proposition will be organized as follows:
$$
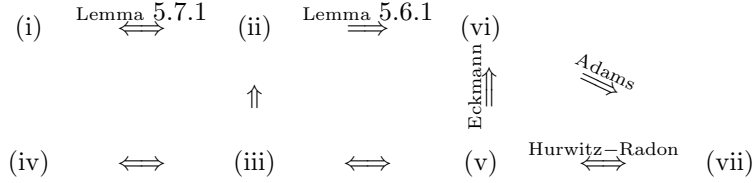
\begin{array}{ccccccc}
{\rm (i)}&\stackrel{\rm Lemma~\ref{la:tangential_biliner}}{\iff}&{\rm (ii)}&\stackrel{\rm Lemma~\ref{la:vector_fields_on_sphere}}{\Longrightarrow}&{\rm (vi)}\\
&&\raise 10pt \hbox{$\Uparrow$}&&
\hbox{\rotatebox{90}{$\stackrel{\text{Eckmann}}{\Longrightarrow}$}} &
\raise 20pt \hbox{\rotatebox{-25}{$\stackrel{\text{Adams}}{\Longrightarrow}$}}
\\
{\rm (iv)}&\iff&{\rm (iii)}&\iff&{\rm (v)}&\stackrel{\rm\hss Hurwitz-Radon\hss}{\iff}&{\rm (vii)}
\end{array}
$$
\noindent
(v)$\iff$(vii): This equivalence was established by
Hurwitz \cite{hurwitz} and Radon \cite{radon}.
See also Eckmann \cite{eckmann42}, \cite{eckmann94}.
\par\smallskip\noindent
(v) $\Longrightarrow$ (vi): This statement was proved in Eckmann
\cite{eckmann42}.
\par\smallskip\noindent
(vi) $\Longrightarrow$ (vii): This statement was proved in Adams \cite{adams}.
\par\smallskip\noindent
(iii) $\Longrightarrow$ (ii): Obvious.
\par\smallskip\noindent
(iii)$\iff$(v): Obvious.
\par\smallskip\noindent
(iii)$\iff$(iv): This equivalence is an easy observation of
linear algebra by the correspondence
$f((x,c),\cdot\ )=F(x)+c I_q$ for $(x,c)\in\R^p\oplus \R$.
\par\smallskip\noindent
(ii) $\Longrightarrow$ (vi): This proof is in the same line of the
argument
(v) $\Longrightarrow$ (vi), but we shall give a proof
in Lemma~\ref{la:vector_fields_on_sphere} for the sake of completeness.

Therefore, the proof of Proposition will be completed if we show
the equivalence (i)$\iff$(ii). This will be carried out in
Lemma~\ref{la:tangential_biliner}.

We shall also give an alternative proof of the implication
(vii) $\Longrightarrow$ (iii) (or (vii)~$\Longrightarrow$~(v)) in
Lemma~\ref{la:id}
as an application of the theory of Clifford algebras associated to
indefinite quadratic forms.

\subsection{Results of Hurwitz-Radon-Eckmann}
\label{subsec:vector_field_on_sphere}
This subsection provides a simple and self-contained 
proof of the implication (ii)$\Rightarrow$(vi)
(see Lemma~\ref{la:vector_fields_on_sphere}) and (vii) $\Rightarrow$ (iii)
(see Lemma~\ref{la:id}).

\begin{la}\label{la:vector_fields_on_sphere}
Let $(p,q)$ be a pair of natural numbers, and assume that there
exists a bilinear map $f:\R^{p+1}\times\R^q\rightarrow \R^q$
satisfying \eqref{eq:division}.
Then, there exists $p$ vector fields on the sphere $S^{q-1}$ which
are linearly independent at every point.
\end{la}
\begin{proof}
Let $\{e_0,\dots,e_p \}$ be the standard basis of $\R^{p+1}$.
Since $f(\,\cdot\ , w):\R^{p+1}\rightarrow\R^q$ is injective
for any $w\in S^{q-1}$, $f(e_0, w), \dots, f(e_p,w)\in \R^q$
are linearly independent. In light that
$g:=f(e_0,\,\cdot\ ): \R^q\rightarrow\R^q$ is invertible,
$t_i(w):=g^{-1}f(e_i, w)$ ($0\leq i\leq p$) are well defined
and linearly independent. We note $t_0(w)=w$.
Now we set
$$Z_i(w):=t_i(w)-\langle t_i(w), w\rangle w\qquad (i=1,\dots,p).$$
Then, $Z_1,\dots,Z_p$ define $p$ vector fields on $S^{q-1}$
which are linearly independent at every point.
\end{proof}

\begin{la}\label{la:id} For any positive integer $q$, there
exists a linear map
$$ f:\R^{\rho(q)}\rightarrow M(q,\R)$$
such that $\tr{f(v)}f(v)=\|v\|^2 I_q$ for any $v\in \R^{\rho(q)}$.
\end{la}
This Lemma is due to Hurwitz \cite{hurwitz} and Radon \cite{radon}. Later,
Eckmann \cite{eckmann42} also gave a simple proof (see also
\cite{eckmann94}).
Here, we give yet another simple proof based on the
Clifford algebra associated to the indefinite quadratic form.
\begin{proof}
It is sufficient to prove Lemma in the case where $q$ is of the form
$q=2^k$ because the Hurwitz-Radon number $\rho(q)$ remains the same if
we multiply $q$ by any odd number. Henceforth, suppose
$q=2^{4\alpha+\beta}$ ($\beta=0$, $1$, $2$, $3$).
According to the value $\beta$, we set
\par\medskip\noindent
\halign{ \hfil\quad Case #&)\quad $\beta$#&$=$#\quad&$(r,s):=(#)$.\cr
1&& $0$, $1$ or $2$. & \beta,8\alpha+\beta\cr
2&& $3$. & 0,8\alpha+6\cr
}

\medskip
In either case, $r-s-1\equiv \pm 1\mod 8$ and $q=2^{(r+s)/2}$. Therefore,
it follows from Proposition~\ref{prop:clifiso} that there exists an algebra isomorphism
$\iota: C(r,s)\simarrow M(q,\R)$,
and we have from Lemma~\ref{la:st_compat}
\begin{equation}
\tr{\iota(x)}=\iota(\st{x}).\label{eq:iota_compat}
\end{equation}
We assert:
\par\medskip\noindent
Claim: For $q=2^{4\alpha+\beta}$, we set $p:=\rho(q)=8\alpha+2^\beta$.
Then, there exist $x_1,\dots,x_{p-1}\in C(r,s)$ with the following
three conditions:
\begin{enumerate}
\ritem[1)]$x_i^2=-1$\qquad ($1\leq i\leq p-1$),
\ritem[2)]$x_i x_j + x_j x_i=0$\qquad  ($i\neq j$, $1\leq i,j\leq p-1$),
\ritem[3)]$\st{x_i}=-x_i$\qquad ($1\leq i\leq  p-1$).
\end{enumerate}
First, let us complete the proof of Lemma by admitting Claim.
The substitution of (3) into (1) and (2) leads to
\begin{equation}
\st{x_i}x_i=1,\qquad \st{x_i}x_j + \st{x_j}x_i=0\qquad \text{($i\neq j$),}
\label{eq:orthogonal}
\end{equation}
for $1\leq i,j\leq p-1$. Clearly the relation \eqref{eq:orthogonal}
still holds for $0\leq i,j\leq p-1$
if we set $x_0:=1$.
We now define a linear map $f$ by
$$ f: \R^p\mapsto M(q,\R),\quad e_i \mapsto \iota(x_i)
\qquad (0\leq i\leq p-1).$$
Then, if $v\in\R^p$ is written as
$\sum_{i=0}^{p-1} c_ie_i$, we have:
$$\tr{f(v)}\cdot f(v)=\left(\sum_{i=0}^{p-1}c_i\tr{\iota(x_i)}\right)
\left(\sum_{i=0}^{p-1}c_i \iota(x_i)\right)=\sum_{i=0}^{p-1}c_i^2 I_q
=\|v\|^2 I_q.$$
Here, the second equality follows from \eqref{eq:iota_compat}
and \eqref{eq:orthogonal}.
Therefore, Lemma~\ref{la:id} has been proved by assuming Claim.

Finally, let us prove Claim.
\par\smallskip\noindent
Case 1) $\beta=0$ or $1$.\qquad $x_i:=v_i^-\quad (i=1,\dots,p-1)$.
\par\smallskip\noindent
Case 2) $\beta=2$ or $3$.\qquad $x_i:=v_i^-\quad (i=1,\dots,p-2)$;
\quad $x_{p-1}:=J_-=v_1^-\dots v_{p-1}^-$.
Then, an easy computation shows the conditions
(1), (2) and (3) in Claim.

Hence, the proof of Lemma completed.
\end{proof}

\subsection{Proof of Theorem~\ref{thm:tangential}}\label{subsec:tangential_bilinear}
This subsection completes the proof of Theorem~\ref{thm:tangential}
by showing the equivalence of the conditions
(i) and (ii) in Proposition~\ref{eq:hurwitz:radon}, namely, the
following:
\begin{la}\label{la:tangential_biliner}
Let $(p,q)$ be a pair of natural numbers.
Then, the following two conditions are equivalent:
\begin{enumerate}
\ritem[i)] The tangential symmetric space of $O(p+1,q)/O(p,q)$
admits a compact Clifford-Klein form.
\ritem[ii)] There exists a bilinear map
$f:\R^{p+1}\times\R^q\rightarrow \R^q$ such that
$$ f(v,w)=0  \text{ only if\/ $v=0$ or $w=0$}.$$
\end{enumerate}
\end{la}
Throughout this subsection, suppose $G/H$ is
the tangential symmetric space of $O(p+1,q)/O(p,q)$.
Here we realize  $O(p+1,q)$ in $GL(p+q+1,\R)$ in the standard
way (see \eqref{eq:df:opq}) and
$$
O(p,q)=\set{g\in O(p+1,q)}{ g e_0=e_0},
$$
where $\{e_0,\dots,e_p\}$ denotes
the standard basis of the Euclidean space $\R^{p+q+1}$.
We take $\theta$ to be the standard Cartan involution of $O(p+1,q)$, given by
$\theta g=\tr{g^{-1}}$ for $g\in O(p+1,q)$. Then, the maximal compact subgroup
$K = O(p+1)\times O(q)$ acts on
\begin{equation}
\p = \left\{
\begin{pmatrix} & \tr{X}\\X&\end{pmatrix} \mid X \in M(q,p+1;\R) \right\}
\simeq M(q,p+1;\R),\label{eq:df:p}
\end{equation}
by $\p\rightarrow\p, X\mapsto k_2Xk_1^{-1}\quad (k_1\in O(p+1),\ 
k_2\in O(q))$.

It follows from Theorem~\ref{thm:tan:simple} that $G/H$ admits a compact
Clifford-Klein form if and only if the following condition (i)$'$
holds. Thus, to show the equivalence (i) $\iff$ (ii),
it is sufficient to prove the equivalence
(i)$'$ $\iff$ (ii)
\begin{enumerate}
\ritem[i)$'$]
There exists a $q$-dimensional subspace $W\subset\p$ satisfying
\begin{equation}
\fa(W)\cap\fa(H)=\{0\}.\label{eq:tangential_space_proper}
\end{equation}
\ritem[ii)]
There exists a bilinear map $f:\R^{p+1}\times \R^q\rightarrow
\R^q$ such that
\begin{equation}
f(v,w)=0 \text{ only if $v=0$ or $w=0$.}\label{eq:division}
\end{equation}
\end{enumerate}
\begin{proof}
Given a $q$-dimensional subspace $W$ in $\p$, we 
fix a linear isomorphism
\begin{equation}
\tilf: \R^q\simarrow W \subset M(q,p+1;\R)\simeq \Hom(\R^{p+1},\R^q),
 \label{eq:ftild}
\end{equation}
and then define a bilinear map
$$f:\R^{p+1}\times \R^q\rightarrow\R^q$$
by using the canonical identification
\begin{equation}
\tilf \leftrightarrow f,\quad \Hom(\R^q,\Hom(\R^{p+1},\R^q))\simeq
\Hom(\R^{p+1}\otimes\R^q,\R^q).\label{eq:tilff}
\end{equation}
Conversely, given a bilinear map $f$ satisfying
\eqref{eq:division}, we define a
subspace $W_f$ of $M(q,p+1;\R)$ by
$$W_f:=\Image\tilf \qquad \text{(see \ref{eq:tilff}).}$$
Then, $\dim W_f=q$ because $\tilf$ is injective.
Now, the equivalence (i) $\iff$ (ii) follows from the next claim.
\par\medskip\noindent
Claim: The condition \eqref{eq:tangential_space_proper}
and \eqref{eq:division} are equivalent
via the correspondence $f\leftrightarrow W=W_f$.

To see Claim, let us find $\fa(H)$. Via the identification \eqref{eq:df:p},
we have
$$ H\cap\p=\set{X\in\p}{X e_0=0}.$$
Then $\fa(H)$ (see \eqref{eq:df:tan:fa}) is computed as follows:
\begin{align}
\fa(H) &= KH K\cap \fa \notag
\\
&=\set{X\in\fa}{(e,X)\in KHK} \notag\\
&=\set{X\in\fa}{k_1Xk_2^{-1}\in H\cap\p\quad \text{for some $(k_1,k_2)
\in K$}} \notag\\
&=\set{X\in\fa}{k_1Xk_2^{-1}e_0=0 \quad \text{for some $(k_1,k_2)\in K$}}
\notag \\
&=\set{X\in\fa}{Xv=0 \quad \text{for some non-zero $v\in\R^{p+1}$}}.
\label{eq:aht}
\end{align}
Suppose that (ii) fails. We shall show (i) fails. Take $v\in\R^{p+1}
\setminus\{0\}$ and $w\in \R^q\setminus\{0\}$ such that
$f(v,w)=0$. We set $X:=\tilf(w)\in W$. Then $X\neq 0$ and
$X v=f(v,w)=0$. In light of $W\subset K\fa(W)K$,
we can find $k_1,k_2\in K$ and a non-zero $X'\in\fa(W)$
such that $X=k_1 X' k_2$. Then $v':=k_2 v\neq 0$ and
$X' v'=k_1^{-1}Xv=0$. Hence $X'\in\fa(H)$ by \eqref{eq:aht}.
Therefore we have $X'\in \fa(W)\cap\fa(H)$. Hence (i) fails.
Likewise, we see if (i) fails then (ii) fails by using
$\fa(W)\subset KWK$.

 \end{proof}

Now, the proof of Theorem~\ref{thm:tangential} has been completed.

\medskip\noindent
{\sc Research Institute for Mathematical Sciences, Kyoto University, Sakyo-ku
Kyoto, 606-8502, Japan}\\
{\it E-mail address: toshi@kurims.kyoto-u.ac.jp, yoshino@kurims.kyoto-u.ac.jp}

\end{document}